\documentclass[12 pt]{article}
\usepackage{amsmath,amsthm,amscd,amssymb}
 \pdfoutput=1
\usepackage[margin=.75 in]{geometry}
\usepackage{subfig}
\usepackage{graphicx}
\usepackage[bookmarks=false]{hyperref}
\usepackage[labelfont=bf]{caption}
\allowdisplaybreaks

\theoremstyle{plain}

\theoremstyle{definition}

\theoremstyle{remark}

\newcommand{\mf}{\mathfrak}
\newcommand{\mc}{\mathcal}
\newcommand{\mb}{\mathbb}
\newcommand{\mbf}{\mathbf}

\newcommand{\Z}{\mathbb{Z}}

\newcommand{\R}{\mathbb{R}}
\newcommand{\C}{\mathbb{C}}
\newcommand{\ba}{\begin{aligned}}
\newcommand{\ea}{\end{aligned}}
\newcommand{\bt}{\begin{thm}}
\newcommand{\et}{\end{thm}}
\newcommand{\bc}{\begin{corollary}}
\newcommand{\ec}{\end{corollary}}
\newcommand{\bl}{\begin{lemma}}
\newcommand{\el}{\end{lemma}}
\newcommand{\bpf}{\begin{proof}}
\newcommand{\epf}{\end{proof}}
\newcommand{\bpb}{\begin{problem}}
\newcommand{\epb}{\end{problem}}
\newcommand{\bd}{\begin{definition}}
\newcommand{\ed}{\end{definition}}
\newcommand{\bn}{\begin{note}}
\newcommand{\en}{\end{note}}
\newcommand{\bp}{\begin{proposition}}
\newcommand{\ep}{\end{proposition}}
\newcommand{\be}{\begin{example}}
\newcommand{\ee}{\end{example}}
\newcommand{\bex}{\begin{exercise}}
\newcommand{\eex}{\end{exercise}}
\theoremstyle{plain}
\newtheorem{thm}{Theorem}[section]

\newtheorem{lemma}[thm]{Lemma}
\newtheorem{corollary}[thm]{Corollary}
\newtheorem{proposition}[thm]{Proposition}
\newtheorem{exercise}[thm]{Exercise}
\newtheorem{problem}[thm]{Problem}

\theoremstyle{definition}
\newtheorem{definition}[thm]{Definition}
\newtheorem{remark}[thm]{Remark}
\newtheorem{example}[thm]{Example}

\newtheorem*{convention}{Important Notational Convention}
\newtheorem*{acknowledgements}{Acknowledgements}

\newtheorem {note}{Note}[section]
\swapnumbers
\theoremstyle{remark}

\theoremstyle{plain}

\title{ A path property of Dyson gaps, \\ Plancherel measures for $Sp(\infty)$,\\ and random surface growth}
\author{Mark Cerenzia %
\thanks{Electronic address: \texttt{cerenzia@princeton.edu}} \\
{\small \emph{Princeton University}} \vspace{-10ex}}
\date{}
\begin{document}
\maketitle

\begin{abstract}
We pursue applications for symplectic Plancherel growth based on a repulsion phenomenon arising in its diffusion limit and on intermediate representation theory underlying its correlation structure. Under diffusive scaling, the dynamics converge to interlaced reflecting Brownian motions with a wall that achieve Dyson non-colliding dynamics. We exhibit non-degeneracy of constraint in this coupled system by deriving a path property that quantifies repulsion between particles coinciding in the limit. We then identify consistent series of Plancherel measures for $Sp(\infty)$ that reflect the odd symplectic groups, despite their non-semisimplicity. As an application, we compute the correlation kernel of the growth model and investigate its local asymptotics: the incomplete beta kernel emerges in the bulk limit, and new variants of the Jacobi and Pearcey kernels arise as edge limits. In particular, we provide further evidence for the universality of the $1/4$-growth exponent and Pearcey point process in the class of anisotropic KPZ with a wall. 
\end{abstract}

\tableofcontents

\section{Introduction} \label{intro}
\begin{figure}[h!] 
\centering
\includegraphics[width=0.24\textwidth]{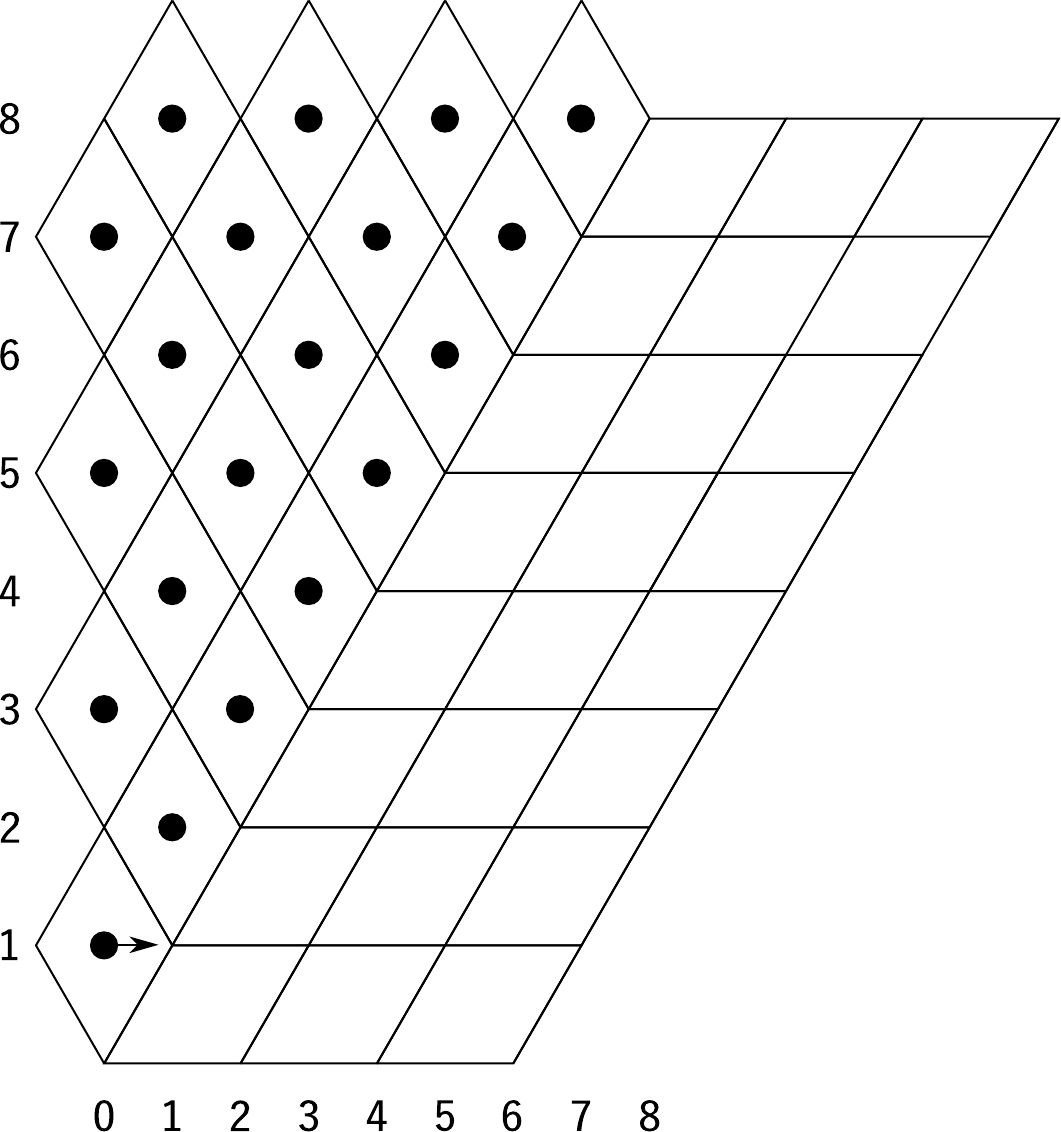} 
\includegraphics[width=0.24\textwidth]{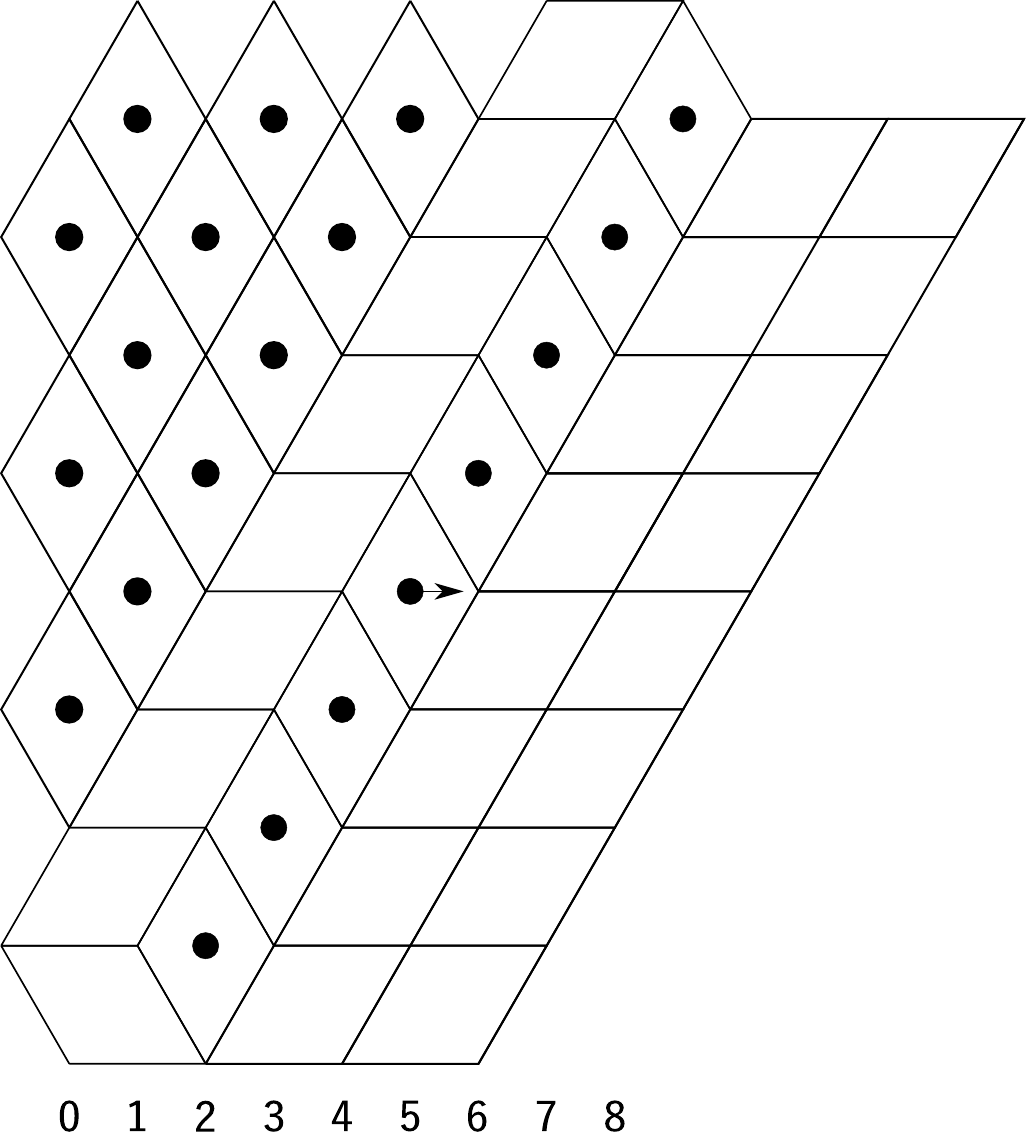} 
\includegraphics[width=0.24\textwidth]{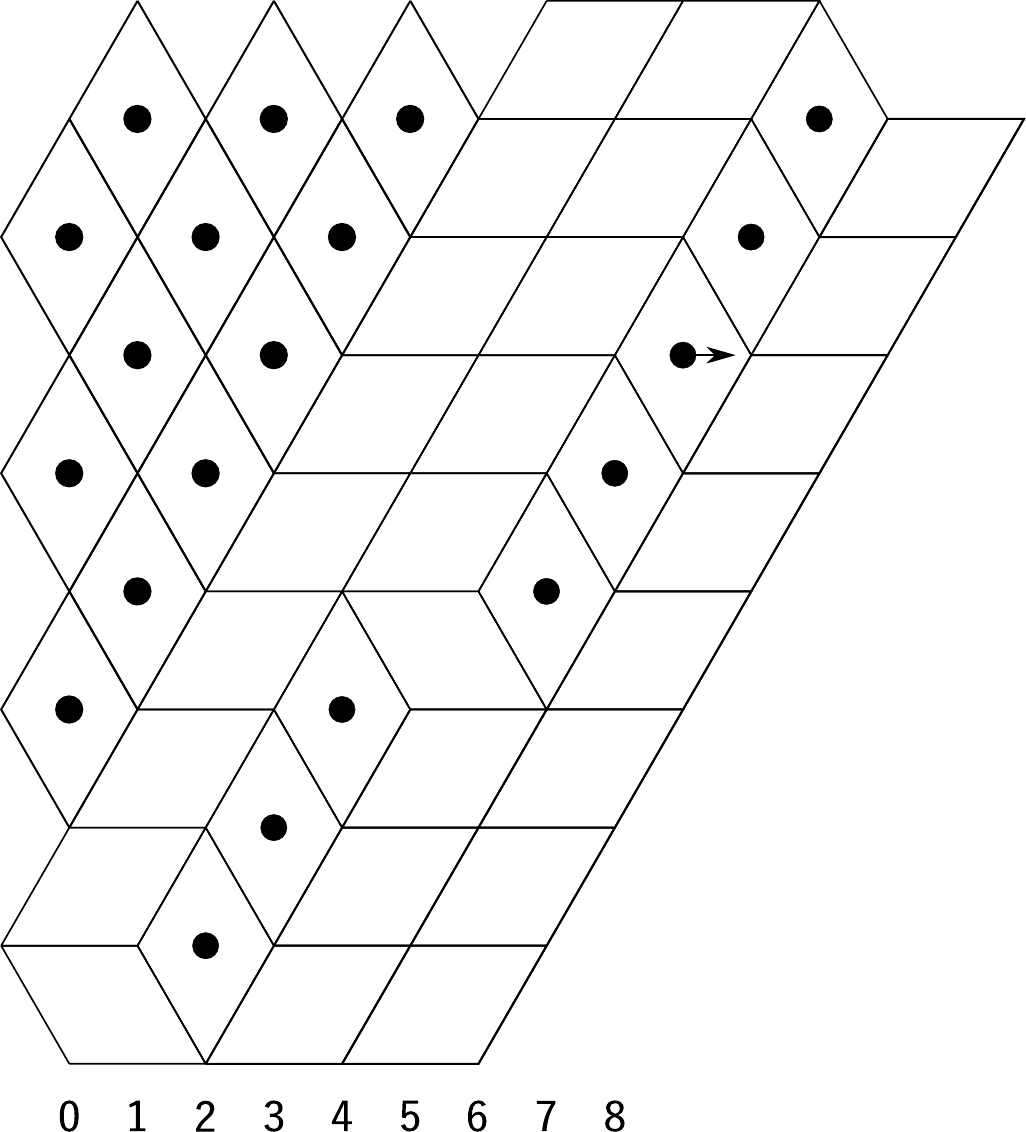}
\includegraphics[width=0.24\textwidth]{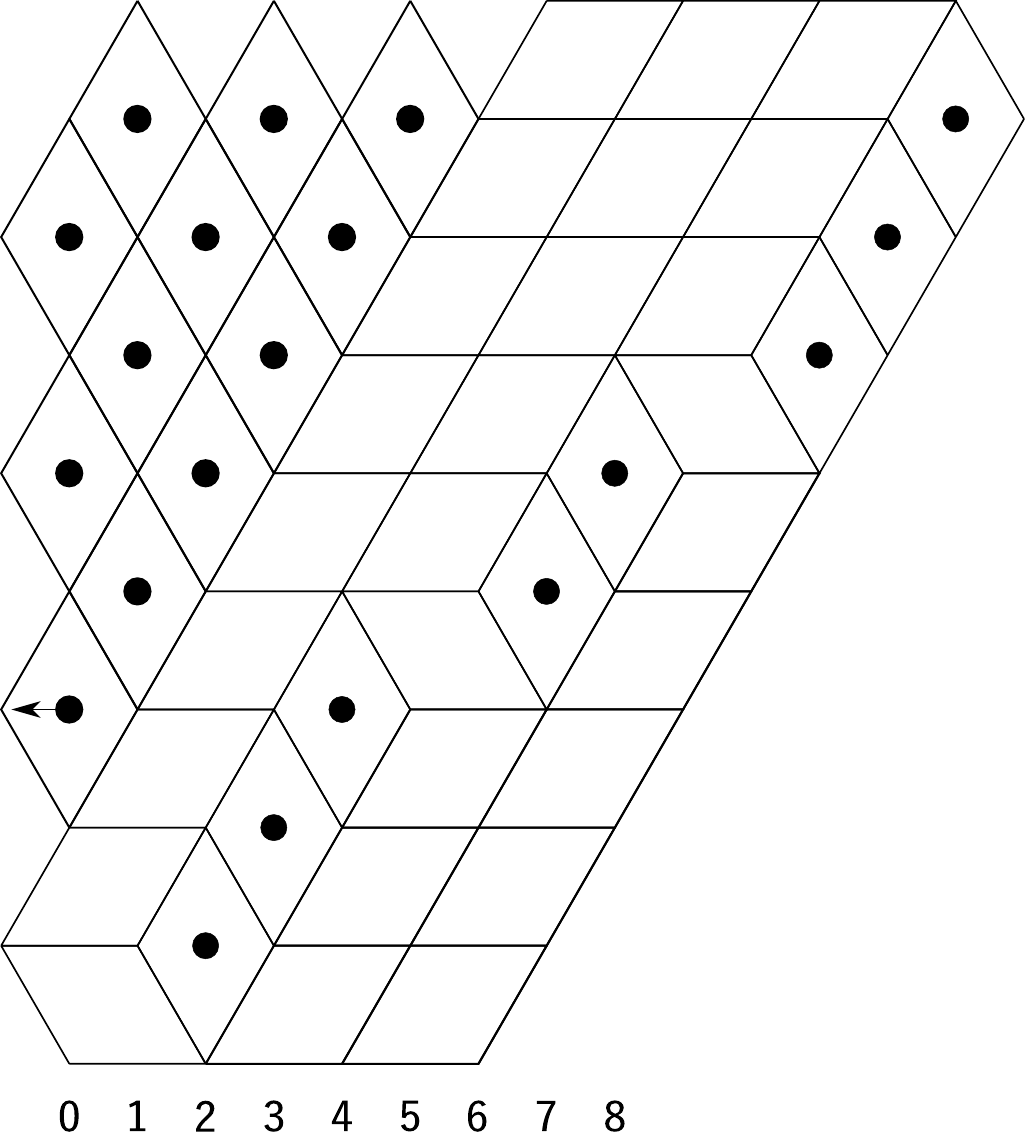} 
\caption{Three possible initial steps of either orthogonal or symplectic Plancherel growth. The arrows indicate the direction that the particle will attempt to jump next; see Figure \ref{modelcomparison} for the fourth step. The tiling determined by their positions suggests a stepped-surface interpretation.} \label{dynamicsteps}
\end{figure}

\emph{Plancherel growth processes} are certain continuous dynamics for Gelfand-Tsetlin patterns that can serve as prototypes for studying probabilistic phenomena in systems involving short-range dependence and simple interactions (for us, blocking and pushing). Figure \ref{dynamicsteps} depicts an instance of the dynamics studied in this paper.
\begin{figure}[h!]
\centering
\hspace{5em}
\includegraphics[width=0.24\textwidth]{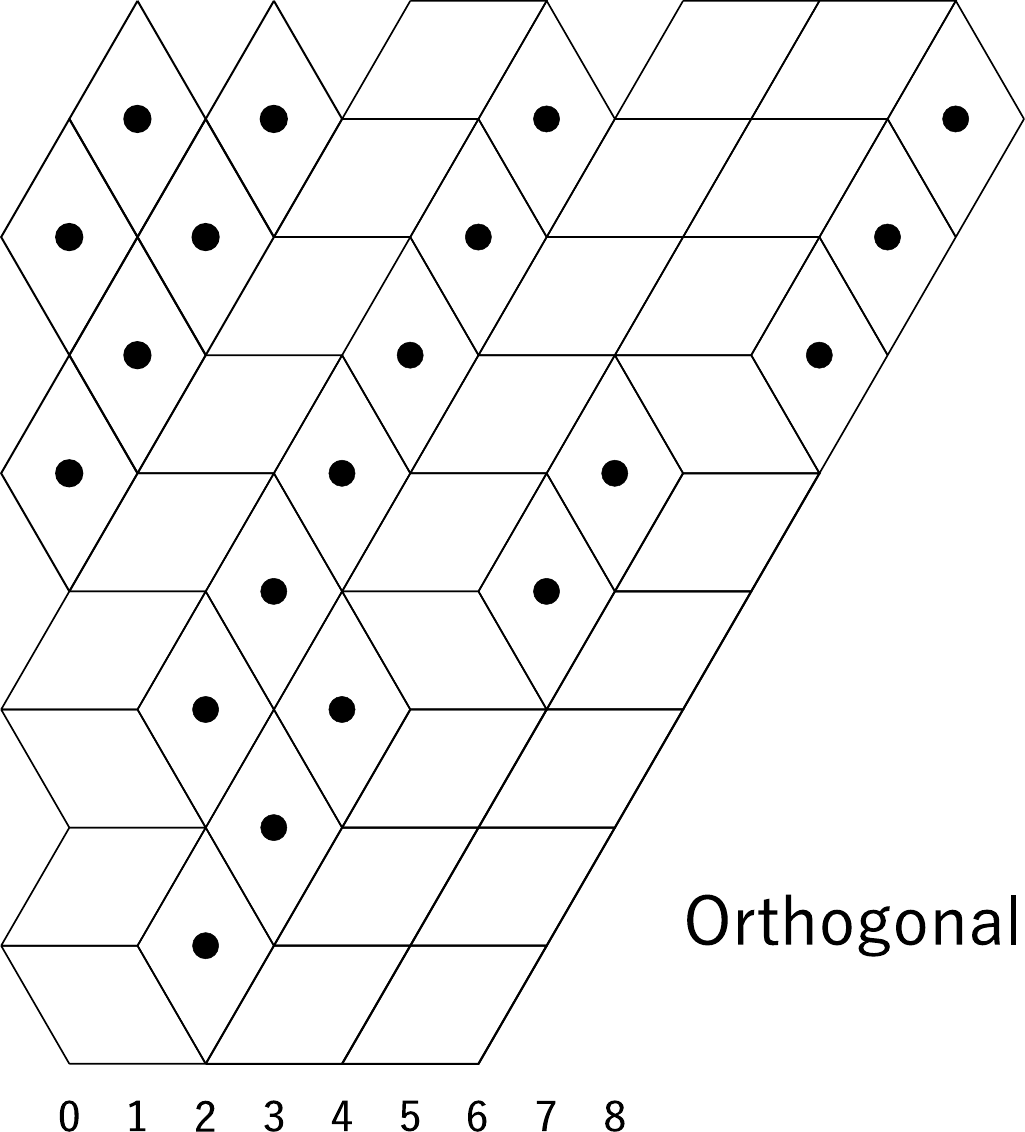}
\hspace{10em}
\includegraphics[width=0.24\textwidth]{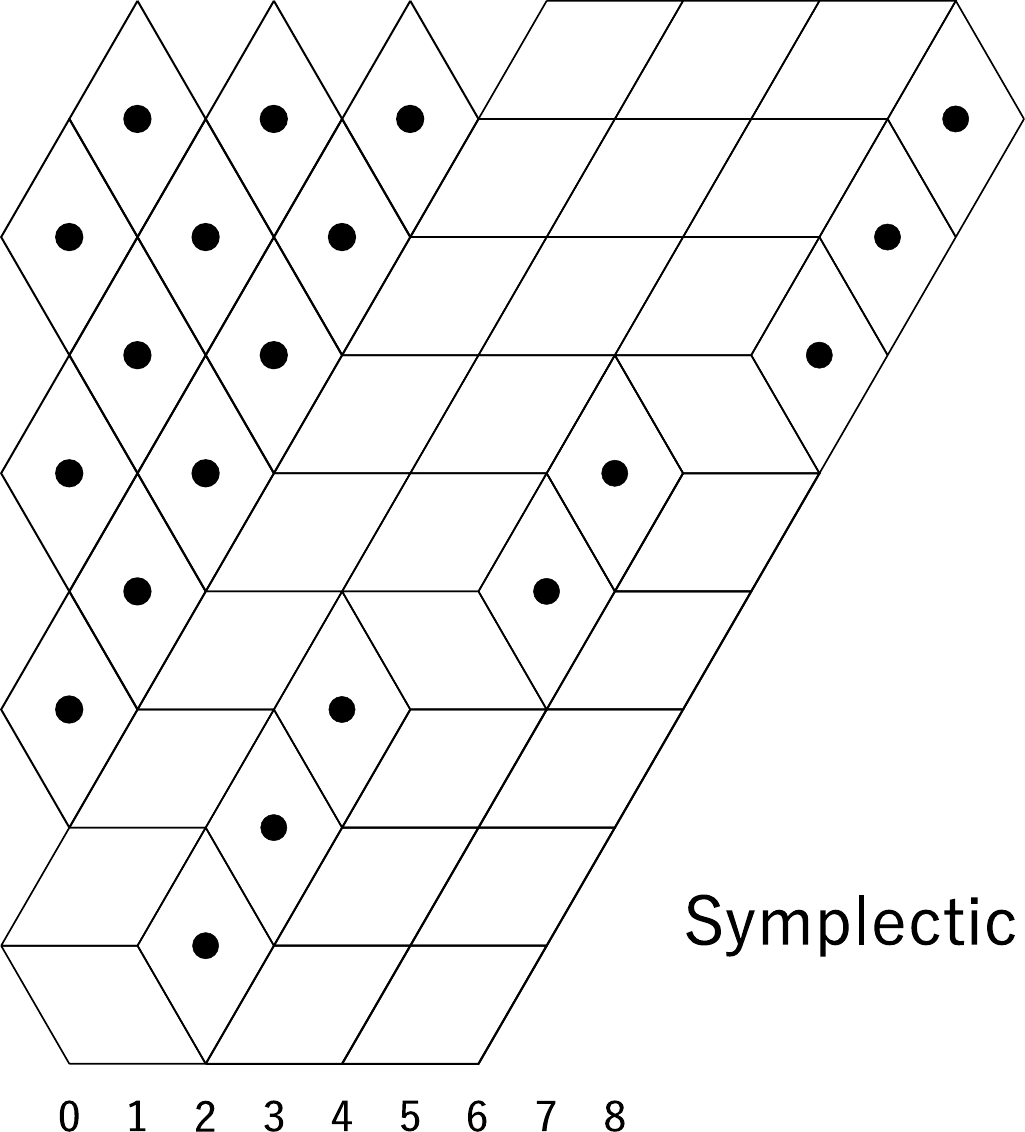}  
\hspace{5em}
\vspace{1.00mm}
\qquad
\newline
\centering
    \includegraphics[width=0.37\textwidth]{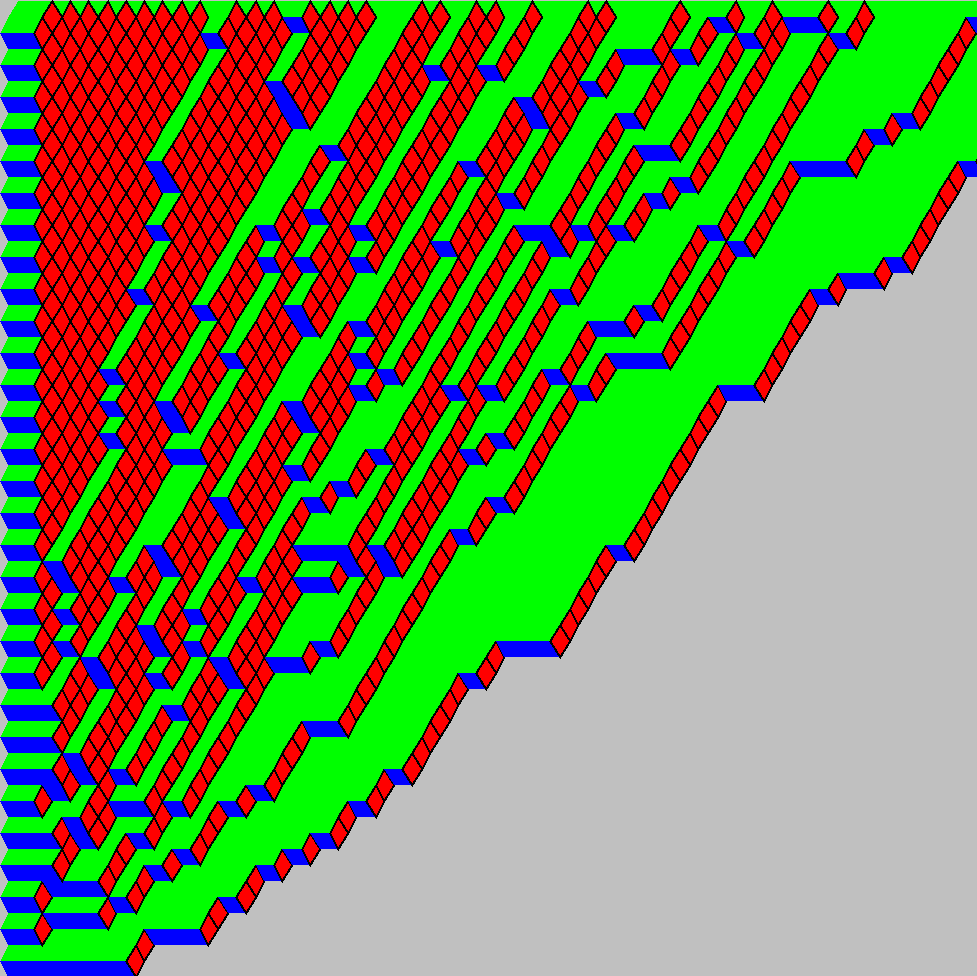}
    \hspace{5em}
    \includegraphics[width=0.37\textwidth]{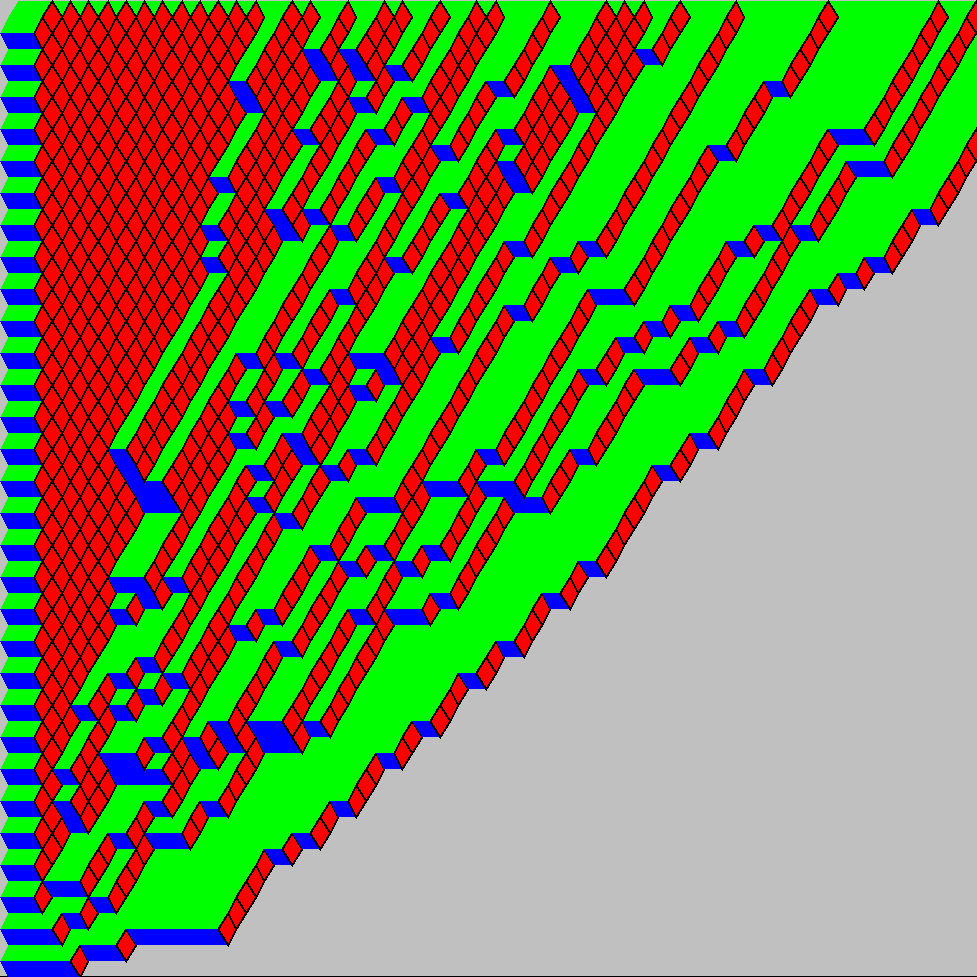}
\caption{The fourth step illustrates the distinctive wall behaviors: reflection (left) and suppression (right). Simulations of the orthogonal and symplectic cases indicate similar average behavior.} \label{modelcomparison}
\end{figure}
In our orthogonal/symplectic patterns, level $n \in \Z_{> 0}:=\{1,2,3 \ldots \}$ has $r_n:=\lfloor (n+1)/2 \rfloor$ particles $y^n_i$, $1 \leq i \leq r_n$, in $\Z_{\geq 0}:=\{ 0,1,2,\ldots \}$. \emph{Symplectic Plancherel growth} features these particles performing independent simple random walks, but a unit jump of a particle $y^n_i$ is suppressed or forced if otherwise it would land outside of the interval (write $y^{n-1}_0 \equiv \infty$)
$$
\begin{cases}
\left [0,y^{n-1}_{i-1} \right ), & i=r_n \ \text{with}  \ n \geq 1 \ \text{odd} \\
\left(y^{n-1}_i, y^{n-1}_{i-1} \right), & \text{otherwise for $1 \leq i \leq r_n$} \\
\end{cases}.
$$
In particular, the ``wall-particles" $y^n_{r_n}$ at odd levels $n \geq 1$ are $\mbf{suppressed}$ by the wall at $0$. This blocking and pushing of a particle by straddling particles on the previous level (including $0$ and $\infty$) serves to maintain \emph{interlacing conditions}, which in the chosen coordinates read as
\begin{equation} \label{coordinatization2}
\begin{cases}
0 \leq y^n_i < y^{n+1}_i, & i = r_n \ \text{with} \ n \geq 1 \ \text{odd} \\
y^{n+1}_{i+1} <  y^n_i < y^{n+1}_i, & \text{otherwise for $1 \leq i \leq r_n$}
\end{cases}.
\end{equation}
\emph{Orthogonal Plancherel growth} studied by Borodin--Kuan \cite{abjk1} differs only by the $\mbf{reflection}$ of jumps into the wall at $0$. Although the difference between the two systems may seem negligible at first, note that a single wall-particle's movement can cause the models to diverge by infinitely many particle positions; see Figure \ref{modelcomparison} for a comparison. Since local limit behavior of such integrable probability models can be sensitive to their basic parameters (notably, to initial conditions; see, e.g., \cite{bfs1}), it is worth investigating to what degree the distinct wall behavior is felt in the correlation structure and in various asymptotic regimes. To explain terminology, when deriving these dynamics via representation theory, the wall behavior originates from the initial conditions of orthogonal polynomials used to express irreducible/indecomposable characters of even-orthogonal/odd-symplectic groups; see Section \ref{continuoustime}. This remark indicates another, a priori significant departure from the orthogonal case: the relevance of the representation theory of non-semisimple Lie groups intermediate in the classical symplectic series. For readers less familiar with these algebraic items, note that all necessary concepts are reviewed and that Section \ref{sectiondiffusionlimit} relies only on modern probability theory.

\subsection{Main goals and motivation}

We have two primary goals:
\begin{enumerate} 

\item To quantify repulsion in the Dyson non-colliding dynamics achieved in the diffusion limit in order to exhibit non-degeneracy of constraint in this regime.

\item To compute series of Plancherel measures for $Sp(\infty)$ that reflect the intermediate representation theory of the non-semisimple odd symplectic groups, and subsequently investigate their role in the correlation structure and local limit behavior of symplectic Plancherel growth.

\end{enumerate}

In pursuing the first goal, we resolve a basic version of the conjecture posed after Proposition 2 of \cite{bfpsw1}, confirm the decompositions $(4),(29)$ in \cite{jw1} when the interlaced system starts at the origin, and thus also lift the restriction in \cite{vgms1} that the limiting system must start in the interior of the chamber. To show that the amount of constraint in the pre-limit system does not degenerate in this regime, we need to quantify the repulsion between particles coinciding in the limit. The path property that assumes this role is interesting in its own right: it delineates a lower envelope at $t=0$ for a multivariate self-similar process with positive components that do \emph{not} drift to infinity almost surely. For such $\R_+$-valued Markov processes drifting to infinity, see \cite{cp1} for rather general techniques and integral tests. The results of Section \ref{sectionpathproperty} at the heart of this development are new in this literature and are inspired by Burdzy--Kang--Ramanan \cite{bkr1}; we otherwise briskly follow the elegant procedures of Gorin-Shkolnikov \cite{vgms1} and Warren \cite{jw1}. For the related case of vicious random walkers, see the series of results by Katori, Tanemura, et.al., whose paper \cite{mkht2} covers convergence of our even levels.

For the second goal, the non-semisimplicity of odd symplectic groups is a potentially serious obstacle for the identification of suitable ``intermediate" Plancherel measures for $Sp(\infty)$ -- such measures can no longer arise as coefficients of \emph{irreducible} characters in the orthogonal decomposition of restricted infinite extreme characters. Circumventing this lack of irreducibility relies on the independent works of Proctor \cite{rp1} and Shtepin \cite{vs1}, and we suspect there are more connections in this vein to be made between probability theory and nonclassical representation theoretic objects. Defosseux \cite{md1} and Warren--Windridge \cite{ww1} have already considered our dynamics for symplectic Gelfand-Tsetlin patterns, with the latter paper anticipating our interest, but the results we derive here are much in the same spirit as Borodin--Kuan \cite{abjk1}, who cover the orthogonal case. Obtaining refined local limit behavior in their manner is well codified, originating with a law of large numbers result for the limit shape of random Young tableaux, derived independently by Vershik--Kerov \cite{vk1} and Logan--Shepp \cite{ls1}. The former pair were studying asymptotics of the maximum dimension of irreducible representations for the symmetric group, and the latter were pursuing a variational problem arising in combinatorics posed by Richard Stanley. Baik--Deift--Johansson \cite{bdj1}, and shortly after Okounkov \cite{ao1}, derived more detailed results concerning local behavior at the edge of the limit shape, inspiring much research on such behavior, including the present work; for a perfunctory list, consider \cite{ bdj1, bb1, bbo1, abic1, abpf2, abpf1, bfs1, abjk2, abjk1, boo1, bo1, abgo1, jk1, jk2, ls1,  ao1, aogo1,aonr1,  vk1}. 

\begin{convention} 
Quantities with indices that overflow are set to $0$ and those with indices that underflow are set to $\infty$. 
\end{convention}

\subsection{Main applications}

\subsubsection{Diffusive scaling limit: non-degeneracy of constrained dynamics} \label{introdiffusion}

The basic description of symplectic Plancherel growth (considered up to level $n \geq 1$) readily leads to the candidate diffusive limiting system $(\mc{W}(t))_{t \geq 0}$ with components $\mc{W}^k_i$ satisfying 
\begin{equation} \label{continuumlimit}
\mc{W}^k_i \ \text{is a Brownian motion $B_i^k$ reflected in the interval} \ 
\begin{cases}
[0,\infty), &\ \  \text{if} \ k=1 \\
[\mc{W}^{k-1}_i, \mc{W}^{k-1}_{i-1}], & \ \ \text{if} \ k \geq 2
\end{cases}
\end{equation}
for $1 \leq k\leq n$, $1 \leq i \leq r_k$, where the $B^k_i$ are independent and started at $0$; note our notational convention puts additional static particles at $0$ (for index overflow) and at $\infty$ (for index underflow). Section \ref{definitionreflection} makes the notion of reflecting Brownian movements in a time-dependent domain rigorous, but we say the coupled dynamics of $(\mc{W}(t))_{t \geq 0}$ are \emph{non-degenerate} if its components uniquely solve the system of equations
\begin{equation} \label{continuumlimit1}
\mathcal{W}^k_i =
\begin{cases}
 B_i^k+ \frac{1}{2} L^0_{\mathcal{W}^k_{i}} - L^0_{(\mathcal{W}^{k-1}_{i-1}-\mathcal{W}^k_i)}, & i=r_k \ \text{with} \ k \ \text{odd} \\
 B_i^k+  L^0_{(\mathcal{W}^k_i - \mathcal{W}^{k-1}_i)} - L^0_{(\mathcal{W}^{k-1}_{i-1}-\mathcal{W}^k_i)} , & \text{otherwise for} \ 1 \leq i \leq r_k
\end{cases},
\end{equation} 
where $L_{Y}^0$ is twice the semimartingale local time at $0$ of a semimartingale $Y$ and where $L^0_\infty \equiv 0$ (see Remark 4.2 of \cite{bkr1} for an explanation of the appropriate local time scalings). Intuitively, these equations indicate the amount of constraint in the system $(\mc{W}(t))_{0 \leq t \leq T}$ is finite for every $T \geq 0$. But demonstrating non-degeneracy of $(\mc{W}(t))_{t \geq 0}$ requires that we resolve a subtle yet concrete difficulty.
\begin{figure}[h!!] 
\centering
\includegraphics[width=0.4\textwidth]{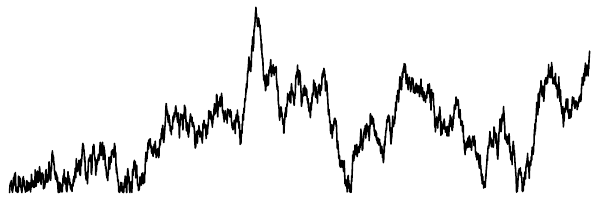} 
\qquad
\includegraphics[width=0.4\textwidth]{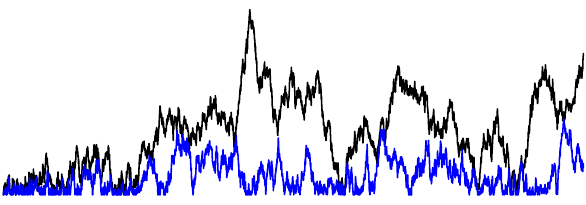} 
\caption{The first path in black is a standard brownian motion reflected at $0$. Reflecting the blue path to stay between $0$ and the black path requires an ``infinite amount" of constraining to achieve.} \label{norepulsion}
\end{figure}
\begin{figure}[h!!]
\centering
\includegraphics[width=0.45\textwidth]{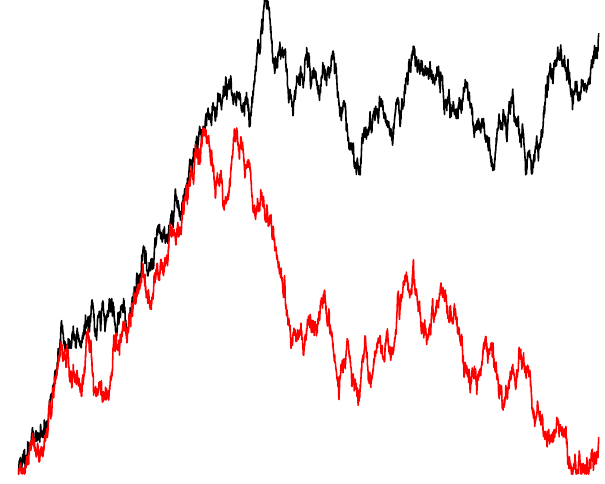} 
\includegraphics[width=0.45\textwidth]{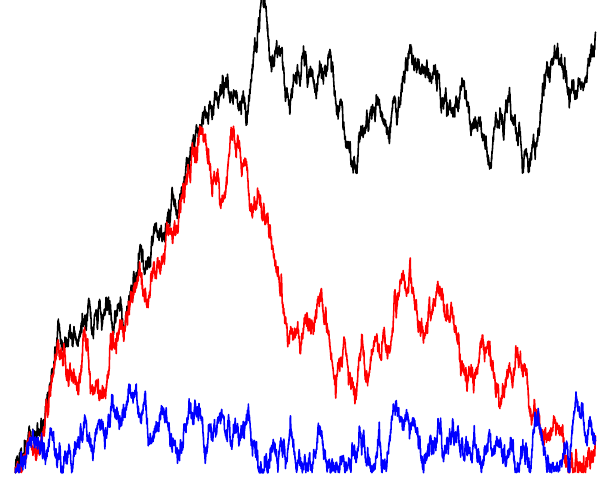} 
\caption{Diffusion limit of first three particles nearest the wall (visualized in the single state space $\R_+$). The level 1 particle in red creates a repulsion between the wall at $0$ and the level 2 particle in black; although this extra room may still not be enough, we show the reflected level 3 wall-particle in blue has non-degenerate constrained dynamics.} \label{repulsion}
\end{figure}

To see the issue in our case, write the dynamics of a typical reflected particle in the basic decomposition $W = B+Y$, where the constraining process $Y$ ``pushes" the Brownian particle $B$ to remain between some independent particles $\ell, r$ on the previous level. If ever the interval $[\ell(t),r(t)]$ collapses, $Y$ might need to ``push an infinite amount" in order to keep $B$ at that point. For example, $Y$ can be finite variation on every period between collapse times of $[\ell(t),r(t)]$, but accumulate infinite variation over infinitely many such periods in finite time; see Figure \ref{norepulsion} for such a scenario.

For a more surprising and relevant example, the interval $[\ell(t),r(t)]$ can ``close too quickly" at $t=0$, say, for some $c>1$, almost surely
\begin{equation} \label{closefast}
\limsup_{t \downarrow 0} \frac{r(t)-\ell(t)}{t^c} = 0;
\end{equation}
in this situation, even if the particles $\ell(t), r(t)$ remain apart for $t>0$, $Y$ will not be of bounded variation. If $Y$ is of infinite variation, then the symbols indicating local times in \eqref{continuumlimit1} do not make sense (even though their increments may). Proposition 4.13 and the proof of Proposition 4.12 of \cite{bkr1} provide more details on these two degenerate examples. 

Fortunately, by the construction of our system \eqref{continuumlimit}, any particles $\ell, r$ that can collide must straddle an independent Brownian particle on a further previous level that could create a sufficiently strong repulsion that $Y$ only needs to ``push a finite amount" to keep $B$ between them; see Figure \ref{repulsion}. Our strategy, then, consists of quantifying this degree of repulsion in order to rule out the degenerate possibility \eqref{closefast} and exhibit non-degeneracy of the constrained dynamics $(\mc{W}(t))_{t \geq 0}$ of \eqref{continuumlimit}. More broadly, Section \ref{sectiondiffusionlimit} will establish the following result.
\bt \label{pathproperty} 
Fix $n \geq 1$. The interlaced system $(\mc{W}(t))_{t \geq 0}$ of reflecting Brownian motions with a wall defined by \eqref{continuumlimit} has levels $\mc{W}^{k}$, $1 \leq k \leq n$, that achieve the Dyson non-colliding dynamics
\begin{equation} \label{dysoncd}
dX^{k}_i = 
\begin{cases}
dB_i + \frac{dt}{X^{k}_i} + \sum_{\substack{j \neq i \\ 1 \leq j \leq r_k}}  \left ( \frac{1}{X_i^{k} - X_j^{k}} +  \frac{1}{X_i^{k} + X_j^{k}}  \right ) dt,  & \text{if} \ \ k \ \ \text{even} \\
dB_i + \frac{1}{2} 1_{(i=r_k)} dL_{X^k_i}^0(t) + \sum_{\substack{j \neq i \\ 1 \leq j \leq r_k}}  \left ( \frac{1}{X_i^{k} - X_j^{k}} +  \frac{1}{X_i^{k} + X_j^{k}}  \right ) dt, & \text{if} \ \ k \ \ \text{odd}
\end{cases},
\end{equation}
for $1 \leq i \leq r_k$, where the $B_i$ are independent standard Brownian motions. Moreover, each particle $\mc{W}^{k+1}_i$, $2 \leq k < n$, $1 < i \leq r_{k+1}$ is reflected in an interval $[\mc{W}^{k}_{i}(t), \mc{W}^{k}_{i-1}(t)]$ satisfying the almost sure path property
$$
\liminf_{t \downarrow 0} \frac{\mc{W}^{k}_{i-1}(t) - \mc{W}^{k}_{i}(t)}{t^c} = \infty
$$
for any $c>1/2$; $(\mc{W}(t))_{t \geq 0}$ is then verifiably the unique, strong solution to \eqref{continuumlimit1} when started at the origin. Hence, under diffusive space-time scaling, symplectic Plancherel growth considered up to level $n$ converges in law to the non-degenerate system \eqref{continuumlimit1}.
\et
The path property above is readily seen to hold in Warren's setting \cite{jw1} of interlaced reflecting Brownian motions (without a wall), which achieve classical Dyson dynamics
\begin{equation} \label{dyson1}
d\lambda^i_t = dB^i_t + \frac{\kappa}{2} \sum_{j \neq i} \frac{dt}{\lambda^i_t - \lambda^j_t}, \ \ \ 1 \leq i \leq N,
\end{equation}
where the $B^i$ are independent. In particular, setting $\kappa =2$ and arguing as in Section \ref{sectionpathproperty} confirms the decompositions $(4),(29)$ in \cite{jw1} when the system starts at the origin and thus also lifts the restriction in \cite{vgms1} that the limiting system must start in the interior of the chamber $\{ x_1 \geq \ldots \geq x_N\}$. Note the path property for \eqref{dyson1} is compelling only for $N \geq 3$; in the case $N=2$, the single gap is proportional to a $3$-dimensional Bessel process and the result is classical.

\subsubsection{Determinantal correlation kernel and local limit behavior} \label{introdeterminantal}

As we saw in the last section, the diffusion limit primarily cares inductively about interactions of three successive levels at a time. To capture the effects of correlations between particles farther apart in the system, we now view symplectic Plancherel growth as determining a point configuration $\mc{X}(\gamma)$, $\gamma \geq 0$, in $\Z_{\geq 0} \times \Z_{> 0}$, i.e., as a random element of $2^{\Z_{\geq 0} \times \Z_{> 0}}$. A fundamental theorem about this point process is that it is determinantal  with an explicitly computable kernel when started from the ``leftmost" initial condition (see the first image of Figure \ref{dynamicsteps}). To state this result and others that follow, let $J_{k,1/2}(x)$, $J_{k,-1/2}(x)$, $k \geq 0$ denote the Chebyshev polynomials of the second and third kind, respectively, (see Section \ref{sectionnotations}) and write
$$
\alpha_k:=
\begin{cases}
1/2, & \text{if} \ \ k \ \text{even} \\
-1/2, & \text{if} \ \ k \ \text{odd}
\end{cases}.
$$
Consider the recoordinatization $\widetilde{\mc{X}}(\gamma)$ defined to have components $\widetilde{\mc{X}}^n_k(\gamma):= \mc{X}^n_k(\gamma)+r_n-(n-k+1)$, $n \geq 1$, $1 \leq k \leq r_n$ (see Section \ref{sectioncentral} for explanation and visualization). The next theorem is a consequence of the similar statement Theorem \ref{generalkernel} concerning the general symplectic Plancherel point process. 
\bt \label{kernel}
The $k$th correlation functions $\{ \widetilde{\rho}^\gamma_k\}_{k \geq 1}$ of $(\widetilde{\mc{X}}(\gamma) )_{\gamma \geq 0}$ are determinantal: for $z_1, \ldots, z_k \in \Z_{\geq 0} \times \Z_{> 0}$,
$$
\widetilde{\rho}^\gamma_k(z_1,\ldots,z_k) := \mathbb{P}(\widetilde{\mc{X}}(\gamma) \supset \{z_1,\ldots, z_k \})= \det[ K^\gamma(z_i, z_j) ]_{i,j=1}^k,
$$
where the (nonsymmetric) kernel $K^\gamma$ is given explicitly, for $(s,n),(t,m) \in \Z_{\geq 0} \times \Z_{> 0}$, by
\begin{equation} \label{explicitkerneleasy}
\ba
K^\gamma((s,n),(t,m)) =  \frac{2^{\alpha_n+1/2}}{\pi} & \frac{1}{2 \pi i}  \int_{-1}^1 \oint \frac{e^{\gamma x}}{e^{\gamma u}} J_{s,\alpha_n}(x) J_{t,\alpha_m}(u) \frac{(1-x)^{r_n + \alpha_n}(1+x)^{1/2}}{(1-u)^{r_m}(x-u)} dudx\\
& +1_{(n \geq m)} \frac{2^{\alpha_n+1/2}}{\pi}  \int_{-1}^1 J_{s,\alpha_n}(x) J_{t,\alpha_m}(x) (1-x)^{r_n-r_m+\alpha_n} (1+x)^{1/2} dx. 
\ea
\end{equation}
Here, the $u$ contour is a positively oriented (i.e., counterclockwise) simple loop around $[-1,1]$.
\et
\begin{figure}[h]
\centering
\includegraphics[width=0.4\textwidth]{symplecticlargetime.pdf} 
\hspace{.5em}
\includegraphics[width=0.327\textwidth]{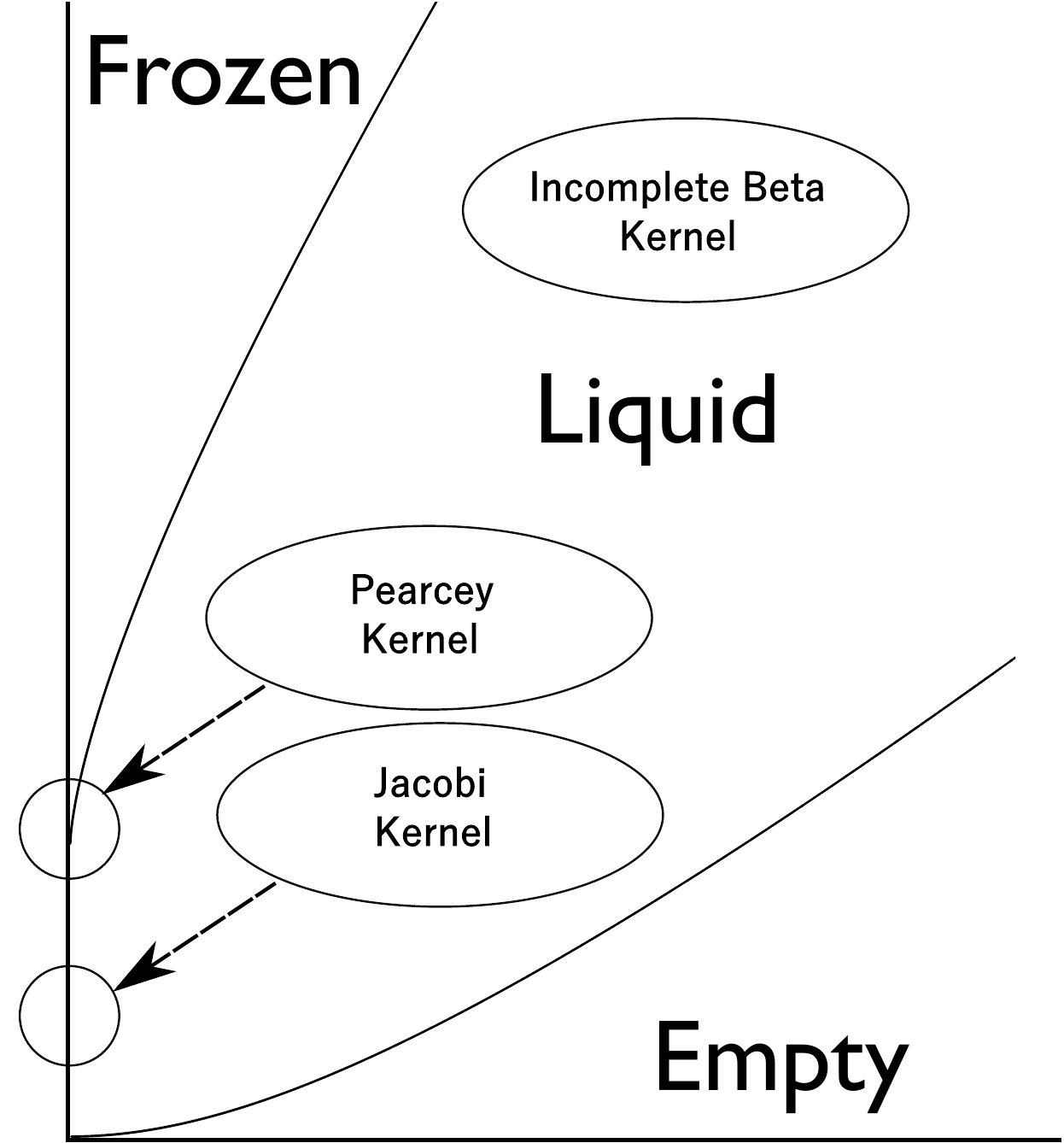} 
\caption{The hydrodynamic limit capturing average behavior indicates a densely-packed region of inactivity (``frozen") and a central region of activity (``liquid"), with the remaining area unreached (``empty"). Notably, the orthogonal and symplectic cases feature the same curves delineating these regions, which are given analytically by Proposition \ref{analyticlines}; compare with Figure 4 of \protect\cite{abjk1}.} \label{averagebehavior}
\end{figure}
Although the correlation kernel \eqref{explicitkerneleasy} of Theorem \ref{kernel} may appear unwieldy at first, its double integral expression allows for asymptotic expansions leading to refined limiting results. Define the region $\mc{D}^{liquid}$ of nontrivial hydrodynamic limiting behavior as the collection of locations and times with a positive chance of either seeing a particle or not:
\begin{equation} \label{liquid1}
\mc{D}^{liquid}: = \left \{ (\tau, \nu, \eta) \in \R_+ \times \R_+ \times \R_+ : 0< \lim_{N \to \infty} \widetilde{\rho}^{[N \tau]}_1(([N\nu],[N\eta])) <1 \right \}.
\end{equation}
Under this hydrodynamic scaling, the method of steepest descent applied directly to the kernel \eqref{explicitkerneleasy} shows that $\mc{D}^{liquid}$ is exactly the region of $(\tau, \nu, \eta) \in \R_+^3$ where the function
\begin{equation} \label{steepestdescent}
S(z) = S_{\tau,\nu,\eta}(z) : = \tau \frac{z+z^{-1}}{2}  - \nu \log z + \eta \log \left ( \frac{z+z^{-1}}{2} -1 \right ).
\end{equation}
has a unique critical point $z_0=z_0(\tau, \nu, \eta)$ in $\mb{H} \setminus \mb{D}$. Similarly define $\mc{D}^{frozen}$ to be the region where we are sure to find particles in this limit, and $\mc{D}^{empty}$ where there can be none; more precisely, these are the regions where the limit in \eqref{liquid1} is $1$ and $0$, respectively. See Figure \ref{averagebehavior} for a depiction of these regions. To identify the nontrivial limit in the liquid region, define the \emph{incomplete beta kernel} by
$$
B(k;l |\zeta): = \frac{1}{2\pi i} \oint_{\overline{\zeta}}^\zeta (1-z)^k z^{-(l+1)} dz,
$$
where the contour crosses $(0,1)$ if $k \geq 0$ and $(-\infty, 0)$ if $k<0$. This kernel serves as a generalization of the sine kernel to describe bulk limits of the Schur process in \cite{aonr1} by Okounkov--Reshetikhin.
\bt{(Bulk Limit: Incomplete Beta Kernel)} \label{bulklimit}
Let $(\tau, \nu, \eta) \in \mc{D}^{liquid}$ and let $z_0=z_0(\tau, \nu, \eta)$ be the unique critical point of \eqref{steepestdescent} in $\mb{H} \setminus \mb{D}$. Assume $\gamma \geq 0$, $(s_i, n_i) \in \Z_{\geq 0} \times \Z_{> 0}$, $1 \leq i \leq k$, depend on $N$ in such a way that $\gamma \sim N \tau >0$, $s_i \sim N \nu>0$, and $r_{n_i} \sim N \eta>0$. Assume the differences $s_i - s_j$, $n_i - n_j$ are of constant order. Then as $N \to \infty$,
$$
\widetilde{\rho}^\gamma_k  ((s_1, n_1), \ldots ,(s_k, n_k))
 \to
\det [ B(n_i - n_j; (s_i-s_j) + (n_i-n_j) - (r_{n_i}- r_{n_j})|z_0) ]_{i,j=1}^k.
$$
\et
At the edges of the liquid region $\mc{D}^{liquid}$, we find variants of the Jacobi and Pearcey kernels of Borodin--Kuan \cite{abjk2, abjk1}. The first of these types arises in the large time limit at a finite distance from the wall, away from the corner of the phase transition; it only seems to have appeared so far in the similar contexts of \cite{abjk1} and \cite{jk1}, though with different Jacobi polynomials.

\bt{(Edge Limit: Discrete Jacobi Kernel)} \label{jacobilimit}
Assume $\gamma \geq 0$, $(s_i, n_i) \in \Z_{\geq 0} \times \Z_{> 0}$, $1 \leq i \leq k$, depend on $N$ in such a way that $\gamma \sim N  \tau >0$, $r_{n_i} \sim N  \eta>0 $ but the $s_i$ are fixed and finite. Assume only the differences $n_i - n_j$ are of constant order. Then as $N \to \infty$,
$$
\widetilde{\rho}^\gamma_k  ((s_1, n_1), \ldots ,(s_k, n_k)) \to
\begin{cases}
\det [ L((s_i, n_i), (s_j, n_j)|1-\eta/\tau) ]_{i,j=1}^k, & \text{if} \ \ 1-\eta/\tau > -1 \\
1, & \text{if} \ \ 1-\eta/\tau \leq -1
\end{cases}
$$
where the discrete Jacobi kernel $L$ is given explicitly, for $(s,n),(t,m) \in \Z_{\geq 0} \times \Z_{> 0}$, by\begin{equation} \label{explicitjacobi}
\ba
L((s,n), (t,m) | \epsilon) : = & \frac{2^{\alpha_n+1/2}}{\pi} \int_{-1}^1 [1_{(n\geq m)}+ 1_{[-1,\epsilon]}(x)] J_{s,\alpha_n}(x)J_{t,\alpha_m}(x)(1-x)^{r_n-r_m+\alpha_n}(1+x)^{1/2} dx\\
\ea.
\end{equation}
\et
Our last result shows the corner of the phase transition (i.e., near where the regions $\mc{D}^{frozen}$ and $\mc{D}^{liquid}$ meet the wall) exhibits a $1/4$-exponent in its fluctuation scale, providing additional rigorous evidence of the universality of this growth exponent in the class of anisotropic KPZ with a wall. The evolving stepped-surface interpretation indicated by Figure \ref{dynamicsteps} allows one to define a height function as ``the number of particles at that point and to the right". The fluctuations about its mean should be described (at least formally) by the \emph{anisotropic KPZ equation}, a stochastic heat equation of the form
$$
\partial_t h =  \nu_x \partial_x^2 h + \nu_y \partial_y^2 h + \frac{1}{2}\lambda_x (\partial_xh)^2 + \frac{1}{2}\lambda_y (\partial_y h)^2 + \xi,
$$
where $\xi$ is a space--time white noise and $\lambda_x, \lambda_y$ have different signs (should they have the same sign, this equation reduces to the classical KPZ equation). We do not explore this equation in this work, we only note its relevance for future research; see, e.g., Borodin--Ferrari \cite{abpf1} for more discussion. Such rigorous results for systems with a wall still seem to be rare; the author only knows of \cite{bnw1, abjk1, jk1, kf1}. The kernel arising in our setting at this critical corner point appears to be an instance of a general family of Pearcey-type posited in Remark 1.5 of \cite{kf1}.
\bt{(Edge Limit: Pearcey Kernel)} \label{pearceylimit}
Let $\widetilde{\rho}^{\gamma, k}_\Delta$ denote the $k$th correlation function of the determinantal process determined by the complementary process $(\widetilde{\mc{X}}(\gamma))^c$. Assume $\gamma \geq 0$, $(s_i, n_i) \in \Z_{\geq 0} \times \Z_{> 0}$, $1 \leq i \leq k$, depend on $N$ in such a way that $\gamma \sim N/2$, $s_i \sim N^{1/4} \nu_i>0$, and $r_{n_i} -N\sim \sqrt{N} \eta_i>0$. Then as $N \to \infty$, 
$$
\ba
 N^{k/4}  \cdot \widetilde{\rho}^{\gamma, k}_\Delta & ((s_1, n_1), \ldots ,(s_k, n_k))
\to
\det [ \mc{K}((\nu_i, \eta_i), (\nu_j, \eta_j)) ]_{i,j=1}^k,
\ea
$$
where the Pearcey Kernel $\mc{K}$ is given explicitly by
$$
\ba
\mc{K}((\nu_1, \eta_1), (\nu_2, \eta_2)) := & 1_{(\eta_1 > \eta_2)}  \frac{1}{\sqrt{\pi (\eta_1 - \eta_2)}} \left ( \exp \left [ - \frac{  (\nu_1 + \nu_2)^2}{\eta_1 - \eta_2} \right ] - \exp \left [ - \frac{ (\nu_1 - \nu_2)^2}{\eta_1 - \eta_2} \right ] \right )  \\
& +  \frac{\sqrt{2}}{ \pi}  \frac{1}{2 \pi i}  \int_{-i\infty}^{i\infty} \int_0^\infty e^{ \frac{(u')^2-(x')^2}{8} + \frac{\eta_2 u' - \eta_1 x'}{2}}\sin[\nu_1 \sqrt{2 x'} ] \sin[\nu_2 \sqrt{2 u'} ] \frac{dx'du'}{\sqrt{u'}(u'-x')}.
\ea
$$
\et

\begin{acknowledgements}
The author would like to thank Ren\'{e} Carmona, Allen Knutson, Daniel Lacker, Ramon Van Handel, and Mykhaylo Shkolnikov for many helpful discussions, as well as Jon Warren for an informative private correspondence. I would especially like to thank Jeffrey Kuan, for early encouragement and consistently valuable guidance, as well as for supporting a simulation of the symplectic growth process at \protect \url{http://www.math.harvard.edu/~jkuan/Cerenzia.html}. 

Part of this research was performed while the author was visiting the Institute for Pure and Applied Mathematics (IPAM), which is supported by the National Science Foundation; the author was also partially supported by NSF: DMS-0806591.
\end{acknowledgements}

\section{Diffusive scaling limit and repulsion of Dyson gaps} \label{sectiondiffusionlimit}

\subsection{Time-dependent Skorokhod problem and Skorokhod reflection map} \label{definitionreflection}
The machinery of the extended Skorokhod map provides a remarkably natural and careful construction of our particle system and its diffusion limit.  A solution to the \emph{Skorokhod problem} (SP) of reflecting a particle $\psi$ in a (potentially time-dependent) interval $[\ell,r]$, $\ell \leq r$, is a pair $(\phi,\eta)$ where $\eta$ ``pushes" the particle $\psi$ to remain in $[\ell,r]$ and the resulting reflected movements $\phi = \psi + \eta \in [\ell,r]$ should be regular; classically, $\eta$ should be of bounded variation and $(\phi, \eta)$ should satisfy the familiar conditions
 \begin{equation} \label{skorokhodconditions}
 \int_0^\infty 1_{(\ell(s)<\phi(s))} d \eta_\ell(s) =  \int_0^\infty 1_{(\phi(s) < r(s))} d \eta_r(s) = 0,
\end{equation}
where $\eta_\ell, \eta_r$ are nondecreasing functions affording the (minimal) decomposition $\eta =  \eta_\ell- \eta_r$. But as explained in Section \ref{introdiffusion}, the conditions \eqref{skorokhodconditions} are too restrictive for our situation. To state a weaker notion of solution, let $D(\R_+, S)$ be the space of right-continuous functions with left limits on $\R_+$ that take values in a separable complete metric space $S$. Write $D := D(\R_+, \R)$, $D^- := D(\R_+, [-\infty, \infty))$, and $D^+ := D( \R_+, (-\infty, \infty])$. 

\bd{(Definition 2.2, \cite{bkr1})} \label{ESP}
Given $(\ell,r, \psi) \in D^- \times D^+ \times D$ with $\ell \leq r$, a pair $(\phi, \eta) \in D \times D$ solves the \emph{extended Skorokhod problem} (ESP) on $[\ell, r]$ for $\psi$ if

\begin{enumerate}

\item $\phi(t) = \psi(t) + \eta(t) \in [\ell(t),r(t)]$ for $t \geq 0$

\item For $0\leq s<t<\infty$, 

\subitem $\eta(t) - \eta(s) \geq 0$, if $\phi < r$ on $(s,t]$; and $\eta(s) - \eta(s-) \geq 0$, if $\phi(s) < r(s)$  

\subitem $\eta(t) - \eta(s) \leq 0$, if $\ell< \phi$ on $(s,t]$; and $\eta(s) - \eta(s-) \leq 0$, if $\ell(s)<\phi(s)$

\end{enumerate}
\ed
Given $(\ell,r, \psi) \in D^- \times D^+ \times D$ with $\ell \leq r$, an explicit candidate for the reflected process $\phi$ is provided by the \emph{extended Skorokhod map} $ \overline{\Gamma}(\ell,r|\psi) : = \psi - \Xi_{\ell,r}(\psi)$, where
$$
\ba
\Xi_{\ell,r} & (\psi)(t) : = \max \bigg (  (\psi(0) - r(0))^+ \wedge \inf_{u \in [0,t]}(\psi(u) - \ell(u)), \sup_{s \in [0,t]} [ (\psi(s) - r(s)) \wedge \inf_{u \in [s,t]}(\psi(u) - \ell(u)) ] \bigg ).
\ea
$$
Note some simplifications: if $r \equiv \infty$, then $\Xi_{\ell,\infty}(\psi)(t) = - \sup_{s \in [0,t]}[\ell(s)-\psi(s)]^+$, $t \geq 0$; also, at the initial time $0$,
\begin{equation} \label{simplification}
\Xi_{\ell,r}(\psi)(0) = [\psi(0) - r(0)]^+ \wedge [\psi(0) - \ell(0)]^{-}.
\end{equation}
In words then, $\overline{\Gamma}(\ell,r|\psi)$ starts by putting the particle $\psi$ at the nearest point in $[\ell,r]$ and maintains its increments except for putting those that land outside of $[\ell,r]$ again at the nearest point inside. Theorem 2.6 and Proposition 2.8 of \cite{bkr1} show $\overline{\Gamma}$ is a jointly continuous function on $D^- \times D^+ \times D$ (each topologized by uniform convergence on compacts) and provides the \emph{unique} solution $$(\phi,\eta):= ( \overline{\Gamma}(\ell,r|\psi) ,  \overline{\Gamma}(\ell,r|\psi)-\psi)$$ to the associated ESP on $[\ell,r]$ for $\psi$. The condition $\inf_{t \in [t_1,t_2]} (r(t) - \ell(t))>0$, $0 \leq t_1 < t_2$, is sufficient for $(\phi, \eta)$ to satisfy the more regular conditions \eqref{skorokhodconditions} on $[t_1,t_2]$; cf. Corollary 2.4 of \cite{bkr1}.

\subsection{First construction of symplectic Plancherel growth and continuum universal limit} \label{firstconstruction}

Fix $n \geq 1$. Consider the collection $X(t) : = \{ X^k_i(t) | 1 \leq k \leq n, \ 1 \leq i \leq r_k \}$, $t \geq 0$, of continuous time simple random walks $X^k_i$ started at $0$, with jumps dictated by independent unit rate poisson processes. For a parameter $N>0$, consider the diffusively scaled versions $\bar{X}(t;N):= X_{tN}/\sqrt{N}$. Construct a process $(\mc{X}(t,N))_{t \geq 0}$ with components $\mc{X}^k_i(t,N)$ \emph{driven} by $\bar{X}$ by setting $\mc{X}^1_1 := \overline{\Gamma}(0, \infty| \bar{X}^1_1)$ and inductively for $2 \leq k \leq n$,
\begin{equation} \label{discretesk}
\mc{X}^k_i :=
\begin{cases}
\overline{\Gamma}  (0,\mc{X}^{k-1}_{i-1}-1/\sqrt{N} \big|  \bar{X}^k_{i} ), & \text{if} \ \ i =r_k \ \ \text{with} \ \ $k$ \ \text{odd} \\
\overline{\Gamma}(\mc{X}^{k-1}_{i}+1/\sqrt{N} ,\mc{X}^{k-1}_{i-1}-1/\sqrt{N} \big|   \bar{X}^k_{i}), & \text{otherwise for} \ 1\leq i \leq r_k
\end{cases},
\end{equation} 
where our convention $\mc{X}^{k-1}_{0} \equiv \infty$ applies. This procedure determines \emph{deterministic} maps $\Phi^{SK}_N(\bar{X})(t) = \mc{X}(t,N)$.
\bd \label{construction1}
\emph{Symplectic Plancherel growth considered up to level $n \geq 1$} is characterized by the interlaced dynamics of the process $\mc{X}(\gamma):= \mc{X}(\gamma,1)$, $\gamma \geq 0$, defined by \eqref{discretesk}. 
\ed
\bn
By the simplification \eqref{simplification} at time $t=0$, $\mc{X}(0)$ is exactly the ``leftmost" initial condition up to level $n$. However, Definition \ref{construction1} refers to the interlaced dynamics and does not depend on this initial condition. In fact, the paper \cite{vgms1} illustrates how the results of this section hold for much more general dynamics and initial conditions, and a forthcoming technical report will show the dynamics of \cite{abjk1, md1, md2, jk1,ww1} all converge to the same system under diffusive scaling. It is only for continuity of focus and simplicity of presentation that we maintain our concrete setting. 
\en

The appropriate continuum state space for the diffusion limit is given by  the collection $\mb{J}^c_{n,paths}$ of arrays $\{ y^k_i \in \R_+ | 1 \leq k \leq n, \ 1 \leq i \leq r_k \}$ of nonnegative real numbers satisfying 
$$
y^{n+1}_{k+1} \leq  y^n_k \leq y^{n+1}_k, \ \ \ \ 1 \leq k \leq r_n.
$$
Consider the collection $B(t):=\{ B^k_i(t) | 1 \leq k \leq n, \ 1 \leq i \leq r_k \}$, $t \geq 0$, of independent Brownian motions started at $0$. Construct now a process $(\mc{W}(t))_{t \geq 0}$ in $\mb{J}^c_{n,paths}$ with components $\mc{W}^k_i$ driven by $(B(t))_{t \geq 0}$ by setting $\mc{W}^1_1 := \overline{\Gamma}(0,\infty| B^1_1)$ and inductively
\begin{equation} \label{continuoussk}
\mc{W}^k_i := \overline{\Gamma}(\mc{W}^{k-1}_{i},\mc{W}^{k-1}_{i-1}| B^k_{i}),
\end{equation}
for $2 \leq k \leq n$, $1 \leq i \leq r_{k}$ (recall our notational convention). This procedure determines another \emph{deterministic} map $\Psi^{SK}(B) : = \mc{W}$. Notice $\Phi^{SK}_{N} \to \Psi^{SK}$ as $N \to \infty$.

\bp \label{diffusionlimit}
Fix $n \geq 1$ and write $d_n := r_n(r_n+1)/2 + r_{n-1}(r_{n-1}+1)/2$. Then considered up to level $n$, the diffusively scaled symplectic Plancherel growth process $(\mc{X}(t;N))_{t \geq 0}$  converges in law on $D([0,\infty), \R^{d_n})$ to the interlaced reflecting system $(\mc{W}(t))_{t \geq 0}$ of \eqref{continuoussk} with a wall started at $\mbf{0} \in \mb{J}^c_{n,paths}$. 
\ep 
\bpf
For $n=1$, the laws of $\bar{X}^1_1$ converge on $D([0,\infty), \R)$ to the law of a standard Brownian motion (see, e.g., the much stronger Theorem 16.14 of \cite{k1}, where convergence in the Skorokhod topology strengthens to convergence in the supremem norm topology since the limit is continuous). In turn, the law of the images $\mc{X}(t;N)=\bar{\Gamma}(0,\infty| \bar{X}^1_1(\cdot, N))(t)$ converge to a standard reflected Brownian motion $\bar{\Gamma}(0,\infty|B^1_1)(t)=B^1_1(t) + \sup_{s \in [0,t]} [ -B^1_1(s)]^+$.

For $n \geq 2$, assume that the statement of the theorem holds for all levels $k$, $1 \leq k \leq n-1$, so that the restriction of $\mc{X}(t;N)$ to these levels converges in law on $D([0,\infty), \R^{d_{n-1}})$. As $N \to \infty$, the rescaled $\bar{X}^n_i(\cdot;N)$ converge in law on $D([0,\infty), \R)$, in the supremum norm topology, to $B_i$, $1 \leq i \leq r_n$, where the $B_i$ are independent standard Brownian motions started at $0$. Independence ensures the component-wise topology of uniform convergence on compacts is respected. Since the Skorokhod map is jointly continuous in this topology (cf. again Theorem 2.6 of \cite{bkr1}), the continuous mapping theorem in the form of Theorem 4.27 of \cite{k1} implies that as $N \to \infty$ the laws of $(\mc{X}^n_i(t;N))_{t \geq 0}$, $1 \leq i \leq r_n$, converge on $D([0,\infty), \R)$ to $(\mc{W}^n_i(t))_{t \geq 0}$, which completes the proof.

\epf

\subsection{Dyson non-colliding dynamics of limiting system levels} \label{sectiondynamics}

To determine the dynamics of the limiting system $(\mc{W}(t))_{t \geq 0}$ constructed in \eqref{continuoussk}, we follow Warren's approach \cite{jw1} of working with two levels at a time. Keep $n \geq 1$ fixed and fix $k$, $1 \leq k < n$. Define
\begin{equation} \label{chambercd}
W^k : = \{ u \in \R^{r_k}_+ : u_1 > \cdots > u_{r_k} \geq_k 0 \}, \ \ \  \geq_k:= 
\begin{cases} 
>, & \text{if} \ \ k \ \ \text{even} \ (\text{Type C}) \\
\geq, & \text{if} \ \ k \ \ \text{odd} \ (\text{Type D}) \\
\end{cases}
\end{equation}
For $u \in W^k$ and $v \in W^{k+1}$, let $u \prec v$ denote \emph{interlacement}:
$$
v_1 \geq u_1 \geq \cdots \geq v_{r_k} \geq u_{r_k} \geq v_{r_k+1},
$$
and additionally define
\begin{equation} \label{doublechamber}
W^{k+1,k}: = \{ (v,u) \in W^{k+1} \times W^k :  u \prec v \},
\end{equation}
For fixed $(v,u) \in W^{k+1,k}$, consider a filtered probability space $\left (\Omega, \mc{F}, \{ \mc{F}_t \}_{t \geq 0}, Q^{k}_{(v,u)}\right )$ supporting a pair $(V_t, U_t)$ of $\R^{r_{k+1}}_+ \times \R^{r_k}_+$--valued interlaced processes, i.e., $U_t \prec V_t$ for $t \geq 0$, with the following dynamics under $Q^{k}_{(v,u)}$. Take $\beta_1, \ldots, \beta_{r_k}$ to be independent $\mc{F}_t$-Brownian motions started at $u_1 , \ldots, u_{r_k}$. Let $U_i(t) := |\beta_i(t \wedge \tau)|$ to be stopped at the collision time
\begin{equation}
\tau := \inf \{ t \geq 0 : |\beta_{i-1}(t)| = |\beta_{i}(t)|,  \ \text{for some}  \ i, \ 1 < i \leq r_{k+1} \}.
\end{equation} 
Then take $\gamma_1, \ldots, \gamma_{r_{k+1}}$ be independent $\mc{F}_t$-Brownian motions started at $v_1, \ldots, v_{r_{k+1}}$, and let the dynamics of $V$ be reflected on $U$ through the extended Skorokhod map by
$$
V_i(t) := \overline{\Gamma}(U_i, U_{i-1}| \gamma_i(\cdot \wedge \tau))(t), \ \ \ 1\leq  i \leq r_{k+1}.
$$
(note our notational convention applies throughout the above definitions). Since $U$ starts at strictly separated components $u \in W^k$, $V$ admits, for $t< \tau$, the decomposition 
\begin{equation} \label{blocked}
V_i(t)= 
\begin{cases}
 \gamma_i(t) + \frac{1}{2} L^0_{V_i}(t) - L^0_{(U_{i-1}-V_i)}(t), & i=r_{k+1} \ \text{with} \ k +1 \ \text{odd} \\
 \gamma_i(t)+  L^0_{(V_i - U_i)}(t) - L^0_{(U_{i-1}-V_i)}(t) , &  \text{otherwise for} \ 1 \leq i \leq r_{k+1}
\end{cases},
\end{equation}
(see Remark 4.2 of \cite{bkr1}; cf. also \cite{kr1}). These expressions allow us to derive the transition probabilities for $(V,U)$ on $W^{k+1,k}$ explicitly; write $\phi_t(x) := \frac{1}{\sqrt{2 \pi t}} e^{-x^2/(2t)}$ and $\Phi_t(x) := \int_{-\infty}^x \phi_t(u) du$. 
\bp \label{explicitdensity}
For $t>0$ and $(v,u), (v',u') \in W^{k+1,k}$, the transition density $$Q^k_{(v,u)}(V_t \in dv', U_t \in du', t < \tau) = q^k_t((v,u), (v',u')) dv' du'$$ for the process $(V,U)$ killed at time $\tau$ is given by
$$
q^k_t((v,u), (v',u')) : = \det
\begin{bmatrix}
A^k_t(v,v') & B^k_t(v,u') \\
C^k_t(u,v') & D^k_t(u,u')
\end{bmatrix},
$$
where
$$
\begin{cases}
(A^k_t(v,v'))_{ij} = \phi_t(v'_j-v_i) + (-1)^{k} \phi_t(v'_j+v_i), & 1 \leq i,j \leq r_{k+1} \\
(B^k_t(v,u'))_{ij} = \Phi_t(u'_j-v_i) + (-1)^{k} \Phi_t(u'_j + v_i) - 1_{(j+1 \leq i)} & 1 \leq i \leq r_{k+1} ,  \ 1 \leq j \leq r_k \\
(C^k_t(u,v'))_{ij} =  \phi_t^{(1)}(v'_j-u_i) + (-1)^{k-1} \phi_t^{(1)}(v'_j+u_i) & 1 \leq i \leq r_{k}, \ 1 \leq j \leq r_{k+1} \\
(D^k_t(u,u'))_{ij} = \phi_t(u'_j-u_i) + (-1)^{k-1} \phi_t(u'_j+u_i), & 1 \leq i,j \leq r_{k} 
\end{cases}.
$$
\ep
The presence of the indicator ``$1_{(j+1 \leq i)}$" in the definition of $B^k_t$ can in part be explained by the calculation required for \eqref{duality}.
\bpf
The proof follows the same line of reasoning as Proposition 2 of \cite{jw1}. Fix $w'=(v',u') \in W^{k+1,k}$ and define $G(t,(v,u)) = G_{(v',u')}(t,(v,u)) : = q^k_t((v,u), (v',u'))$. We need to check that $G$ satisfies the heat equation on $(0,\infty) \times \R^{r_{k+1}} \times \R^{r_k}$ and the appropriate boundary conditions. We focus on the latter; the former follows since all functions of $(v,u)$ involved in the expression above solve the heat equation.

In the reflecting wall case $k$ odd, we compute $\partial_{u_{r_k}} |_{u_{r_k} = 0} (D^k_t(u,u'))_{r_kl}=\partial_{u_{r_k}} |_{u_{r_k} = 0} (C^k_t(u,v'))_{r_kl} =  0$, for any $1 \leq l \leq r_k$, so that $\partial_{u_{r_k}} |_{u_{r_k} = 0} G(t,(v,u))= 0$. For the absorbing wall case $k$ even, we have $G(t,(v,u)) = 0$ when $u_{r_k} = 0$. For either parity of $k$, if $v_i=u_i$, $1 \leq i \leq r_k$, we have $\partial_{v_i}|_{v_i = u_i} A_t(v,v')_{ij} = -C_t(u,v')_{ij}$ and $\partial_{v_i}|_{v_i = u_i} B_t(v,u')_{ij} = -D_t(u,u')_{ij}$, so that $\partial_{u_{r_k}} |_{u_{r_k} = 0} G(t,(v,u)) = 0$. The case $v_{i+1} =u_{i}$ is the same. Now if $u_{i-1} = u_{i}$, $1 < i \leq r_k$, then the $u_i$ and $u_{i+1}$ rows of $[C^k_t(u,v'), D^k_t(u,u')]$ are equal, giving $G(t,(v,u)) = 0$. Finally, one readily checks the correct initial condition $\lim_{t \downarrow 0} G(t,(v,u)) = \Pi_{i=1}^{r_k} \delta_{(u_i'-u_i)}$.

Now fix $f: W^{k+1,k} \to \R$ bounded and continuous. Note that the function 
$$
F(t,w) : = \int_{W^{k+1,k}} q_t^k(w, w') f(w') dw'.
$$
solves the heat equation on $(0,\infty) \times W^{k+1,k}$ and satisfies the same boundary conditions as above. Hence, for $T,\epsilon>0$, Ito's Formula shows that the process $(F(T-t+\epsilon,(V_t,U_t)), \mc{F}_t )_{t \in [0,T]}$ is a bounded local martingale, and thus a true martingale. This property, the bounded convergence theorem, and the regularity of $q^k_t$ (cf. Lemma 1 of \cite{jw1}) together give
$$
\ba
\int_{W^{k+1,k}} q_T^k(w, w') f(w') dw' & = \lim_{\epsilon \downarrow 0} F(T+\epsilon,(v,u)) \\
& =  \lim_{\epsilon \downarrow 0}  \mb{E}_{Q^k_{(v,u)}} [ F(\epsilon, (V_T,U_T))1_{(T<\tau)} ]= \mb{E}_{Q^k_{(v,u)}} [ f(V_T,U_T) 1_{(T<\tau)} ]
\ea
$$
completing the proof.
\epf
Let $P^k_u$ be the law of the stopped process $U_t$ started at $u \in W^k$ with the components $U_i(t) = |\beta_i(t \wedge \tau)|$ from above. Proposition \ref{explicitdensity} implies that for $u' \in W^k$,
$$
P^k_u(U_t \in du' , t < \tau) = p_t^k(u,u') du' := \det[ \phi_t(u'_j-u_i) + (-1)^{k-1} \phi_t(u'_j+u_i) ]_{i,j=1}^{r_k} du'.
$$
Now define
\begin{equation} \label{hk}
h_k(u): = 
\begin{cases}
\prod_{1 \leq i < j \leq r_k} (u_i^2 - u_j^2)\cdot \prod_{i=1}^{r_k} u_i,  & \text{if} \ \ k \ \ \text{even} \\
 \prod_{1 \leq i < j \leq r_k} (u_i^2 - u_j^2), & \text{if} \ \ k \ \ \text{odd} 
\end{cases},\ \ \ \ 
\alpha_k:= 
\begin{cases}
1/2, & \text{if} \ \ k \ \ \text{even} \\
-1/2, & \text{if} \ \ k \ \ \text{odd}
\end{cases},
\end{equation}
Note the functions $h_k$, $k \geq 1$, are positive harmonic in $W^k$ with respect to the generator of $U_t$, and $u \in W^{k}$ has strictly separated components by definition. We may then define the Doob $h_k$--transform $P^{k,+}_u$ of $P^k_u$ (cf. \cite{dg1}) by
$$
P^{k,+}_u : = \frac{h_k(U(t))}{h_k(u)} \cdot P^k_u \ \ \ \text{on} \ \ \mc{F}_t, 
$$
with $h_k$--transformed density
\begin{equation} \label{dynamicssource}
P^{k,+}_u(U_t \in du', t < \tau) = p_t^{k,+}(u,u') du' := \frac{h_k(u')}{h_k(u)} p_t^k(u,u') du'.
\end{equation}
Under this measure, $\tau$ is infinite and expression \eqref{dynamicssource} implies, by Lemma 3 of \cite{mkht1} and Ito's formula for convex functions (Theorem 7.1.(v), \cite{ikss1}), that $U_t$ satisfies the non-colliding dynamics \eqref{dysoncd}, corresponding to the case $k$. We further use Lemma 4 of \cite{mkht1} (cf. also \cite{pbpbno1}) to arrive at the law $P^{k,+}_{0}$ of $U_t$ issued from the origin:
\begin{equation} \label{entrancelaw}
P^{k,+}_0(U_t \in du' ) = \mu_t^{k}(u') du' := \frac{t^{-r_k(r_k+\alpha_k)}}{C_k} \exp \left \{ - \frac{ | u' |^2}{2t}  \right \} \cdot  h_k(u')^2 du'.
\end{equation}

Now for $w=(v,u) \in W^{k+1,k}$, the marginal $u \in W^k$ has strictly separated components, so we may define a measure $Q^{k,+}_w$ to be the Doob $h_k$-transform of $Q^k_w$:
$$
Q^{k,+}_w : = \frac{h_k(U(t))}{h_k(u)} \cdot Q^k_w \ \ \ \text{on} \ \ \mc{F}_t, 
$$
with $h_k$--transformed density
$$
Q^{k,+}_w ((V_t, U_t) \in dw') =  q_t^{k,+}(w,w') dw' := \frac{h_k(u')}{h_k(u)} q_t^k(w,w') dw'
$$
for $w' = (v',u') \in W^{k+1,k}, \ t>0.$ Note for $v \in W^{k+1}$,
$$
\int_{\{u \in W^k: u \prec v\}} h_k(u) du = \frac{h_{k+1}(v)}{k!!}
$$
so that $\lambda^k(v,u) := k!! \frac{h_k(u)}{h_{k+1}(v)}$ is a transition kernel from $W^{k+1}$ down to $W^{k}$; in particular, $\lambda^k(v, \cdot)$ is a probability density on $\{u \in W^k: u \prec v\}$. This density together with the entrance laws $\{\mu_t^{k}\}_{k \geq 1}$ furnish an entrance law $\nu_t^k((v,u)) := \mu_t^{k+1}(v)\lambda^k(v,u)$ of a measure $Q^{k,+}_{(0,0)}$ for the process $(V, U)$ issued from the origin: 
$$
Q^{k,+}_{(0,0)}(U_t \in dw) = \nu_t^k(w) dw, \ \ \ \nu_{t+s}^k(w) = \int_{W^{k+1,k}} \nu_s^k(w')q_t^{k,+}(w',w) dw'.
$$
The proof that $\nu_t^k((v,u)) =\mu_t^{k+1}(v)\lambda^k(v,u)$ is an entrance law is a consequence of an \emph{intertwining relation} that holds between the families $\{p_t^{k,+}\}_{k \geq 1}$ and $\{q_t^{k,+}\}_{k \geq 1}$, $t>0$: 
\begin{equation} \label{intertwining}
\int_{\{u \in W^k: u \prec v\} } \lambda^k(v,u) q_t^{k,+}((v,u), (v',u')) du' = p_t^{k+1,+}(v,v') \lambda^k(v',u'), \  \  v \in W^{k+1}, \  (v',u') \in W^{k+1,k}.
\end{equation}
This intertwining readily follows by directly ``integrating $u$ out" in the expression for $q_t^k$ in Proposition \ref{explicitdensity}:
\begin{equation} \label{duality}
\int_{\{u \in W^k: u \prec v\} } q_t^k((v,u),(v',u')) du = p^{k+1}_t(v,v').
\end{equation}
We now use these properties to determine the projected dynamics of $V$.

\bp \label{marginal}
Under $Q^{k,+}_{(0,0)}$, the marginal $(V_t)_{t \geq 0}$ has distribution $P^{k+1, +}_0$, so the process $V_t$ satisfies the system \eqref{dysoncd}, corresponding to the case $(k+1)$.
\ep
\bpf
Take $0<t_1 < \cdots < t_m$. Then iteratively using the intertwining relation \eqref{intertwining}, compute
$$
\ba
Q^{k,+}_{(0,0)} &(V_{t_1} \in A_1, V_{t_2} \in A_2, \ldots , V_{t_m} \in A_m)  \\
& = \int_{A_1} dv_1 \int_{\{u_1 \prec v_1\}} du_1 \cdots \int_{A_m} dv_m \int_{\{u_m \prec v_m\}} du_m \cdot \nu^k_{t_1}(v_1,u_1) \prod_{i=2}^{m} q_{t_i -t_{i-1}}^{k,+}((v_{i-1}, u_{i-1}), (v_i, u_i)) \\
& =  \int_{A_1} dv_1  \cdots \int_{A_m} dv_m \cdot \mu^{k+1}_{t_1}(v_1) \prod_{i=2}^{m} p_{t_i -t_{i-1}}^{k+1,+}(v_{i-1}, v_i) \cdot  \int_{\{u_m \prec v_m\}} \lambda^k(v_m, u_m ) du_m. \\
\ea
$$
Noting $\lambda^k(v_m, \cdot)$ is a probability density on $\{u_m \in W^k: u_m \prec v_m\}$ completes the proof.
\epf

\subsection{A path property of Dyson gaps and non-degeneracy of coupled dynamics}  \label{sectionpathproperty}

Proposition \ref{marginal} does not explain how the process $U_t$ is involved in the non-colliding dynamics of $V_t$. In this section, we confirm that the non-degenerate constraining dynamics \eqref{blocked} coupling $V_t$ with $U_t$ survive sending the starting point to the boundary of the chamber; thus, under the given setup $\left (\Omega, \mc{F}, \{ \mc{F}_t \}_{t \geq 0}, Q^{k,+}_{(0,0)} \right )$, the coupled system $(V_t,U_t)$ is a genuine semimartingale.

Besides being realizable as a Doob transform, $U_t$ is also the eigenvalue process of a certain matrix-valued process $H_t^k$; see Section IV of \cite{az1} or Section III.B of \cite{mkht1}. In particular, the process inherits the Brownian scaling property $U_i(t) \overset{d}{=} \sqrt{t} U_i(1)$ and we may cite the Hoffman--Weilandt inequality (Lemma 2.1.19, \cite{rm}): 
\begin{equation} \label{hoffmanweilandt}
\sum_{i=1}^{r_k}(U_i(t)-U_i(s))^2 \leq \text{tr}(H_t^k- H_s^k)^2.
\end{equation}
\bp \label{holder}
Fix $k \geq 1$, $1 < i \leq r_{k+1}$, and $h<1$. Write $\delta^i_t:=U_{i-1}(t)-U_i(t)$. Then for $z < e^{-8r_k^2}$, 
$$
\mb{P} \left ( \inf_{1\leq t \leq h^{-1}} \delta^i_t \leq z \right ) \lesssim \frac{1}{(\log z)^2}.
$$
\ep
\bpf
For $x$ in the chamber $W^k$ of \eqref{chambercd}, write
\begin{equation} \label{functionsforlya}
u_{i,1}(x) := \sum_{j \neq i} \left (\frac{1}{x_i - x_j} +\frac{1}{x_i + x_j} \right )+ \frac{1_{(k \ even)}}{x_i} , \ \ \ \ u_{i,2}(x) := -\partial_{x_i} u_{i,1}(x) \geq 0.
\end{equation}
Let $\lambda(t)=(\lambda_1(t) \geq \ldots \geq \lambda_{r_k}(t))$ denote the unique strong solution of 
$$
d\lambda_i(t)= d\beta_i(t) + u_{i,1}(\lambda(t)) dt, \ \ \ 1 \leq i \leq r_k,
$$ 
started at the origin, and notice $\lambda$ agrees with the process $U$ under $P^{k,+}_{0}$, except that $U_{r_k} = |\lambda_{r_k}|$ for $k$ odd. Note the identities (see pg. 252 of \cite{rm})
$$
\sum_{i=1}^{r_k} u_{i,1}(x)x_i = 2\binom{r_k}{2} + r_k 1_{(k \ even)}, \ \ \ \sum_{i=1}^{r_k} [u_{i,2}(x) - u_{i,1}(x)^2] = 0.
$$
Consider the Lyapunov function
\begin{equation} \label{lyafunction}
f(x):= \frac{2}{r_k} \sum_{i=1}^{r_k} x_i^2 - \frac{2}{r_k^2} \left [ \sum_{1 \leq i < j \leq r_k} \left [ \log (x_i - x_j) + \log (x_i + x_j)\right ] + 1_{(k \ even)}  \sum_{i=1}^{r_k} \log x_i \right],
\end{equation}
and note the inequality (see pg. 251 of \cite{rm})
\begin{equation} \label{lyapunov}
\frac{1}{x_i \pm x_j} \leq e^{r_k^2(f(x) +8)} \ \ \ \text{for} \ \ \ 1 \leq i<j \leq r_{k+1},
\end{equation}
where our notational convention applies, i.e., $x_{r_{k+1}} \equiv 0$ if $k$ even. We can also readily calculate 
$$
\partial_{x_i} f(x) = \frac{2}{r_k} \left [ 2x_i - \frac{1}{r_k}u_{i,1}(x) \right], \ \ \ \partial_{x_i}^2 f(x) = \frac{2}{r_k} \left [ 2 + \frac{1}{r_k}u_{i,2}(x) \right].
$$
Putting these facts together with Ito's formula gives for $t \in [1,h^{-1}]$
\begin{equation} \label{submartingale}
f(\lambda(t)) = f(\lambda(1)) + 2N_k\cdot (t-1) - \frac{1}{r_k^2} \sum_{i=1}^{r_k} \int_1^t u_{i,2}(\lambda(s)) ds + \frac{2}{r_k} \sum_{i=1}^{r_k} \int_1^t \left [ 2\lambda_i(s) - \frac{1}{r_k} u_{i,1}(\lambda(s)) \right ] d\beta_i(s).
\end{equation}
where $N_k:=(2r_k-1)+2\cdot1_{(k \ even)}$. In particular, $Z_t:=f(\lambda(t)) + \frac{1}{r_k^2} \sum_{i=1}^{r_k} \int_1^t  u_{i,2}(\lambda(s)) ds$ is an $\mc{F}_t$-submartingale for $t \in [1,h^{-1}]$. Hence, using \eqref{lyapunov}, we have for $z<e^{-8r_k^2}$,
$$
\ba
\mb{P} \left ( \inf_{1\leq t \leq h^{-1}} \delta^i_t \leq z \right ) & \leq \mb{P} \left ( \sup_{1\leq t \leq h^{-1}} Z_t \geq \frac{1}{r_k^2} \log \left ( \frac{1}{z} \right) - 8 \right )  \leq  \frac{ \mb{E} \left [f(\lambda_1) + 2N_k\cdot(h^{-1}-1) +M(h^{-1}) \right ]^2}{\left (\frac{1}{r_k^2} \log z+ 8 \right )^2 }
\ea
$$
where $M(t)$ is the martingale term of \eqref{submartingale} and where we have used Doob's submartingale tail inequality with parameter $p=2$ conditional on $\mc{F}_1$. To see the expectation is finite, first note the term $\mb{E}f(\lambda_1)^2$ is finite by applying the inequality $\log (1/x) \leq a (1/x)^{1/a}$ for $a, x >0$: for example, we can estimate
$$
\ba
\mb{E} (\log \delta_t^i )^{2} & = \mb{E}[ (\log(1/\delta^i_t))^{2}1_{( \delta_t^i \leq 1)}] +  \mb{E}[ (\log (\delta_t^i))^{2}1_{( \delta_t^i > 1)}] \leq 4 \left ( \mb{E}\left [ \frac{1}{\delta_t^i}\right] + \mb{E}[\delta_t^i ] \right ) \leq C_{k}'< \infty,
\ea
$$
where the upper bound $C_{k}'$ is independent of $t$ since $[1,h^{-1}]$ is compact and bounded away from $0$. Moreover, by Ito's isometry and independence of the $\mc{F}_t$-Brownian motions $\beta_i$, we have
$$
\mb{E} M(h^{-1})^2=\frac{4}{r_k^2} \sum_{i=1}^{r_k} \int_1^{h^{-1}} \mb{E}  \left [ 2\lambda_i(s) - \frac{1}{r_k} u_{i,1}(\lambda(s)) \right ]^2 ds<\infty,
$$
which is finite by the form \eqref{entrancelaw} of the density and by definition \eqref{functionsforlya} of the function $u_{i,1}(x)$. This completes the proof.
\epf
\bn
To establish this proposition in the previous version of this paper, we first used the integrability of $\mb{E}[\delta^i_s\delta^i_t]^{-q}$, $q<3/2$, along with basic chaining (cf. \cite{mt1}) to prove that $(\log \delta_t^i)^2$ is H{\"o}lder continuous on $[1,h^{-1}]$ of parameter $\beta \in (0,\frac{1}{6})$, or when applied to the classical system \eqref{dyson1}, of parameter $\beta \in (0, \frac{1}{2}(1-\frac{2}{\kappa+1}))$. This earlier approach thus fails for $\kappa =1$, but our new method captures a special probabilistic structure of the Dyson gaps with the Lyapunov function \eqref{lyafunction} and (with some extra effort) is sharp enough to be adapted to the critical GOE case.
\en
\bpf[Proof of Theorem \ref{pathproperty}]
We proceed by first demonstrating the path property of the consecutive Dyson gaps $\delta^i_t:=U_{i-1}(t)-U_i(t)$ and then exploiting it to deduce non-degeneracy of the constrained dynamics \eqref{blocked} for $V$. Fix $1 < i \leq r_{k+1}$ for $k \geq 2$, and fix $C>0$. Then
$$
\ba
\mb{P} \left ( \liminf_{t \downarrow 0} \frac{\delta^i_t}{t^c}<C \right )  & = \mb{P} \left (  \inf_{h^n \leq t  \leq h^{n-1}} \frac{\delta^i_t}{t^c}<C, \ \ \text{i.o.} \ n\right ) =  \mb{P} \left (\limsup_{n \to \infty} \left \{ \inf_{h^n \leq t  \leq h^{n-1}} \frac{\delta^i_t}{t^c}<C \right \} \right ) \\
\ea
$$
It is sufficient to show the probabilities of these distinguished events are summable, so we have
$$
\mb{P} \left ( \inf_{h^n \leq t  \leq h^{n-1}} \frac{\delta^i_t}{t^c}<C \right ) = \mb{P} \left ( \inf_{1 \leq t  \leq h^{-1}} \frac{\delta^i_{ t}}{t^c}<C(h^n)^{c-1/2} \ \right ) \leq \mb{P} \left ( \inf_{1 \leq t  \leq h^{-1}} \delta^i_{t}<Ch^{-c} (h^{c-1/2})^{n} \ \right )
$$
where we have used Brownian scaling for each fixed $t$ in the infimum and implicitly the fact $\delta^i_t \geq 0$. By Proposition \ref{holder} and the fact $h^{c-1/2}<1$, we may invoke Borel-Cantelli to conclude the path property.

Now, to reduce notational burden, write $W:=V_i$, $B:=\gamma_i$, $\ell(t) := U_{i}(t), \ r(t):=U_{i-1}(t)$, and consider the decomposition $W= B+Y$. The non-colliding dynamics \eqref{dysoncd} satisfy $\inf_{t \in [t_1,t_2]} (r(t)-\ell(t))>0$ for $0 < t_1 < t_2$, so the variation $\mc{V}_{[t_1,t_2]}Y$ of $Y$ on any such interval $[t_1,t_2]$ is finite by Corollary 2.4 of \cite{bkr1}. It therefore suffices to check the same is true for some interval $[0, \eta)$, $\eta \in (0,1)$, at $0$. We have just shown that for any $c \in (1/2, 1)$, almost surely
\begin{equation}
\liminf_{t \downarrow 0} \frac{r(t) - \ell(t)}{t^c} = \infty.
\end{equation}
As a result, $f(t):=r(t) - \ell(t)$ satisfies $f(t) \geq g(t):=t^{c}$ for small $t \geq 0$. Following the \emph{parabolic box} approximation of Theorem 4.8 of \cite{bkr1}, let $\{ s_k \}_{k \in \Z_{<0}}$ enumerate, decreasing as $k \downarrow -\infty$, the family of points $\{2^{-p} + q \cdot g^{1/\alpha}(2^{-p}) \}$ indexed by $q = 0, \ldots , \lfloor 2^{-p}/g^{1/\alpha}(2^{-p}) \rfloor -1$ and $p \geq 1$, where $\alpha \in (0,1/2)$ is to be determined. It is straightforward to see that for $s_k=2^{-p_k} + q_k g^{1/\alpha}(2^{-p_k})$,
$$
g^{1/\alpha}(2^{-p_k}) \leq s_{k+1} - s_{k} \leq 2 g^{1/\alpha}(2^{-p_k}), 
$$
To estimate the variation of $Y$ on the intervals $[s_{k}, s_{k+1}]$, define 
$$
b_k: = \inf_{t \in [s_{k}, s_{k+1}]} r(t), \ \ \ a_k: = \sup_{t \in [s_{k}, s_{k+1}]} \ell(t), \ \ \ m_k:= (a_k+b_k)/2.
$$
Note that $r, \ell$, and thus also $f$, are almost surely $\beta$--H{\"o}lder continuous for any $\beta \in (0,1/2)$ by, say, the Hoffman--Weilandt inequality \eqref{hoffmanweilandt}. Let $c_\beta$ be the H{\"o}lder constant for $r, \ell$. Assuming we choose $\beta>\alpha$, there exists some random $k_{\ell, r, \beta}$ so that for $k<k_{\ell, r, \beta}$, the path property $f(s_k) \geq g(s_k) = (s_k)^c$ holds along with the estimate $ 1/2^\alpha - 2c_\beta | s_k - s_{k+1}|^{\beta-\alpha} \geq c_{\alpha, \beta}>0$, for some $c_{\alpha, \beta}$. Then for such $k<k_{\ell, r, \beta}$
\begin{equation} \label{parabolicestimate}
\ba
2c_\beta |s_k|^\beta \geq f(s_k) \geq b_k - a_k & \geq f(s_k) -  \sup_{t \in [s_{k}, s_{k+1}]} | r(s_k)-r(t)| - \sup_{t \in [s_{k}, s_{k+1}]}| \ell(t) - \ell(s_k)| \\
& \geq g(s_k) - 2c_\beta | s_k - s_{k+1}|^\beta \\
& \geq g(2^{-p_k}) - 2c_\beta | s_k - s_{k+1}|^\beta \\
& \geq \frac{1}{2^\alpha}|s_k - s_{k+1}|^\alpha - 2c_\beta | s_k - s_{k+1}|^\beta \\
& \geq c_{\alpha, \beta}|s_k - s_{k+1}|^\alpha
\ea
\end{equation}
Now write $\tau_1 := s_k$ and for $j \geq 1$ define stopping times 
$$
\sigma_j: = \inf \{ t \in [\tau_j, s_{k+1}] : W(t) = m_k \}, \ \ \ \tau_{j+1}: = \inf \{ t \in [\sigma_j , s_{k+1}] : |W(t) -m_k| \geq (b_k-a_k)/4\}.
$$
Let  $\mb{P}^{\ell, r}, \mb{E}^{\ell, r}$ denote the probability and expectation conditional on $\ell, r$ and note $B(t)$ is independent of $\ell, r$. Notice for all $j \geq 1$, $\mc{V}_{[\sigma_j,\tau_{j+1}]}Y = 0$ (by Definition \ref{ESP}) and we can estimate
$$
\ba
\mb{E}^{\ell, r} \left [ | W(\sigma_j) - W(\tau_j) | \right ] & \leq \mb{E}^{\ell, r} \left [ 1_{(\tau_j<s_{k+1})}1_{(\sigma_j=s_{k+1})}|W(\sigma_j) - m_k|  \right ]  + \mb{E}^{\ell, r} \left [1_{(\tau_j<s_{k+1})} | m_k - W(\tau_j) | \right ] \\
& \leq \mb{P}^{\ell, r}(\tau_j<s_{k+1}) [ f(s_{k+1}) + f(s_k) ].
\ea
$$
Hence by Doob's $L^2$--submartingale inequality conditional on $\mc{F}_{\tau_j}$, for all $j \geq 1$,
$$
\ba
\mb{E}^{\ell, r} \left [ | Y(\sigma_j) - Y(\tau_j) | \right ] & \leq \mb{E}^{\ell, r} \left [ | W(\sigma_j) - W(\tau_j) | \right ] + \mb{E}^{\ell, r} \left [ \sup_{t \in [\tau_j,\sigma_j]} | B(t) - B(\tau_j) | 1_{(\tau_j<s_{k+1})} \right ] \\
& \leq [f(s_{k+1}) + f(s_k) + 2(s_{k+1} - s_k)^{1/2}] \cdot \mb{P}^{\ell, r}(\tau_j<s_{k+1}).
\ea
$$
Using now the Doob $L^{1/\alpha}$--submartingale tail bound, we expect over the period $[s_k,s_{k+1}]$ only finitely many Brownian oscillations greater than $(b_k-a_k)/4 \geq c_{\alpha, \beta}|s_k-s_{k+1}|^\alpha/4$ (here, we used estimate \eqref{parabolicestimate}):
$$
\ba
\mb{E}^{\ell, r} \left [ \sum_{j \geq 2} 1_{(\tau_j<s_{k+1})} \right ] = \sum_{j \geq 2} \mb{P}^{\ell, r}(\tau_j<s_{k+1}) & \leq \sum_{j \geq 2} \mb{P}^{\ell, r} \left (  \sup_{t \in [\sigma_{j-1}, \tau_j]} | B(t) - B(\sigma_{j-1}) |   \geq \frac{c_{\alpha, \beta}|s_k - s_{k+1}|^\alpha}{4}    \right )  \\
& \leq \sum_{j \geq 2}  \mb{E}^{\ell, r} \left [ 4 \frac{(\tau_j-\sigma_{j-1})^{1/2}}{c_{\alpha, \beta}|s_k - s_{k+1}|^\alpha}   \right ]^{1/\alpha} \\
& \leq  \frac{4^{1/\alpha}\mb{E}^{\ell, r} \left [\sum_{j \geq 2} (\tau_j-\sigma_{j-1}) \right ]}{(c_{\alpha, \beta})^{1/\alpha} |s_k - s_{k+1}|} \\
& \leq \frac{4^{1/\alpha}}{(c_{\alpha, \beta})^{1/\alpha}} = : C_{\ell, r, \alpha, \beta} < \infty, \ \ \ \ \text{a.s.}\ \text{for  $k<k_{\ell, r, \beta}$},
\ea
$$
where the third inequality uses the assumption $1/(2\alpha)>1$ with $\tau_j-\sigma_{j-1} \leq 1$. Since $Y$ is monotone on $[\tau_j,\sigma_j]$, $j \geq 1$ (by Definition \ref{ESP}), our work so far yields, for $k< k_{\ell, r, \beta}$,
$$
\ba
\mb{E}^{\ell, r} \mc{V}_{[s_k,s_{k+1}]}Y & \leq \sum_{j \geq 1} \mb{E}^{\ell, r} \left [ | Y(\sigma_j) - Y(\tau_j) | \right ]  \\
& \leq (C_{\ell, r, \alpha, \beta}+1)[f(s_{k+1})+f(s_k) +2 |s_k-s_{k+1}|^{1/2}] \\
& \leq (C_{\ell, r, \alpha, \beta}+1)[2c_\beta(s_{k+1}^\beta + s_k^\beta) +2^{1+1/2} (2^{-p_k})^{c/2\alpha} ] \\
& \leq (C_{\ell, r, \alpha, \beta}+1)[2\cdot 4^\beta c_\beta (2^{-p_k})^\beta+2^{1+1/2} (2^{-p_k})^c]\\
& \leq C_{\ell, r, \alpha, \beta}' (2^{-p_k})^\beta
\ea
$$
by the $\beta$--H\"{o}lder continuity of $f$ and the facts $f(0)=0$ and $c/(2\alpha)>c>\beta$. Summing this estimate over the relevant intervals $[s_{k'},s_{k'+1}]$ contained in $[2^{-p_k}, 2^{-p_k+2}]$ gives
$$
\sum_{k':[s_{k'},s_{k'+1}] \subset [2^{-p_k}, 2^{-p_k+2}]} \mb{E}^{\ell, r} \mc{V}_{[s_{k'},s_{k'+1}]}Y \leq C_{\ell, r, \alpha, \beta}'' \frac{2^{-p_k} (2^{-p_k})^\beta}{g^{1/\alpha}(2^{-p_k})} = C_{\ell, r, \alpha, \beta}'' (2^{-p_k})^{\beta+1 - c/\alpha}.
$$
Hence, as long as we take $\alpha, \beta, c$ so that $1/2<c < \alpha(\beta + 1), \alpha < \beta < 1/2$ (e.g., take $\alpha = 6/15, \beta=7/15$), summing over $p = p_k \geq 1$ yields a finite number, as required.

To summarize, by repeatedly applying Proposition \ref{marginal} and the conclusions of this proof, the pathwise constructed system $(\mc{W}(t))_{t \geq 0}$ of \eqref{continuoussk} is the unique strong solution to the non-degenerate interlaced system \eqref{continuumlimit1} of reflecting Brownian motions with a wall started at the origin. Further, for all $k \geq 1$, the $k$th level $\mc{W}^k$ is distributed as $P^{k,+}_0$, each pair $(\mc{W}^{k+1},\mc{W}^k)$ distributed as $Q^{k,+}_{(0,0)}$, and conditional on $\mc{W}^k$, the dynamics $\mc{W}^{k+1}$ are independent of $\mc{W}^1, \ldots, \mc{W}^{k-1}$ for $k \geq 2$. Convergence of symplectic Plancherel growth $\mc{X}(t,N)$ to this system is the content of Proposition \ref{diffusionlimit}, which completes the proof of Theorem \ref{pathproperty}.
\epf

\section{Representation theory of symplectic groups} \label{sectionrep}

\subsection{Notations and ancillary results} \label{sectionnotations}

For $n \geq 1$, we write $r_n:= \lfloor (n+1)/2 \rfloor$. Let $\mathbb{J}_n$ denote the set of \emph{partitions} $\lambda=( \lambda_1 \geq \ldots \geq \lambda_{r_n} \geq 0)$ of nonnegative integers of length at most $r_n$ (the \emph{length} of a partition $\lambda$ is the number of nonzero terms). For $\lambda \in \mathbb{J}_n, \mu \in \mathbb{J}_{n+1}$, write $ \lambda \prec \mu$ if $\lambda$ \emph{interlaces} $\mu$ in the sense
$$
\mu_1 \geq \lambda_1 \geq \mu_2 \geq \lambda_2 \geq \ldots  \geq \mu_{r_n} \geq \lambda_{r_n} \geq \mu_{r_n+1}.
$$ 
For $\lambda \in \mb{J}_n$, the transformation $\widetilde{\lambda}_i : = \lambda_i +r_n - i$ arises naturally and frequently in representation theoretic formulas. For the construction of Section \ref{sectioncentral}, define the collection $\mathbb{J}_{n,seq}$ of all finite sequences $u = (u^1, \cdots, u^n)$ of partitions $u^i \in \mb{J}_i$, $1 \leq i \leq n$, of length $n$. Let $\mb{J}_{n,paths} \subset \mb{J}_{n,seq}$ denote the subset of \emph{Gelfand-Tsetlin patterns}, i.e., finite interlaced sequences $u = (u^1 \prec \cdots \prec u^n)$. Let $\mb{J}_{\infty, seq} $, $\mb{J}_{\infty, paths} $ denote infinite versions of each of these sets. Recall our notational convention:

\begin{convention} 
Quantities with indices that overflow (e.g. $\lambda_{r_n +1}$ for $\lambda \in \mb{J}_n$) are set to $0$ and those with indices that underflow (e.g. $\lambda_{0}$ for $\lambda \in \mb{J}_n$) are set to $\infty$. 
\end{convention}

For each $r = r_n$, consistently fix a basis $\{e_{-1}, e_{1}, \ldots, e_{-r}, e_{r} \}$ for $\C^{2r}$ and define an inner product by $[x,y]:= \sum_{k=1}^r x_k y_{-k} - x_{-k} y_{k}$ for $x,y \in \C^{2r}$. The even symplectic group $Sp(2r,\C)$ is the set of linear transformations $T$ preserving this form: $[Tx,Ty] = [x,y]$ for all $x,y \in \C^{2n}$. The odd symplectic group $Sp(2r-1,\C)$ is defined as the closed subgroup of transformations stabilizing $e_r$, and we realize $Sp(2r-2,\C)$ as the further subgroup stabilizing the pair $e_{-r},e_{r}$. This furnishes a fixed choice of embeddings $Sp(2r-2,\C) \subset Sp(2r-1,\C) \subset Sp(2r,\C)$. Proposition 1.1 of \cite{vs1} shows $Sp(2r-1,\C)$ admits the (Jordan) semidirect product decomposition $H_{2r-1}(\C) \rtimes Sp(2r-2,\C)$ into unipotent and semisimple parts, where $H_{2r-1}(\C)$ is the complex Heisenberg group. Define also the even compact symplectic group $Sp(2r): \cong U(2r) \cap Sp(2r,\C)$ and define $Sp(2r-1)$ to be the subgroup of $Sp(2r)$ stabilizing $e_r$.

\begin{remark}
An intrinsic, though different, definition for odd symplectic group is the set of transformations preserving a skew-symmetric bilinear form of \emph{maximal possible rank}; see section 9 of Proctor \cite{rp1} for a discussion of candidate definitions. We have adopted the definition of Gelfand-Zelevinsky \cite{gz1} to exploit the above inclusions and to account for intermediate levels of the Gelfand-Tsetlin scheme (see \cite{vs1} and Section \ref{sectioncharacters}).
\end{remark}

If $n=2r$, the finite-dimensional irreducible representations $V_\lambda$ of $Sp(n,\mathbb{C})$ are parametrized by $\lambda \in \mb{J}_{2r}$. In the less common odd case $n=2r-1$, the same partitions $\lambda \in \mathbb{J}_{2r-1} = \mathbb{J}_{2r}$ now parametrize certain \emph{indecomposable} $Sp(2r-1,\C)$--modules $\mc{L}(\lambda)$ that exhaust the nonisomorphic representations of $Sp(2r-1,\C)$ (see \cite{vs1}). For either parity of $n$ and for $\lambda \in \mb{J}_n$, denote the corresponding character by $\chi^\lambda_{n}$, and write $\text{dim}_{2r}\lambda := \text{dim}_{Sp(2r,\C)} V_\lambda$, $\text{dim}_{2r-1}\lambda := \text{dim}_{Sp(2r-1,\C)} \mc{L}(\lambda)$. Section \ref{sectioncharacters} shows these characters are expressible in terms the Chebyshev polynomials of the second and third kind $J_{k,1/2}(x)$, $J_{k,-1/2}(x)$, $k \geq 0$, which are defined by the initial conditions $J_{0,1/2}=1, J_{1,1/2}(x) = 2x$ and $J_{0,-1/2}=1, J_{1,-1/2}(x) = 2x-1$, while both satisfy the same three term recurrence 
\begin{equation} \label{recurrence}
xp_k(x) = \frac{1}{2} p_{k+1}(x) + \frac{1}{2} p_{k-1}(x), \ \ \ k \geq 1.
\end{equation}
They also satisfy, for $z$ on the unit circle,
\begin{equation} \label{identity1}
J_{k,1/2}\left ( \frac{z + z^{-1}}{2} \right ) = \frac{z^{k+1} - z^{-(k+1)}}{z- z^{-1}}, \ \ \ J_{k,-1/2}\left ( \frac{z + z^{-1}}{2} \right ) = \frac{z^{k+1/2} + z^{-(k+1/2)}}{z^{1/2}+z^{-1/2}},
\end{equation}
Our choice of notation comes from the relationships
$$
J_{k,1/2}:=\frac{(2k+2)!!}{2 \cdot (2k+1)!!} \cdot \mbf{J}^{(1/2,1/2)}_k, \ \ \ J_{k,-1/2}:=\frac{(2k)!!}{ (2k-1)!!} \cdot \mbf{J}^{(-1/2,1/2)}_k,
$$ 
where $\mbf{J}^{(\alpha,\beta)}_k$ are the Jacobi polynomials of parameter $(\alpha,\beta)$ defined to be orthogonal with respect to the weight $w_{(\alpha,\beta)}(x):=(1-x)^\alpha (1+x)^\beta 1_{[-1,1]}(x)$. We also have $\mbf{J}^{(\alpha,\beta)}_k(1) = \binom{k+\alpha}{k}$ for $k\geq1$ and $=1$ for $k=0$, so that $J_{k,-1/2}(1) = 1$ for $k \geq 0$, a fact we will need to exploit many times. Lastly, define an inner product with respect to a \emph{normalized} weight  
\begin{equation} \label{innerproduct}
\langle f, g \rangle_\alpha:= \frac{2^{\alpha+1/2}}{\pi}   \int_\R f(x)g(x) w_{(\alpha,1/2)}(x) dx.
\end{equation}
Using the relations in \eqref{identity1}, one can compute directly $\langle J_{l,\alpha}, J_{k,\alpha} \rangle_{\alpha}= \delta_{lk},$ which explains the choice of normalization in \eqref{innerproduct}. We will frequently make use of the orthogonal decomposition
\begin{equation} \label{orthogonaldecomposition}
T(x) = \sum_{k=0}^\infty \langle J_{k,\alpha_n}, T \rangle_{\alpha_n}  J_{k,\alpha_n}(x)
\end{equation} for $T \in C^1[-1,1]$ (cf. Lemma 2 of \cite{abjk1}). See Szego \cite{gs1} for more details on this discussion.

We now collect some results that are fundamental for computing the correlation kernel in Section \ref{sectiondeterminantal} and for uncovering the intertwining relationship in Section \ref{sectionMarkov}. 

\bl  \label{identities}
The Chebyshev polynomials satisfy the following identities
\begin{enumerate}

\item

\subitem $\sum_{r=0}^{s} J_{r,-1/2} = J_{s,1/2}$, for all $s \geq 0$. 

\subitem $\sum_{r=0}^{s-1} 2 J_{r,1/2} = \frac{J_{s,-1/2} - 1}{x-1}$, for all $s \geq 1$.

\item For any $T \in C^1[-1,1]$, 

\subitem $\sum_{r=s+1}^\infty \langle J_{r,-1/2}, T \rangle_{-1/2} =  \langle J_{s,1/2}, T(1)-T \rangle_{-1/2}$

\subitem  $\sum_{r=s}^\infty \langle J_{r,1/2}, T \rangle_{1/2} = \langle J_{s,-1/2}, T \rangle_{-1/2}$

\end{enumerate}
\el

\bpf

For $x \in [-1,1]$, let $x = (z+z^{-1})/2$ for some $z$ on the unit circle. The first part then follows from using the identities in \eqref{identity1} and exploiting the resulting finite geometric sums. For the second part, the orthogonality relations and the first part together give
$$
\langle J_{s,1/2}, 1 \rangle_{-1/2} = \langle J_{s,1/2}, J_{0,-1/2} \rangle_{-1/2} =  \left \langle \sum_{r=0}^{s} J_{r,-1/2}, J_{0,-1/2} \right \rangle_{-1/2} =\langle J_{0,-1/2}, J_{0,-1/2} \rangle_{-1/2}  = 1.
$$
Orthogonal decomposition \eqref{orthogonaldecomposition} gives
$$
\langle J_{s,1/2}, T(1) \rangle_{-1/2} = T(1) = \sum_{r =0}^\infty \langle J_{r, -1/2}, T \rangle_{-1/2} J_{r,-1/2}(1) = \sum_{r =0}^\infty \langle J_{r, -1/2}, T \rangle_{-1/2},
$$
where we have used the fact $J_{r,-1/2}(1) = 1$. Then subtract $\langle J_{s,1/2}, T \rangle_{-1/2} = \sum_{r =0}^{s} \langle J_{r, -1/2}, T \rangle_{-1/2}$ from both sides to prove the first identity of the second part. For the other identity, note
$$
\sum_{r=0}^{s-1} \langle J_{r,1/2}, T \rangle_{1/2} = \frac{1}{2} \left \langle \frac{J_{s,-1/2} - 1}{x-1}, T  \right\rangle_{1/2} \overset{s \to \infty}{\rightarrow}  - \frac{1}{2}  \langle (x-1)^{-1}, T \rangle_{1/2} = \langle 1 , T \rangle_{-1/2}
$$
where the term disappears necessarily from convergence of orthogonal decompositions \eqref{orthogonaldecomposition} and the last equality uses the normalization of our inner products \eqref{innerproduct}. Subtracting the finite sum $\sum_{r=0}^{s-1} \langle J_{r,1/2}, T \rangle_{1/2} = \langle 1-J_{s,-1/2}, T \rangle_{-1/2}$ from the right side completes the proof.

\epf
Define
$$
\alpha_n:=
\begin{cases}
1/2, & \text{if} \ n \ \text{even} \\
-1/2, & \text{if} \ n \ \text{odd}
\end{cases}
 \ \ \ \ \phi_n(s,t):=
\begin{cases}
2\cdot 1_{(s<t)},&\text{if} \ n \ \text{even} \\
1_{(s \leq t)},& \text{if} \ n \ \text{odd}
\end{cases}.
$$ 
For smooth $E \in C^\infty[-1,1]$, denote the $m$th Taylor Remainder of $E$ about $1$ by
$$
R^E_m(x): = 
\begin{cases}
E(x), & \ \ m \leq 0 \\
E(x) - \sum_{k=0}^{m-1} \frac{E^{(k)}(1)}{k!}(x-1)^{k}, & \ \ m \geq 1 
\end{cases}.
$$
and let $\Psi^{n}_{r_n-l}(s): = 
\langle J_{s,\alpha_n},  (x-1)^{r_n-l} R^{E}_{l-r_n} \rangle_{\alpha_n}$.
 
\bp \label{comprule}
The functions $\Psi^{n}_{r_n-l}(s)$, $n \geq 1$, $l \in \Z$, satisfy the composition rule
$$
(\phi_{n-1}*\Psi^{n}_{r_n-l})(s):= \sum_{t\geq 0} \phi_{n-1}(s,t)\Psi^{n}_{r_n-l}(t) = \Psi^{n-1}_{r_{n-1}-l}(s).
$$
\ep
\bpf
Fix $n$ odd so $r_{n-1} = r_n-1$ and $\alpha_{n} = -1/2$. Then using part $(3)$ of Lemma \ref{identities} with $T_m(x) := (x-1)^{-m} R^E_{m}(x)$ and the fact $T_m(1) := \lim_{x \to 1} T_m(x) = E^{(m)}(1)/m!$, we get
$$
\ba
\sum_{t\geq 0} \phi_{n-1}(s,t)\Psi^{n}_{r_n-l}(t) & = 2\sum_{t=s+1}^\infty \langle J_{t,-1/2}, (x-1)^{r_n-l} R^E_{l-r_n} \rangle_{-1/2} \\
&= \left \langle J_{s,1/2}, \frac{T_{r_n-l}-T_{r_n-l}(1)}{x-1} \right \rangle_{1/2} \\
&  =  \left \langle J_{s,1/2}, \frac{R^E_{l-r_{n-1}}(x)}{(x-1)^{l-r_{n-1}}} \right \rangle_{1/2}  =\Psi^{n-1}_{r_{n-1}-l}(s),
\ea
$$
where the second equality uses our inner product normalization \eqref{innerproduct}.
The case $n$ even instead involves $r_n = r_{n-1}$, $\alpha_{n} = 1/2$, and part $(4)$ of Lemma \ref{identities}:
$$
\ba
\sum_{r\geq 0} \phi_{n-1}(s,r)\Psi^{n}_{r_n-l}(r)& =  \sum_{r=s}^\infty \langle J_{r,1/2}, (x-1)^{r_n-l} R^E_{l-r_n} \rangle_{1/2} \\
&= \langle J_{s,-1/2}, (x-1)^{r_{n-1}-l} R^E_{l-r_{n-1}} \rangle_{-1/2} = \Psi^{n-1}_{r_{n-1}-l}(s).
\ea
$$ 
\epf

We collect one last result that will help us express our probability measures on partition sequences in determinantal form.
\bl \label{indication}
Indication of interlacement between two partitions $\lambda \in \mb{J}_{n}, \mu \in \mb{J}_{n+1}$ takes the determinantal form
$$
1_{(\lambda \prec \mu)} = \det[ \phi_{n}(\widetilde{\lambda}_i, \widetilde{\mu}_j)]_{i,j=1}^{r_{n+1}} \times 
\begin{cases}
2^{-r_{n+1}}, & \text{if} \ \ n \ \ \text{even} \\
1, & \text{if} \ \ n \ \ \text{odd}
\end{cases},
$$
where $\widetilde{\lambda}_i : = \lambda_i +r_n - i$.
\el
\bpf
First suppose $\lambda \prec \mu$. If $n$ is even, then $r_n = r_{n+1} - 1$ while if $n$ is odd, $r_n = r_{n+1}$. Under the transformation $\lambda, \mu \to \widetilde{\lambda}, \widetilde{\mu}$, the interlacing condition becomes
$$
\begin{cases}
 \widetilde{\mu}_1 > \widetilde{\lambda}_1 \geq \widetilde{\mu}_2 > \widetilde{\lambda}_2 \geq \ldots \geq \widetilde{\mu}_{r_{n+1}} > \widetilde{\lambda}_{r_{n+1}} \equiv -1, & \text{if} \ \ n \ \ \text{even} \\
 \widetilde{\mu}_1 \geq \widetilde{\lambda}_1 \geq \widetilde{\mu}_2 \geq \widetilde{\lambda}_2 \geq \ldots \geq \widetilde{\mu}_{r_{n+1}} \geq \widetilde{\lambda}_{r_{n+1}} \geq 0, & \text{if} \ \ n \ \ \text{odd}
\end{cases}
$$
Since $\phi_n(s,t) = 2\cdot 1_{(s<t)}$ for $n$ even, the strict inequalities $\widetilde{\mu}_i > \widetilde{\lambda}_i$ imply $\phi_n(\widetilde{\lambda}_i, \widetilde{\mu}_j) = 2$ for $1 \leq j \leq i \leq r_{n+1}$, and similarly if $n$ odd. Hence, the $r_{n+1} \times r_{n+1}$ matrix $[\phi_{n}(\widetilde{\lambda}_i, \widetilde{\mu}_j)]_{i,j=1}^{r_{n+1}}$ is triangular filled with $2$'s if $n$ even and with $1$'s if $n$ odd, which proves the statement.

If $\lambda \nprec \mu$, then let $k^*$ be the largest index of $\lambda$ such that the interlacing condition fails, i.e., one of $\mu_{k^*+1}>\lambda_{k^*}$ or $\lambda_{k^*} > \mu_{k^*}$ holds. If $\lambda_{k^*} > \mu_{k^*}$, then $\lambda_i > \mu_{k^*}\geq \mu_j$ for $1 \geq i \geq k^* \geq j \geq r_{k+1}$, which implies that the block of $[ \phi_{n}(\widetilde{\lambda}_i, \widetilde{\mu}_j)]_{i,j=1}^{r_{n+1}}$ with bottom left corner $(k^*,k^*)$ is filled with $0$'s. If $k^*=r_{n+1}$, then the $k^*$ column is a zero vector and ensures the determinant is $0$. Otherwise, by maximality of $k^*$, the $k^*$ and $k^*+1$ columns are the same (the matrix is triangular), and so again the determinant is $0$. The case $\mu_{k^*+1}>\lambda_{k^*}$ is handled the same way.
\epf

\subsection{Branching rules and character formulas} \label{sectioncharacters}

Branching rules determine the decomposition of a representation when restricted to the action of a subgroup. But restricting an irreducible representation $V_\lambda$, $\lambda \in \mathbb{J}_{2r}$, of $Sp(2r,\C)$ to the action of $Sp(2r-1,\C)$ does not readily yield a decomposition, at least not into irreducible representations. Fortunately, the restricted $V_\lambda$ does decompose into a \emph{semi-direct} sum of the indecomposable $Sp(2r-1,\C)$-modules $\{\mc{L}(\mu) \}_{\mu \in \mathbb{J}_{2r-1}}$. More precisely, Shtepin \cite{vs1} establishes the following (multiplicity-free) branching rules: for any irreducible representation $V_\lambda$, $\lambda \in \mathbb{J}_{2r}$, of $Sp(2r,\C)$ and any indecomposable representation $\mc{L}(\lambda)$, $\lambda \in \mathbb{J}_{2r-1}$, of $Sp(2r-1,\C)$, we have
\begin{equation} \label{branching}
V_\lambda |_{Sp(2r-1,\C)} = \boxplus_{\substack{\mu \in \mb{J}_{2r-1} \\ \mu \prec \lambda}} \mc{L}(\mu), \ \ \ \mc{L}(\lambda) |_{Sp(2r-2,\C)} = \oplus_{\substack{\mu \in \mb{J}_{2r-2} \\ \mu \prec \lambda}} V_\mu.
\end{equation}
where $\boxplus$ denotes semidirect sum.

Fix $\lambda \in \mb{J}_n$ and recall the notation $\widetilde{\lambda}_i:=\lambda_i + r_n - i$. Proposition 5.1 of \cite{vs1} shows that the indecomposable odd character $\chi^\lambda_{2r-1}$ depends only on its semisimple part $Sp(2r-2,\C)$. Using this fact with the branching rules \eqref{branching} gives
\begin{equation} \label{oddcharacter}
\chi^\lambda_{2r-1} = \chi^\lambda_{2r-1}|_{Sp(2r-2,\C)} = \sum_{\substack{\mu \in \mb{J}_{2r-2} \\ \mu \prec \lambda}} \chi^\mu_{2r-2}.
\end{equation}
Note also the determinantal identities 
\begin{equation} \label{detidentity}
\ba
\det [ z_j^{i-1+1} - z_j^{-(i-1+1)} ]_{i,j=1}^{r}  &= \det[ (z_j+z_j^{-1})^{i-1}]_{i,j = 1}^r
 \cdot \prod_{k=1}^r (z_k - z_k^{-1}) \\
\det [ z_j^{i-1+1/2} + z_j^{-(i-1+1/2)}]_{i,j=1}^{r} &= \det[ (z_j+z_j^{-1})^{i-1}]_{i,j = 1}^r \cdot \prod_{k=1}^r (z_k^{1/2} + z_k^{-1/2}).
\ea
\end{equation}
Combining \eqref{identity1}, \eqref{oddcharacter}, \eqref{detidentity} with the Weyl character formula for the even case (see Chapter 24, \cite{fh1}), we can compute for either parity of $n$ and for $M \in Sp(n)$
\begin{equation} \label{characterformula}
\chi^\lambda_{n}(M)=\chi^\lambda_{n}(z_1, \ldots, z_r) = \frac{ \det  \left [ J_{\widetilde{\lambda}_i,\alpha_n} \left ( \frac{z_j + z_j^{-1}}{2} \right ) \right ]_{i,j=1}^{r}}{\det [ (z_j + z_j^{-1})^{r-i}]_{i,j=1}^{r}},
\end{equation}
where the spectrum $\{z_1,z_1^{-1}, \ldots, z_{r_n},z_{r_n}^{-1} \}$ of $M$ has $z_i$ on the unit circle and $z_{r_n} = z_{r_n}^{-1} = 1$ if $n$ odd. Finally, only for the sake of explicitness, the dimension formulas can be obtained either from the branching rules \eqref{branching} or from evaluating $\lim_{Sp(n) \ni g \to e} \chi^\lambda_n(g)$, where $e$ is the identity:
\begin{equation} \label{dimension}
\text{dim}_{2r}\lambda = \prod_{1\leq i<j \leq r} \frac{l_i^2 - l_j^2}{m_i^2 - m_j^2} \cdot \prod_{1 \leq i \leq r} \frac{l_i}{m_i}, \ \ \ \ 
\text{dim}_{2r-1} \lambda = \prod_{1 \leq i<j \leq r} \frac{l_i'^2 - l_j'^2}{m_i'^2 - m_j'^2},
\end{equation}
where $l_i : = \lambda_i + (r_n - i) +1$, $m_i :=(r_n - i) +1$, and $l_i':= l_i - 1/2, \ m'_i := m_i - 1/2$ for $1 \leq i \leq r_n$.

\subsection{Plancherel measures for $Sp(\infty)$} \label{sectioninfinite}

For a general topological group $G$, a function $f:G \to \C$ is a \emph{class function} if $f(hgh^{-1}) = f(g)$ for all $g,h \in G$, and is \emph{positive definite} if the matrix $[f(g_j^{-1}g_i)]_{i,j = 1}^n$ is positive definite for any $g_1, \ldots, g_n \in G$. Consider the inductive limit $Sp(\infty):= \cup_{r=1}^\infty Sp(2r)$. A general characterization problem of Okounkov-Olshanski (Theorem 5.2, \cite{aogo1}; see also \cite{rpb1}) shows that the extreme points of the convex set of positive definite class functions $\chi : Sp(\infty) \to \C$ normalized to have $\chi(e) = 1$ at the identity $e$ are parametrized by the collection of triples $(\alpha, \beta, \delta)$ with $\alpha, \beta$ nonnegative decreasing summable sequences satisfying $\sum_1^\infty (\alpha_i +\beta_i) \leq \delta$ and $\gamma:= \delta - \sum_1^\infty (\alpha_i +\beta_i) \geq 0$. By slight abuse, we express these triples with either parameter $\delta$ or $\gamma$. For $\omega := (\alpha, \beta, \gamma)$ and $M \in Sp(n)$ with spectrum $\{z_1,z_1^{-1}, \ldots, z_{r_n},z_{r_n}^{-1} \}$ (again, the $z_i$ are on the unit circle and $z_{r_n} = z_{r_n}^{-1} = 1$ if $n$ odd), we have
$$
\chi^\omega (M) = \chi^\omega (z_1, \ldots, z_{r_n}) = \prod_{j=1}^{r_n} E^{(\alpha,\beta,\gamma)} \left ( \frac{z_j + z_j^{-1}}{2} \right )
$$ 
where
\begin{equation} \label{restrictionfunctions}
E^{(\alpha,\beta,\gamma)} \left ( \frac{z + z^{-1}}{2} \right ) := e^{\frac{\gamma}{2}(z+z^{-1} - 2)} \prod_{i=1}^\infty \frac{ (1 + \frac{\beta_i}{2}(z-1))(1+\frac{\beta_i}{2}(z^{-1} -1))}{(1 - \frac{\alpha_i}{2}(z-1))(1-\frac{\alpha_i}{2}(z^{-1} -1))}
\end{equation}
(the odd case follows by observing the inclusion $Sp(2r-1) \subset Sp(2r)$). Since the irreducible characters $\left \{ \chi^\lambda_{2r} \right \}_{\lambda \in \mb{J}_{2r}}$ form a complete orthonormal basis for class functions on $Sp(2r)$, we may write the restriction $\chi^\omega |_{n}:=\chi^\omega |_{Sp(n)}$ as a convex combination
\begin{equation} \label{evenrestriction}
\chi^\omega |_{2r} = \sum_{\lambda \in \mathbb{J}_{2r}} P^\omega_{2r}(\lambda) \frac{\chi^\lambda_{2r}}{\dim_{2r} \lambda}, 
\end{equation}
where convergence holds with respect to the inner product  $( \cdot, \cdot )_{2r}$ for characters. Evaluating both sides at $e$ yields $\sum_{\lambda \in \mathbb{J}_{2r}} P^\omega_{2r}(\lambda) = 1$ and positive definiteness of $\chi^\omega$ ensures that $P^\omega_{2r} \geq 0$. The $P^\omega_{2r}$ are thus probability measures on $\mathbb{J}_{2r}$. In the odd case, we cannot a priori rely on any classical theory for restrictions $\chi^\omega |_{2r-1}$ to $Sp(2r-1)$, so instead we use \eqref{evenrestriction} and the consequence \eqref{oddcharacter} of the branching rules to get
\begin{equation} \label{oddrestriction}
\ba
\chi^\omega |_{2r-1} & = \left ( \chi^\omega |_{2r} \right ) |_{2r-1} = \sum_{\lambda \in \mathbb{J}_{2r}} \frac{P^\omega_n(\lambda)}{\dim_{2r} \lambda} \chi^\lambda_{2r}|_{2r-1} = \sum_{\lambda \in \mathbb{J}_{2r}} \frac{P^\omega_n(\lambda)}{\dim_{2r} \lambda} \sum_{\mu \prec \lambda} \chi^\mu_{2r-1} \\
& = \sum_{\mu \in \mathbb{J}_{2r-1}}  \left (  \sum_{\lambda \succ \mu}\frac{P^\omega_{2r}(\lambda) \dim_{2r-1} \mu }{\dim_{2r} \lambda} \right ) \frac{ \chi^\mu_{2r-1}}{\dim_{2r-1} \mu} =: \sum_{\mu \in \mathbb{J}_{2r-1}} P^\omega_{2r-1}(\mu) \frac{\chi^\lambda_{2r-1}}{\dim_{2r-1} \mu}.
\ea
\end{equation}
The implicitly defined $P^\omega_{2r-1}$ are of course positive, and we will deduce unit mass on $\mb{J}_{2r-1}$ directly from their explicit expression, which we derive next. For each $\omega$, we call the series $P^\omega_n$, $n \geq 1$, the \emph{Plancherel measures for the infinite-dimensional symplectic group} $Sp(\infty)$ and we refer to the $P^\omega_{2r-1}$, $r \geq 1$, as \emph{intermediate}.
\bt \label{partitionmeasure}
Choose a triple $\omega = (\alpha,\beta,\gamma)$ to determine $E^\omega$ as in \eqref{restrictionfunctions}. The associated series of Plancherel measures $P_n^\omega$, $n \geq 1$, for $Sp(\infty)$ is explicitly given by
\begin{equation} \label{plancherelmeasures}
\ba
P_n^\omega(\lambda)= 2^{\binom{r_n}{2}} \cdot \det \left [  \langle J_{\widetilde{\lambda}_j,\alpha_n} , (x-1)^{r_n-i} E^\omega \rangle_{\alpha_n} \right]_{i,j=1}^{r_n} \dim_n \lambda, \ \ \ \ \lambda \in \mb{J}_n,
\ea
\end{equation}
where $\binom{1}{2} \equiv 0$ and we continue the notation $\widetilde{\lambda}_j: = \lambda_j + r_n-j$ for $\lambda \in \mb{J}_n$.
\et

\bpf
First assume $n$ is even. Note that $E^\omega \in C^1[-1,1]$ and write $x_i := (z_i + z_i^{-1})/2 \in [-1,1]$, where $z_1, \ldots, z_{r_n}$ provide the eigenvalues for $M \in Sp(n)$. By orthogonal decomposition \eqref{orthogonaldecomposition} and an infinite extension of the Cauchy-Binet Formula (cf. Lemma 1 of \cite{abjk1}), we have 
$$
\ba
\det \left [  (x_j-1)^{i-1}E^\omega(x_j) \right ]_{i,j=1}^{r_n} & = \det \left [  \sum_{k \geq 0}  \langle J_{k,\alpha_n} , (x-1)^{r_n-i} E^\omega \rangle_{\alpha_n} J_{k,\alpha_n}(x_j)\right ]_{i,j=1}^{r_n} \\
&= \sum_{k_1 > \ldots > k_r \geq 0} \det [ \langle J_{k_j,\alpha_n} , (x-1)^{r_n-i} E^\omega \rangle_{\alpha_n} ]_{i,j=1}^{r} \cdot \det[ J_{k_i,\alpha_n}(x_j) ]_{i,j=1}^{r_n}.
\ea
$$
Making the substitution $k_i =\widetilde{\mu_i}$, which takes a \emph{nonstrict} partition $\mu_1 \geq \ldots \geq \mu_{r_n} \geq 0$ to the strict one in the last summation above, and rearranging yields
$$
\ba
\chi^\omega |_{n}(M)=\prod_{k=1}^r E^\omega(x_k) & = \sum_{\mu \in \mb{J}_n} \det \left [  \langle J_{\widetilde{\mu}_j,\alpha_n} , (x-1)^{r_n-i} E^\omega \rangle_{\alpha_n} \right]_{i,j=1}^{r_n} \cdot \frac{ \det[J_{\widetilde{\mu}_i,\alpha_n}(x_j)]_{i,j=1}^{r}}{\det[(x_j-1)^{r-i}]_{i,j=1}^{r_n}} \\
& = \sum_{\mu \in \mb{J}_n} \det \left [  \langle J_{\widetilde{\mu}_j,\alpha_n} , (x-1)^{r_n-i} E^\omega \rangle_{\alpha_n} \right]_{i,j=1}^{r_n} \cdot \frac{ \det[J_{\widetilde{\mu}_i,\alpha_n}(x_j)]_{i,j=1}^{r_n}}{\det[x_j^{r-i}]_{i,j=1}^{r_n}}\\
&= 2^{\binom{r}{2}} \sum_{\mu \in \mb{J}_n} \det \left [ \langle J_{\widetilde{\mu}_j,\alpha_n} , (x-1)^{r_n-i} E^\omega \rangle_{\alpha_n} \right]_{i,j=1}^{r_n} \cdot \chi^\mu_{n}(z_1, \ldots, z_r).
\ea
$$
where we have used $\det[(x_j-1)^{r_n-i}]_{i,j=1}^{r_n} = \det[x_j^{r_n-i}]_{i,j=1}^{r_n}$. Applying the inner product $(\cdot, \chi^\lambda_n )_{n}$ for characters to both sides and recalling $P_n^\omega(\lambda):=
\left (\chi^\omega |_{n} , \chi^\lambda_{n} \right )_{n} \cdot \dim_{n} \lambda$ proves the result for $n$ even.

To compute $P^\omega_{n-1}(\mu)$ for $n-1$ odd, we use the formula just derived for $P^\omega_{n}(\lambda)$ to get (note $\alpha_{n-1} = -1/2$ and $r_n = r_{n-1}$) 
\begin{equation} \label{consistencycalc}
\ba
\frac{P^\omega_{n-1}(\mu)}{\dim_{n-1} \mu} := \sum_{\lambda \succ \mu}\frac{P^\omega_{n}(\lambda) }{\dim_{n} \lambda} & = \sum_{\lambda \in \mb{J}_{n}}  \det \left [  \langle J_{\widetilde{\lambda}_j,1/2} , (x-1)^{r_n-i} E^\omega \rangle_{1/2} \right]_{i,j=1}^{r_{n}} 1_{(\mu \prec \lambda)}\\
  & =   \sum_{\lambda \in \mb{J}_{n}}  \det \left [  \langle J_{\widetilde{\lambda}_j,1/2} , (x-1)^{r_n-i} E^\omega \rangle_{1/2} \right]_{i,j=1}^{r_{n}} \det \left [ 1_{( \widetilde{\mu}_j \leq \widetilde{\lambda}_i)} \right ]_{i,j=1}^{r_n} \\
  & = \det \left [ \sum_{k=\widetilde{\mu}_j}^\infty \langle J_{k,1/2} , (x-1)^{r_{n-1}-i} E^\omega \rangle_{1/2} \right]_{i,j=1}^{r_{n-1}} \\
  & = \det \left [ \langle J_{\widetilde{\mu}_j,-1/2} , (x-1)^{r_{n-1}-i} E^\omega \rangle_{-1/2} \right]_{i,j=1}^{r_{n-1}},
\ea
\end{equation}
where the third equality follows from our Lemma \ref{indication}, the fourth from the generalized Cauchy-Binet, and the last equality follows from part (4) of Lemma \ref{identities}. 

For either parity of $n$, recall $\dim_n \mu =  \lim_{Sp(n) \ni g \to e} \chi^\mu_n(g)$ and write $y_j(g) := (z_j(g) + z_j(g)^{-1})/2$ with the $z_j(g)$ giving the eigenvalues of the limiting $g \in Sp(n)$. Using \eqref{characterformula}, we then conclude unit mass by computing
\begin{align}  \label{unitmass}
\sum_{\mu \in \mb{J}_{n}} P^\omega_{n}(\mu) & = \sum_{\mu \in \mb{J}_{n}} 2^{\binom{r}{2}} \cdot \det \left [  \langle J_{\widetilde{\mu}_j,\alpha_n} , (x-1)^{r_n-i} E^\omega \rangle_{\alpha_n} \right]_{i,j=1}^{r_{n}}  \cdot \dim_{n} \mu \notag \\
& = \lim_{Sp(n) \ni g \to e}  \sum_{\mu \in \mb{J}_{n}} \frac{ \det \left [  \langle J_{\widetilde{\mu}_j,\alpha_n} , (x-1)^{r_n-i} E^\omega \rangle_{\alpha_n} \right]_{i,j=1}^{r_{n}}\cdot \det[J_{\widetilde{\mu}_i,\alpha_n}(y_j(g))]_{i,j=1}^{r_n} }{\det [ y^{r_n -i}_j(g) ]_{i,j=1}^{r_n}} \notag \\
& = \lim_{Sp(n) \ni g \to e}   \frac{ \det \left [ \sum_{k=0}^\infty \langle J_{k,\alpha_n} , (x-1)^{r_n-i} E^\omega \rangle_{\alpha_n} J_{k,\alpha_n}(y_j(g)) \right ]_{i,j=1}^{r_n} }{\det [ y^{r_n -i}_j(g) ]_{i,j=1}^{r_n}}  \\
& = \lim_{Sp(n) \ni g \to e}   \frac{ \det \left [ (y_j(g)-1)^{r_n-i} E^\omega(y_j(g)) \right ]_{i,j=1}^{r_n} }{\det [ (y_j(g)-1)^{r_n-i} ]_{i,j=1}^{r_n} } \notag \\
& = \lim_{Sp(n) \ni g \to e}  \prod_{j=1}^{r_n} E^\omega(y_j(g)) = 1, \notag
\end{align}
where Cauchy-Binet and orthogonal decomposition \eqref{orthogonaldecomposition} were used in the third and fourth equalities. This completes the proof.
\epf

\section{Determinantal correlation structure} \label{sectiondeterminantal}

\subsection{Central probability measures on partition paths} \label{sectioncentral}

The branching rules \eqref{branching} imply that interlaced sequences $u \in \mb{J}_{n,paths}$ ending in $u^n=\lambda$ exactly parametrize the basis elements of either $V_\lambda$ or $\mc{L}(\lambda)$. We thus have the elementary, yet crucial, relation
\begin{equation} \label{weightidentity}
\sum_{ \substack{u \in \mb{J}_{n,paths} \\ u^n = \lambda}} w_u = \dim_n \lambda.
\end{equation} 
where the weight $w_u :=  1_{(u^1 \prec u^2 \prec \cdots \prec u^n)}=\prod_{k=1}^{n-1} 1_{(u^k \prec u^{k+1})}$ indicates interlacement. For $u \in \mb{J}_{n,seq}$, consider the cylinder sets
$
C_u: = \{ t \in \mb{J}_{\infty, seq} | t^1 = u^1, \ldots, t^n = u^n \}.
$
We now realize the Plancherel measures $P^\omega_n$ on $\mathbb{J}_n$ of Theorem \ref{partitionmeasure} as embedded in a single probability measure $P^\omega$  on $\mb{J}_{\infty, seq}$ by the prescription
$$
P^\omega(C_u): = \frac{w_u P^\omega_{n}(u^n)}{\dim_n u^n}, \ \ \ u \in \mb{J}_{n, seq}.
$$
The weight $w_u$ ensures $P^\omega$ is supported on $\mathbb{J}_{\infty,paths} \subset \mb{J}_{\infty, seq}$. The \emph{consistency relations} 
\begin{equation} \label{consistency}
\frac{P^\omega_{n-1}(\mu)}{ \dim_{n-1} \mu } = \sum_{\substack{ \lambda \in \mb{J}_n \\ \lambda \succ \mu}}\frac{P^\omega_{n}(\lambda)}{\dim_{n} \lambda}
\end{equation}
established as in \eqref{consistencycalc} guarantee $P^\omega$ is well-defined. Lastly given $u,v \in \mb{J}_{n,seq}$, we obviously have the \emph{centrality condition} $P^\omega(C_u) w_v = P^\omega(C_v) w_u.$

The mapping
$$
\mb{J}_{\infty, seq} \ni u \mapsto \widetilde{\mc{P}}(u) := \{ (x_k^m,m)=(\widetilde{u}^m_k, m) \ | \ m \geq 1, \ 1 \leq k \leq r_m \}.
$$ 
pushes $P^\omega$ forward to a probability measure  $\widetilde{\xi}^\omega: =P^\omega \circ \widetilde{\mc{P}}^{-1}$ on $2^{\Z_{\geq 0} \times \Z_{> 0}}$, which determines a simple point process $\widetilde{\mc{X}}^\omega$ in $\Z_{\geq 0} \times \Z_{> 0}$, i.e., a random element of $2^{\Z_{\geq 0} \times \Z_{> 0}}$. We call this process simply the \emph{Plancherel point process}, though this title does not depend on the choice of coordinates. In these coordinates, the interlacing condition becomes 
$$
\begin{cases}
x^{n+1}_{k+1} \leq  x^n_k < x^{n+1}_k, \ \ \ \text{if} \ \ n \ \ \text{even} \\
x^{n+1}_{k+1} \leq  x^n_k \leq x^{n+1}_k, \ \ \ \text{if} \ \ n \ \ \text{odd}
\end{cases}, \ \ \ 1\leq k \leq r_n.
$$
and symplectic Plancherel growth can be visualized as
\begin{equation} \label{coordinates1}
 \ba
&| x^5_3 \ \ \ x^5_2 \ \  \  x^5_1 \\
&| x^4_2 \ \ \  x^4_1 \\
&| x^3_2 \ \ \  x^3_1\\
&|  x^2_1 \\
&| x^1_1 \to
\ea
 \ \ \  \longrightarrow \ \ \ \ 
 \ba
&| x^5_3 \ \ \ x^5_2 \to \    x^5_1 \\
&|  x^4_2 \ \ \  \ \ \ \ x^4_1 \\
&| x^3_2 \ \ \  \ \ \  \ x^3_1\\
&|  \ \ \ \ \ x^2_1 \\
&|  \ \ \  \ \ x^1_1
\ea
 \ \ \  \longrightarrow \ \ \ \ 
 \ba
&| \overset{\gets}{x^5_3} \ \ \  \ \ \  \ x^5_2 \  \   x^5_1 \\
&|  x^4_2 \ \ \  \ \ \ \ x^4_1 \\
&| x^3_2 \ \ \  \ \ \  \ x^3_1\\
&|  \ \ \ \ \ x^2_1 \\
&|  \ \ \  \ \ x^1_1
\ea
 \ \ \  \longrightarrow \ \ \ \ 
 \ba
&| x^5_3 \ \ \  \ \ \  \ x^5_2 \  \   x^5_1 \\
&|  x^4_2 \ \ \  \ \ \ \ x^4_1 \\
&| x^3_2 \ \ \  \ \ \  \ x^3_1\\
&|  \ \ \ \ \ x^2_1 \\
&|  \ \ \  \ \ x^1_1
\ea
\end{equation} 
Define the \emph{$k$th correlation function} for $z_1, \ldots, z_k \in \Z_{\geq 0} \times \Z_{> 0}$ by
$$
\widetilde{\rho}^\omega_k(z_1,\ldots,z_k): =\widetilde{\xi}^\omega(\{ E \in 2^{\Z_{\geq 0} \times \Z_{> 0}} | E \supset \{ z_1, \ldots, z_k \} \}) =  \mathbb{P}(\widetilde{\mc{X}}^\omega \supset \{z_1,\ldots, z_k \}).
$$
The coordinatization of Figure \ref{dynamicsteps} is similarly determined by the mapping
$$
\mb{J}_{\infty, seq} \ni u \mapsto \mc{P}(u) := \{ (y^m_k,m) = \left ( u^m_k + m-2k+1, m \right ) \ | \ m \geq 1, ,  \ 1 \leq k \leq r_m \},
$$ 
which as before determines a simple point process $\mc{X}^\omega$ in $\Z_{\geq 0} \times \Z_{> 0}$ via the probability measure $\xi^\omega: =P^\omega \circ \mc{P}^{-1}$ on $2^{\Z_{\geq 0} \times \Z_{> 0}}$. Its correlation functions $\rho_k^\omega$ are defined the same way.

\subsection{Computation of correlation kernel} 

Theorem \ref{kernel} of the introduction is a consequence of the following general statement. 

\bt \label{generalkernel}
Choose a triple $\omega = (\alpha,\beta,\gamma)$ as in Section \ref{sectioninfinite} so that there exists a loop in $\C$ around the interval $[-1,1]$ containing no zeros of the $E^\omega$ of \eqref{restrictionfunctions} (this is true if $\beta_1 <1$). Then the correlation functions $\{ \widetilde{\rho}^\omega_k\}_{k \geq 1}$ of the general Plancherel point process $\widetilde{\mc{X}}^\omega$ are determinantal: for $z_1, \ldots, z_k \in \Z_{\geq 0} \times \Z_{> 0}$,
$$
\widetilde{\rho}^\omega_k(z_1,\ldots,z_k) := \mathbb{P}(\widetilde{\mc{X}}^\omega \supset \{z_1,\ldots, z_k \})= \det[ K^\omega(z_i, z_j) ]_{i,j=1}^k,
$$
where the (nonsymmetric) kernel $K^\omega$ is given explicitly by
\begin{equation} \label{explicitkernel}
\ba
K^\omega((s,n),(t,m)) =  \frac{2^{\alpha_n+1/2}}{\pi} \frac{1}{2 \pi i} & \int_{-1}^1 \oint \frac{E^\omega(x)}{E^\omega(u)} J_{s,\alpha_n}(x) J_{t,\alpha_m}(u) \frac{(1-x)^{r_n + \alpha_n}(1+x)^{1/2}}{(1-u)^{r_m}(x-u)} dudx\\
& +1_{(n \geq m)} \frac{2^{\alpha_n+1/2}}{\pi}  \int_{-1}^1 J_{s,\alpha_n}(x) J_{t,\alpha_m}(x) (1-x)^{r_n-r_m+\alpha_n} (1+x)^{1/2} dx. 
\ea
\end{equation}
for $(s,n),(t,m) \in \Z_{\geq 0} \times \Z_{> 0}$, where the complex integral is a positively oriented (i.e., counterclockwise) simple loop around $[-1,1]$ containing no zeros of $E^\omega$.
\et

For any $n \geq 1$, a point configuration $\mb{X}_n := \{ x^m_k \in \Z | \ 1 \leq m \leq n, \  1 \leq k \leq r_m \}$ in $\Z_{\geq 0} \times [n]$ determines the cylinder set $C_u : = \widetilde{\mc{P}}^{-1}(\mb{X}_n)$ with $u \in \mb{J}_{n,seq}$ satisfying $\widetilde{u}^m_k = x^m_k$, $1 \leq m \leq n$. Using Lemma \ref{indication}, our push-forward measure $\widetilde{\xi}^\omega = P^\omega \circ \widetilde{\mc{P}}^{-1}$ takes the determinantal form
\begin{equation} \label{detmeasure}
\ba 
\widetilde{\xi}^\omega(\mb{X}_n)=P^\omega(C_u)& = w_u \frac{P^\omega_n(u^n)}{\dim u^n} = 2^{\binom{r_n}{2}} \prod_{m=1}^{n-1} 1_{(u^m \prec u^{m+1})} \cdot \det \left [ \langle J_{\widetilde{u^n_j},\alpha_n}, (x-1)^{r_n-i}E^\omega \rangle_{\alpha_n} \right]_{i,j=1}^{r_n} \\
& = \text{const} \left ( \prod_{m=1}^{n-1}  \det \left[ \phi_m(x^m_i, x^{m+1}_{j} ) \right]_{i,j=1}^{r_{m+1}} \right) \cdot \det \left [ \Psi^n_{r_n-j}(x^{n}_i) \right]_{i,j=1}^{r_n}
\ea
\end{equation}
where $\Psi^n_{r_n-i}(s) := \langle J_{s,\alpha_n}, (x-1)^{r_n-i} R^{E^\omega}_{i-r_n} \rangle_{\alpha_n}$ was defined in Proposition \ref{comprule}. Also, define convolution over $\Z_{\geq 0}$ by $(f*g)(x,y): = \sum_{z \geq 0} f(x,z) g(z,y)$ for bivariate functions $f,g$ and by $(f*g)(x):= \sum_{z \geq 0} f(x,z) g(z)$ if $g$ is univariate.
\bp \label{checkconditions}
For any $(s,n), (t,m) \in \Z_{\geq 0} \times \Z_{> 0}$ and $k \in \Z$, define the functions
$$
\ba
 \Phi^{m}_{r_m-k}(t) &:=  \frac{1}{2\pi i} \oint \frac{ J_{t,\alpha_m}(w) }{E(w)(w-1)^{r_m-k+1}} dw \\
\phi^{[n,m)}(s,t) &:=-\frac{1}{2 \pi i}  \oint \left \langle J_{s,\alpha_n}, \frac{J_{t,\alpha_m}(u)(u-1)^{r_n-r_m} }{x-u} \right \rangle_{\alpha_n} du, \ \ \  \text{for} \ \  n < m, 
\ea
$$
where the contours are positively oriented (i.e., counterclockwise) simple loops around $[-1,1]$ containing no zeroes of $E$. Then the simple point process $\widetilde{\mc{X}}^\omega$ determined by $\widetilde{\xi}^\omega$ of \eqref{detmeasure} has determinantal correlation functions $\widetilde{\rho}_k^\omega$ with kernel
\begin{equation} \label{prekernel}
K^\omega((s,n), (t, m)) = - \phi^{[n,m)}(s,t)1_{(n<m)} + \sum_{k=1}^{r_m} \Psi^{n}_{r_{n} - k}(s) \Phi^{m}_{r_{m} - k}(t).
\end{equation}
\ep
\bpf
The result is a restatement of Proposition \ref{algebraicargument}, so it suffices to identify the functions used in Appendix \ref{appendix} with the functions in the statement above. By orthogonal decomposition \eqref{orthogonaldecomposition}, we have for $1 \leq k,l \leq r_n$ 
$$
\ba
\sum_{s \geq 0} \Psi^{n}_{r_n - k}(s)  \Phi^{n}_{r_n-l}(s) &= \frac{1}{2 \pi i} \oint \sum_{s \geq 0} \langle J_{s,\alpha_n}, (x-1)^{r_n-k} E \rangle_{\alpha_n} J_{s,\alpha_n}(w) \frac{1}{E(w)(w-1)^{r_n-l+1}} dw \\
& = \frac{1}{2 \pi i} \oint \frac{1}{(w-1)^{k-l+1}} dw = \delta_{kl}.
\ea
$$
Note also that $\Phi^{m}_{r_m-k}(t)$ is a polynomial in $t$ of the same degree as $(\phi_{q_k-1} * \phi^{[q_k,m)})(-1,t)$, where $q_k = 2k-1$. These items confirm that $\{ \Phi^{m}_{r_m-k}(t) \}_{k=1}^{r_m}$ is the unique basis of the linear span of $\{ (\phi_{q_k-1} * \phi^{[q_k,m)})(-1,t) \}_{k=1}^{r_m}$ that is biorthogonal to the $\{\Psi^{n}_{r_n - k}(s)\}_{k=1}^{r_n}$.

We also need to show $\phi^{[n,m)}= \phi_{n} * \cdots * \phi_{m-1}$ for $n<m$. First assume $m = n+1$. If $n$ is odd, $r_n = r_{n+1}$, $\alpha_n = -1/2$, $\alpha_{n+1} = 1/2$, and $\phi_n(s,t) = 1_{(s \leq t)}$, so that
$$
\ba
-\frac{1}{2 \pi i}  \oint \left \langle J_{s,-1/2} , \frac{J_{t,1/2}(u)}{x-u} \right \rangle_{-1/2} du & = \langle J_{s,-1/2}, J_{t,1/2} \rangle_{-1/2} \\
& = \left \langle J_{s,-1/2}, \sum_{r=0}^t J_{r,-1/2} \right \rangle_{-1/2}  = \phi_n(s,t),
\ea
$$
where the first equality computes the residue at $u = x$, the second follows from (1) of Lemma \ref{identities}, and the last from orthogonality. Similarly, if $n$ is even, $\phi_n(s,t) = 2\cdot 1_{(s<t)}$, $r_n = r_{n+1}-1$, and the $\alpha_\cdot$'s switch, so that
$$
\ba
-\frac{1}{2 \pi i}  \oint \left \langle J_{s,1/2} , \frac{J_{t,-1/2}(u)}{x-u} \right \rangle_{1/2} \frac{1}{u-1} du & =  \left \langle J_{s,1/2}, \frac{J_{t,-1/2} - 1}{x-1} \right \rangle_{1/2} \\
& = 2 \left \langle J_{s,1/2}, \sum_{r=0}^{t-1} J_{r,1/2} \right \rangle_{1/2}  = \phi_n(s,t),
\ea
$$
where now the first equality computes residues at \emph{both} $u = x$ and $u = 1$, and uses the fact that $J_{t,-1/2}(1) = 1$. The full statement for general $n<m$ then follows by induction.
\epf

From our current \eqref{prekernel}, a few more calculations are required to arrive at our main expression \eqref{explicitkernel} for $K^\omega$ in Theorem \ref{kernel}. 
\bl For any $n,m \geq 1$,
\begin{equation} \label{nearlykernel}
\ba
\sum_{k=1}^{r_m} \Psi^{n}_{r_{n} - k}(s) \Phi^{m}_{r_{m} - k}(t)
 =  \frac{1}{2 \pi i} \oint &  \left \langle J_{s,\alpha_n}, \frac{(u-1)^{r_n}(R^E_{r_m-r_n}(u) - E(u)) +(x-1)^{r_n}E}{(x-u)E(u)(u-1)^{r_m}} \right \rangle_{\alpha_n} J_{t,\alpha_m}(u) du\\
& + \left \langle J_{s,\alpha_n}, \frac{ (x-1)^{r_n-r_m}R^E_{r_m-r_n}}{E} J_{t,\alpha_m} \right \rangle_{\alpha_n}. \\
\ea
\end{equation}
\el
Note how the terms in \eqref{nearlykernel} simplify if $r_n \geq r_m$, since in this case $R^E_{r_m-r_n} \equiv E$. 

\bpf
First note the identity
\begin{equation} \label{geometricidentity}
\sum_{k=1}^q \left ( \frac{u-1}{x-1} \right )^k = \left ( \frac{u-1}{x-u} \right ) \left ( 1- \left ( \frac{u-1}{x-1}\right )^q \right).
\end{equation}
Then if $r_m \leq r_n$, 
\begin{equation} \label{geometric}
\ba
\sum_{k=1}^{r_m} \Psi^{n}_{r_{n} - k}(s) \Phi^{m}_{r_{m} - k}(t) = \frac{1}{2 \pi i} \oint \frac{J_{s,\alpha_m}}{E(u)} \left \langle J_{s, \alpha_n}, \frac{(x-1)^{r_n}}{(u-1)^{r_m}} \frac{E(x)}{x-u} \left ( 1- \left ( \frac{u-1}{x-1} \right )^{r_m}  \right ) \right \rangle_{\alpha_n} du.
\ea
\end{equation}
and taking the residue at $u=x$ of the second summand in the inner product yields \eqref{nearlykernel}. If instead $r_m>r_n$, write
$$
\sum_{k=1}^{r_m} \Psi^{n}_{r_{n} - k}(s) \Phi^{m}_{r_{m} - k}(t) = \sum_{k=1}^{r_n} \Psi^{n}_{r_{n} - k}(s) \Phi^{m}_{r_{m} - k}(t) + \sum_{k=r_n+1}^{r_m} \Psi^{n}_{r_{n} - k}(s) \Phi^{m}_{r_{m} - k}(t).
$$
The first summand is treated as in \eqref{geometric}:
$$
\ba
\sum_{k=1}^{r_n} \Psi^{n}_{r_{n} - k}(s) \Phi^{m}_{r_{m} - k}(t) = \frac{1}{2 \pi i} \oint \frac{J_{s,\alpha_m}}{E(u)} \left \langle J_{s, \alpha_n}, \frac{(x-1)^{r_n}}{(u-1)^{r_m}} \frac{E(x)}{x-u} \left ( 1- \left ( \frac{u-1}{x-1} \right )^{r_n}  \right ) \right \rangle_{\alpha_n} du.
\ea
$$
For the second summand, we use the identity \eqref{geometricidentity} repeatedly to get
$$
\ba
 \sum_{k=r_n+1}^{r_m}  (E(x) -E(1)) & \frac{(u-1)^k}{(x-1)^{k-r_n}}  - 1_{(r_m \geq r_n+1)}  \sum_{k=r_n+2}^{r_m} \sum_{r=1}^{k-r_n-1} \frac{E^{(r)}(1)}{r!} \frac{(u-1)^k}{(x-1)^{k-r_n-r}}\\
 =  (E(x) &-E(1)) \frac{(u-1)^{r_n+1}}{x-u}  \left ( 1- \left ( \frac{u-1}{x-1}\right )^{r_m-r_n} \right) \\
  &- 1_{(r_m \geq r_n+1)} \sum_{r=1}^{r_m-r_n-1} \frac{E^{(r)}(1)}{r!} \frac{(u-1)^{r_n+r+1}}{x-u} \left ( 1- \left ( \frac{u-1}{x-1}\right )^{r_m-r_n-r} \right).
\ea
$$
and then arrive at
$$
\ba
\sum_{k=r_n+1}^{r_m} \Psi^{n}_{r_{n} - k}(s) \Phi^{m}_{r_{m} - k}(t) &= - \frac{1}{2 \pi i } \oint J_{t, \alpha_m}(u) \left \langle J_{s, \alpha_n}, \frac{(x-1)^{r_n-r_m} R^E_{r_m-r_n}}{E(u)(x-u)} \right \rangle_{\alpha_n} du  \\
& + \frac{1}{2 \pi i } \oint J_{t, \alpha_m}(u) \left \langle J_{s, \alpha_n}, \frac{(u-1)^{r_n-r_m} (R^E_{r_m-r_n}(u)-E(u)+E(x))}{E(u)(x-u)} \right \rangle du
\ea
$$
Taking the residue at $u=x$ of the first term and combining with the first summand completes the proof.
\epf
Now plug the expression \eqref{nearlykernel} into \eqref{prekernel}. If $r_n \geq r_m$, the expression \eqref{nearlykernel} simplifies, quickly leading to the main expression \eqref{explicitkernel} for $K^\omega$ in Theorem \ref{kernel} (note we have multiplied by the conjugating factor ``$(-1)^{r_n-r_m}$", which vanishes in the determinant). But if $r_n<r_m$, we get \eqref{explicitkernel} along with the additional term
$$
\ba
\frac{1}{2 \pi i} \oint & \left \langle J_{s,\alpha_n}, J_{t,\alpha_m}(u) \frac{(u-1)^{r_n-r_m}R^E_{r_m-r_n}(u)}{E(u)(x-u)} \right \rangle_{\alpha_n} du+ \left \langle J_{s,\alpha_n}, J_{t,\alpha_m} \frac{(x-1)^{r_n-r_m}R^E_{r_m-r_n}}{E} \right \rangle_{\alpha_n}
\ea
$$
The residue of the first term at $u=x$ exactly cancels the second, and since $$(u-1)^{r_n-r_m}R^E_{r_m-r_n}(u) = \sum_{k=r_m-r_n}^\infty \frac{E^{(k)}(1)}{k!}(u-1)^{k+r_n-r_m},$$ there is no residue at $u=1$ (recall by assumption $E=E^\omega$ does not have a pole at $1$ either). This completes the proof of Theorem \ref{generalkernel}.

\section{Multilevel Markov process} \label{sectionMarkov}

We state and prove the results of this section without coordinatizing $\mb{J}_{\infty,seq}$.

\subsection{Single level dynamics}

The point processes $\mc{X}^\omega$ arising from our general construction in Section \ref{sectiondeterminantal} are parametrized by the triples $\omega = (\alpha, \beta,\gamma)$ of Section \ref{sectioninfinite}. The next proposition dictates a natural way to transition between the Plancherel measures $P^\omega_n$ on a single level $\mb{J}_n$. Define $P^\psi_n$ by replacing $E^\omega$ with $\psi$ in the expression for $P^\omega_n$ of Theorem \ref{partitionmeasure}. For any fixed $\psi \in C^1[-1,1]$, define a $\mb{J}_n \times \mb{J}_n$--matrix $T^\psi_n$ with entries
$$
T^\psi_{n}(\mu, \lambda) := \det \left [ \left \langle  J_{\widetilde{\mu}_i,\alpha_n} , J_{\widetilde{\lambda}_j,\alpha_n}\psi \right \rangle_{\alpha_n} \right]_{i,j=1}^{r_n} \cdot \frac{\dim_n\lambda}{\dim_n \mu}.
$$ 

\begin{proposition} \label{Markov1}
Let $\psi, \psi_1, \psi_2 \in C^1[-1,1]$ be nonzero at $1$.
\begin{enumerate}

\item
$
\sum_{\mu \in \mathbb{J}_n} P^{\psi_1}_n(\mu) T^{\psi_2}_{n}(\mu,\lambda) = P^{\psi_1\cdot \psi_2}_n(\lambda),
$ 
and similarly $T^{\psi_1}_n T^{\psi_2}_n=T^{\psi_1\psi_2}_n$

\item $\sum_{\lambda \in \mb{J}_n} T^\psi_n(\mu,\lambda) = \psi(1)^{r_n}$

\item If $\psi_m \in C^1[-1,1]$ converge uniformly to $\psi$, then $T^{\psi_m}_{n}(\mu, \lambda) \to T^\psi_{n}(\mu, \lambda)$ as $m \to \infty$.

\item $T^\psi_n$ has nonnegative entries if $\psi(x) = p_0 + p_1x$ with $p_0>p_1 \geq 0$ or if $\psi_\gamma(x) = e^{\gamma(x-1)}$, $\gamma \geq 0$.

\item The (stochastic) semigroup $\{ T^{\psi_\gamma}_n \}_{ \gamma \geq 0}$ operating on the Banach space of absolutely summable functions $l^1(\mb{J}_n)$ on $\mb{J}_n$ is Feller: $ \Vert T^{\psi_\gamma}_n - Id \Vert_{l^1(\mb{J}_n)} \overset{\gamma \to 0}{\to} 0$.

\end{enumerate}
\end{proposition}
For example, set $\omega_0 = (0,0,0)$, and for $k \geq 1$, let $\omega_k$ be given by $\alpha = (q, q, \ldots, q, 0, \ldots )$ (i.e., $k$ q's), $\beta = (0,0, \ldots, )$, and $\gamma = 0$. Since $E^{\omega_k}\cdot E^{\omega_1} = E^{\omega_{k+1}}$ for $k \geq 0$, this choice along with the proposition yields a discrete-time Markov chain $X_t$, $t \geq 0$ with state space $\mb{J}_n$ and distribution $P^{\omega_k}_n$ at time $k$.
\bpf
Both items of the first point follow readily from the generalized Cauchy-Binet and orthogonal decomposition \eqref{orthogonaldecomposition}; for example,
$$
\ba
\sum_{\mu \in \mathbb{J}_n} & P^{\psi_1}_n(\mu) T^{\psi_2}_{n}(\mu,\lambda)  \\
&= 2^{\binom{r_n}{2}} \dim_n\lambda \sum_{\mu \in \mathbb{J}_n} \det \left [ \left \langle  J_{\widetilde{\mu}_j,\alpha_n} , (x-1)^{r_n - i} \psi_1 \right \rangle_{\alpha_n} \right]_{i,j=1}^{r_n} \cdot \det \left [ \left \langle  J_{\widetilde{\mu}_i,\alpha_n} , J_{\widetilde{\lambda}_j,\alpha_n}\psi_2 \right \rangle_{\alpha_n} \right]_{i,j=1}^{r_n} \\
& = 2^{\binom{r_n}{2}} \dim_n\lambda \cdot \det \left [ \left \langle \sum_{k=0}^\infty  \left \langle  J_{k,\alpha_n} , J_{\widetilde{\lambda}_j,\alpha_n}\psi_2 \right \rangle_{\alpha_n} , J_{k,\alpha_n} , (x-1)^{r_n - i} \psi_1 \right \rangle_{\alpha_n} \right]_{i,j=1}^{r_n}  \\
& = 2^{\binom{r_n}{2}} \dim_n \lambda \cdot \det \left [ \left \langle J_{\widetilde{\lambda}_j,\alpha_n}\psi_2 , (x-1)^{r_n - i} \psi_1 \right \rangle_{\alpha_n} \right]_{i,j=1}^{r_n} = P^{\psi_1\psi_2}_n(\lambda).
\ea
$$
The second point is similar to \eqref{unitmass}: use $\dim_n \lambda = \lim_{g \to e} \chi^\lambda_n(g)$ and \eqref{characterformula} to compute
\begin{align*}
\sum_{\lambda \in \mb{J}_{n}} T^\psi_n(\mu, \lambda) & = \frac{\lim_{Sp(n) \ni g \to e} }{\dim_n \mu} \sum_{\lambda \in \mb{J}_{n}}  \det \left [  \langle J_{\widetilde{\mu}_i,\alpha_n} , J_{\widetilde{\lambda}_j,\alpha_n} \psi \rangle_{\alpha_n} \right]_{i,j=1}^{r_{n}} \frac{\det[ J_{\widetilde{\lambda}_i}(y_j(g)) ]_{i,j=1}^{r_{n}}}{\det[(z_j(g)+z_j(g)^{-1})^{r_n-i}]_{i,j=1}^{r_n}}   \\
& =\frac{\lim_{Sp(n) \ni g \to e} }{\dim_n \mu} \frac{\det \left [ \sum_{k=0}^\infty \langle J_{\widetilde{\mu}_i,\alpha_n} \psi,  J_{k,\alpha_n} \rangle_{\alpha_n}J_{k}(y_j(g)) \right]_{i,j=1}^{r_{n}}}{\det[(z_j(g)+z_j(g)^{-1})^{r_n-i}]_{i,j=1}^{r_n}}  \\
& =\frac{\lim_{Sp(n) \ni g \to e} }{\dim_n \mu} \frac{\det \left [ J_{\widetilde{\mu}_i,\alpha_n}(y_j(g)) \right]_{i,j=1}^{r_{n}}}{\det[(z_j(g)+z_j(g)^{-1})^{r_n-i}]_{i,j=1}^{r_n}}\prod_{k=1}^{r_n} \psi(y_k(g)) = \psi(1)^{r_n}
 \end{align*}
where $z_j(g)$ are the eigenvalues of the limiting $g \in Sp(n)$, with $y_j(g) := (z_j(g) + z_j(g)^{-1})/2$, and where Cauchy-Binet and orthogonal decomposition \eqref{orthogonaldecomposition} were used in the second and third equalities. The third point of the proposition is a straightforward application of the dominated convergence theorem. 

Now for the fourth point, first take $\psi(x) = p_0 + p_1x$. Using the three term recurrence \eqref{recurrence}, orthogonality relations, and $xJ_{0,\alpha_n} = \frac{1}{2}[J_{1,\alpha_n}+J_{0,\alpha_n}(1-2\alpha_n)/2 ]$, we compute 
\begin{equation} \label{computetransitionentries}
 \langle  J_{\widetilde{\mu}_i,\alpha_n} , J_{\widetilde{\lambda}_j,\alpha_n}\psi  \rangle_{\alpha_n} = p_0 \delta_{\widetilde{\mu}_i,\widetilde{\lambda}_j} + \frac{p_1}{2} \left [ \delta_{\widetilde{\mu}_i,\widetilde{\lambda}_j+1} + \delta_{\widetilde{\mu}_i,\widetilde{\lambda}_j-1}1_{(\widetilde{\lambda}_j \geq 1)} + \frac{(1-2\alpha_n)}{2} \delta_{\widetilde{\mu}_i,0}\delta_{\widetilde{\lambda}_j,0}\right ].
\end{equation}
If $\widetilde{\lambda}_i > \widetilde{\mu}_i+1$ for some $1 \leq i \leq r_n$, then $\widetilde{\lambda}_k > \widetilde{\mu}_l+1$ for $1 \leq k \leq i \leq l \leq r_n$, which implies $T^\psi_{n}(\mu,\lambda) = 0$ (the resulting matrix admits a $2 \times 2$ block form with an off-diagonal block of $0$'s and a diagonal block with a zero vector). The same conclusion of course holds if $\widetilde{\mu}_i>\widetilde{\lambda}_i +1$, so assume $|\widetilde{\lambda}_i - \widetilde{\mu}_i|\leq 1$ for all $1 \leq i \leq r_n$. If $|\widetilde{\mu}_i - \widetilde{\mu}_{i+1}|>1$ for some $1 \leq i < r_n$, then $|\widetilde{\mu}_i - \widetilde{\lambda}_{i+1}| \leq 1$ implies $|\widetilde{\lambda}_{i}-\widetilde{\mu}_{i+1}| > 1$ and similarly $|\widetilde{\lambda}_{i}-\widetilde{\mu}_{i+1}| \leq 1$ implies $|\widetilde{\mu}_i - \widetilde{\lambda}_{i+1}| > 1$. In either case, $T^\psi_{n}(\mu,\lambda)$ breaks into a product of determinants, one of which is the $i \times i$ upper left corner minor of $[\langle  J_{\widetilde{\mu}_i,\alpha_n} , J_{\widetilde{\lambda}_j,\alpha_n}\psi  \rangle_{\alpha_n} ]_{i,j=1}^{r_n}$. Iterating this argument reduces consideration to the case where $|\widetilde{\mu}_i - \widetilde{\mu}_{i+1}|\leq1$ for all $1 \leq i < r_n$, which means $\mu_i$ are all equal to some $p \in \Z_{\geq 0}$. Now the two blocks corresponding to $\{i : \lambda_i=p\pm1\}$ are triangular with nonnegative entries \eqref{computetransitionentries} and straddle a tridiagonal block corresponding to $\{ i:\lambda_i=p = \mu_i \}$. Thus, writing $q^*:= \max\{ i:\lambda_i=p = \mu_i \}$, $q_*:= \min\{ i:\lambda_i=p = \mu_i \}$, $q:=q^*-q_*+1$, we have further reduced consideration to the determinant of the $q \times q$ tridiagonal matrix
\begin{equation} \label{tridiagonal1}
A_q(\lambda):=
\begin{bmatrix}
p_0 & p_1/2 & 0 & 0 & \cdots & 0  \\
 p_1/2& p_0  &  p_1/2  & 0 & \cdots & 0  \\
0 & p_1/2 & p_0  &  p_1/2 & \cdots & 0  \\
\vdots & \vdots  & \ddots  & \ddots & \vdots & 0  \\
0 & \cdots  & \cdots  & \cdots &  p_0 &  p_1/2 \\
0 & \cdots  & \cdots  & \cdots &   p_1/2& p_0 + p_{2,q^*}  \\
\end{bmatrix}.
\end{equation}
where $p_{2,q^*}:=\frac{p_1(1-2\alpha_n)}{4} \delta_{\widetilde{\lambda}_{q^*},0}$. Writing $r_{\pm} := (p_0 \pm \sqrt{p_0^2-p_1^2})/2$ and $h_n(x) := p_0x^n-\frac{p_1^2}{4}x^{n-1}$, we can compute, using the tridiagonal determinant recurrence relations and initial conditions,
\begin{equation} \label{lowerbound}
\det A_q(\lambda) = \frac{1}{\sqrt{p_0^2-p_1^2}} \cdot \left[ h_{q}(r_+) - h_{q}(r_-) + p_{2,q^*}(r_+^n-r_-^n)\right].
\end{equation}
Notice that $h_n(\frac{p_1^2}{4p_0}) = 0$ and $h$ is increasing on $[\frac{p_1^2}{4p_0}, \infty)$. Since $r_+ > r_- \geq \frac{p_1^2}{4p_0}$, we have $h_{q}(r_+) - h_{q}(r_-)>0$ and so the expression \eqref{tridiagonal1} is positive, as required. The case $\psi_\gamma(x) = e^{\gamma(x-1)}$ now follows from the previous parts by noting
\begin{equation} \label{limittransition}
T^{\psi_\gamma}_n =  \lim_{m \to \infty} T^{(1+\gamma(x-1)/m)^m}_n = \lim_{m \to \infty} T^{1+\gamma(x-1)/m}_n \cdots T^{1+\gamma(x-1)/m}_n. 
\end{equation}

For the final point, note that stochasticity of $T_n^{\psi_\gamma}$ implies
\begin{equation} \label{equivform}
\Vert T^{\psi_\gamma}_n - Id \Vert_{l^1(\mb{J}_n)} = \sup_{\lambda \in \mb{J}_n} \sum_{\mu \in \mb{J}_n} (\delta_{\lambda \mu} - T^{\psi_\gamma}_n(\lambda,\mu)) = \sup_{\lambda \in \mb{J}_n} (2 - 2T_n^{\psi_\gamma}(\lambda,\lambda))
\end{equation}
and \eqref{lowerbound} implies
$$
\ba
\inf_{\lambda \in \mb{J}_n} T^{1+\gamma(x-1)/m}_n(\lambda,\lambda) &=  \frac{h_n(r_+) - h_n(r_-)}{\sqrt{p_0^2-p_1^2}} \\
&= p_0 \left ( \sum_{k=0}^{n-1} r_+^k r_-^{(n-1)-k} \right ) - \frac{p_1^2}{4}\left ( \sum_{k=0}^{n-2} r_+^k r_-^{(n-2)-k} \right )  \overset{\gamma \to 0}{\to} 1,
\ea
$$
where we have taken $p_0 = 1-\gamma/m>\gamma/m = p_1$ in definitions that depend on these quantities. By \eqref{limittransition}, $T^{\psi_\gamma}_n(\lambda,\lambda) \geq \lim_{m \to \infty} [T^{(1+\gamma(x-1)/m)}_n(\lambda,\lambda)]^m$, so letting $\gamma \to 0$ in \eqref{equivform} completes the proof.
\epf

\subsection{Intertwined multilevel dynamics: discrete steps}

Define \emph{cotransition probabilities} from $\mb{J}_{n+1}$ down to level $\mb{J}_n$ by 
$$
T_n^{n+1}(\mu, \lambda) = \frac{\dim_{n}\lambda}{\dim_{n+1}\mu} \cdot 1_{(\lambda \prec \mu)};
$$
cf. the consistency relations \eqref{consistency}. Identity \eqref{weightidentity} implies $\dim_{n+1} \mu$ is the number of interlaced sequences of length $n+1$ ending in $\mu$, so summing $T_n^{n+1}(\mu, \lambda)$ over $\lambda \in \mb{J}_n$ evaluates to $1$. These stochastic operators ensure that the interlacing condition is preserved when we link single level dynamics to a multilevel evolution, but the following intertwining relationship is key.

\begin{proposition} \label{interlacingrelations}
Fix $\psi \in C^1([-1,1])$ with $\psi(1)=1$. For $n \geq 1$, the stochastic operators $T^\psi_n$ and $T_n^{n+1}$ satisfy the intertwining relations 
$$
\Delta_n^{n+1} : =T^{n+1}_n T^\psi_n = T^\psi_{n+1}T^{n+1}_n.
$$
\end{proposition}
\bpf
First assume $n$ is even, so that $\alpha_n = 1/2$, $\alpha_{n+1} = -1/2$, $r_n = r_{n+1} -1$, and $\phi_n(s,t) = 2 \cdot 1_{(s<t)}$. For $\mu \in \mb{J}_{n+1}$ and $\lambda \in \mb{J}_n$, we compute directly
$$
\ba
(T^{n+1}_nT^\psi_n)(\mu, \lambda) & = \sum_{z \in \mb{J}_n} T_n^{n+1}(\mu,z) T^\psi_n(z, \lambda) = \sum_{z \in \mb{J}_n}\frac{\dim_{n}z}{\dim_{n+1}\mu} \cdot 1_{(z \prec \mu)} T^\psi_n(z, \lambda)  \\
& = \frac{\dim_n \lambda}{\dim_{n+1} \mu} 2^{-r_{n+1}} \sum_{z \in \mb{J}_n} \det \left [ \phi_n(\widetilde{z}_j,\widetilde{\mu}_i) \right ]_{i,j=1}^{r_{n+1}} \det \left [  \left \langle J_{\widetilde{z}_j, 1/2}, J_{\widetilde{\lambda}_j, 1/2} \psi \right \rangle_{1/2} \right ]_{i,j=1}^{r_{n}} =\cdots,
\ea
$$
where we have used Lemma \ref{indication}. Expand the first determinant along the $r_{n+1}$-column, where the convention $z_{r_{n+1}} \equiv -1$ applies. Omitting constants, the $l$th resulting summand, $1 \leq l \leq r_{n+1}$, is 
\begin{align*}
 \sum_{z \in \mb{J}_n} \det \left [ \phi_n(\widetilde{z}_j, \widetilde{\mu}_i) \right ]_{\substack{1 \leq i \neq l \leq r_{n+1} \\ 1 \leq j \leq r_n}} & \det \left [  \left \langle J_{\widetilde{z}_j, 1/2}, J_{\widetilde{\lambda}_j, 1/2} \psi \right \rangle_{1/2} \right ]_{i,j=1}^{r_{n}} \\
  & = \det \left [ \sum_{k \geq 0} \phi_n(k,\widetilde{\mu}_i)  \left \langle J_{k, 1/2}, J_{\widetilde{\lambda}_j, 1/2} \psi \right \rangle_{1/2} \right ]_{\substack{1 \leq i \neq l \leq r_{n+1} \\ 1 \leq j \leq r_n}} \\
& = \det \left [ \left \langle  \sum_{k = 0}^{\widetilde{\mu}_i - 1} 2 J_{k, -1/2}, J_{\widetilde{\lambda}_j, 1/2} \psi \right \rangle_{1/2} \right ]_{\substack{1 \leq i \neq l \leq r_{n+1} \\ 1 \leq j \leq r_n}}  \\
& = \det \left [ \left \langle \frac{J_{\widetilde{\mu}_i, -1/2}-1}{x-1}, J_{\widetilde{\lambda}_j, 1/2} \psi \right \rangle_{1/2} \right ]_{\substack{1 \leq i \neq l \leq r_{n+1} \\ 1 \leq j \leq r_n}}  \\
& =  \det \left [ 2 \left \langle (J_{\widetilde{\mu}_i, -1/2}\psi)(1) - J_{\widetilde{\mu}_i, -1/2}\psi, J_{\widetilde{\lambda}_j, 1/2} \right \rangle_{-1/2} \right ]_{\substack{1 \leq i \neq l \leq r_{n+1} \\ 1 \leq j \leq r_n}}  \\
 & =  \det \left [ 2 \left \langle J_{\widetilde{\mu}_i, -1/2}\psi, \sum_{k = \widetilde{\lambda}_j+1}^\infty J_{k, -1/2} \right \rangle_{-1/2} \right ]_{\substack{1 \leq i \neq l \leq r_{n+1} \\ 1 \leq j \leq r_n}} 
 \\
 & = \det \left [ \left \langle J_{\widetilde{\mu}_i, -1/2}\psi, \sum_{k =0}^\infty \phi_n(\widetilde{\lambda}_j, k) J_{k, -1/2} \right \rangle_{-1/2} \right ]_{\substack{1 \leq i \neq l \leq r_{n+1} \\ 1 \leq j \leq r_n}},
\end{align*}
where we have used Cauchy-Binet in the first equality, the first part of Lemma \ref{identities} for the third, multilinearity/skew-symmetry of determinants and the assumption $(J_{\widetilde{\mu}_i, -1/2}\psi)(1) =1$ for the fourth, and the second part of Lemma \ref{identities} for the fifth. Continuing our calculations, we sum over $l$ to get
\begin{align*}
\cdots & = \frac{\dim_n \lambda}{\dim_{n+1} \mu} 2^{-r_{n+1}}   \sum_{l=1}^{r_{n+1}} (-1)^{r_{n+1}+l} \det \left [ \left \langle J_{\widetilde{\mu}_i, -1/2}\psi, \sum_{k =0}^\infty \phi_n(\widetilde{\lambda}_j, k) J_{k, -1/2} \right \rangle_{-1/2} \right ]_{\substack{1 \leq i \neq l \leq r_{n+1} \\ 1 \leq j \leq r_n}} \\
& = \frac{\dim_n \lambda}{\dim_{n+1} \mu} 2^{-r_{n+1}} \det \left [ \left \langle J_{\widetilde{\mu}_i, -1/2}\psi, \sum_{k =0}^\infty \phi_n(\widetilde{\lambda}_j, k) J_{k, -1/2} \right \rangle_{-1/2} \right ]_{i,j=1}^{r_{n+1}} \\
& = \frac{\dim_n \lambda}{\dim_{n+1} \mu} 2^{-r_{n+1}} \sum_{z \in \mb{J}_{n+1}} \det \left [ \left \langle J_{\widetilde{\mu}_i, -1/2}\psi,  J_{z_j, -1/2} \right \rangle_{-1/2} \right ]_{i,j=1}^{r_{n+1}} \cdot \det  \left [ \phi_n(\widetilde{\lambda}_j, z_i)  \right ]_{i,j=1}^{r_{n+1}} \\
& = (T^\psi_{n+1} T^{n+1}_n)(\mu,\lambda),
\end{align*}
where for the first equality we note that the entries with $j=r_{n+1}$ are all $1$, which follows from $\widetilde{\lambda}_{r_{n+1}} \equiv -1$, $(J_{\widetilde{\mu}_i, -1/2}\psi)(1) =1$, and orthogonal decomposition \eqref{orthogonaldecomposition}. The much more straightforward $n$ odd case involves the other components of Lemmas \ref{identities} and \ref{indication} and is left to the reader.

\epf
For $u,t \in  \mathbb{J}_{n,paths}$ define
\begin{equation} \label{conditionalprob}
L^{k+1, \psi}_k(u,t):=\frac{T^\psi_{k+1}(u^{k+1}, t^{k+1}) T_k^{k+1}(t^{k+1}, t^k)}{\Delta^{k+1}_k(u^{k+1},t^k)} 1_{(\Delta^{k+1}_k(u^{k+1},t^k) \neq 0 )}.
\end{equation}
In words, it is the probability of the transition ``$u^{k+1} \to t^{k+1} \to t^k$" by a jump and cotransition conditional on ``$u^{k+1} \to t^k$" occurring by such steps; note this quantity only depends on $u^{k+1},t^k,t^{k+1}$. Let $\psi$ be either $(1-p_1)+p_1x$, $1/2>p_1 \geq 0$, or  $\psi_\gamma(x):=e^{\gamma(x-1)}$. Just as the stochastic operators $T_n^\psi$ account for transitions between the probability measures $P_n^\psi$ on $\mb{J}_{n}$, the stochastic (by Proposition \ref{Markov1}) transition operator
$$
A_n^{\psi}(u, t) :=  T^\psi_1(u^1,t^1) \cdot \prod_{k=1}^{n-1} L^{k+1, \psi}_k(u,t),
$$
supplies an evolution of probability measures on $\mathbb{J}_{n,paths}$ of the form
\begin{equation} \label{seqmeasure}
P_{n,paths}^{\psi}(t): = P^{\psi}_n(t^n) \cdot \prod_{k=1}^{n-1} T^{k+1}_k(t^{k+1},t^k), \ \ \ t \in \mb{J}_{n,paths}
\end{equation}
in the usual sense $ (P_{n,paths}^{\psi_1} A^{\psi_2}_n)(t) = P_{n,paths}^{\psi_1\psi_2}$ for $\psi_1, \psi_2 \in C^1[-1,1]$, $\psi_1(1),\psi_2(1) \neq 0$. To prove this, compute, for $t \in \mb{J}_{n,paths}$,
\begin{align} \label{pathtrans}
(P_{n,paths}^{\psi_1} & A^{\psi_2}_n)(t) = \sum_{u \in \mb{J}_{n,paths}}  P_{n,paths}^{\psi_1}(u)A^{\psi_2}_n(u,t) \notag \\
& = \sum_{u^2 \prec \cdots \prec u^n} P^{\psi_1}_n(u^n) \prod_{k=2}^{n-1} T^{k+1}_k(u^{k+1},u^k) \left ( \sum_{u^1 \prec u^2} T^{2}_1(u^{2},u^1)T^{\psi_2}_1(u^1,t^1) \right ) \prod_{k=1}^{n-1} L^{k+1, \psi_2}_k(u,t) \notag \\
& = \sum_{u^3 \prec \cdots \prec u^n} P^{\psi_1}_n(u^n)  \prod_{k=3}^{n-1} T^{k+1}_k(u^{k+1},u^k) \times  \notag \\
& \phantom{ {}= \sum_{u^3 \prec \cdots \prec u^n} P^{\psi_1}_n(u^n) }  \times \left (  \sum_{u^2 \prec u^3} T^{3}_2(u^{3},u^2)T^{\psi_2}_2(u^2,t^2) \right )  \prod_{k=2}^{n-1} L^{k+1, \psi_2}_k(u,t) \cdot T^{2}_1(u^{2},u^1)= \cdots  \notag \\
 \cdots & = \left ( \sum_{u^n \in \mb{J}_n} P^{\psi_1}_n(u^n)  T^{\psi_2}_k(u^{n},t^{n} ) \right )  \prod_{k=1}^{n-1} T^{k+1}_k(t^{k+1},t^k) = P^{\psi_1\psi_2}_n(t^n)  \prod_{k=1}^{n-1} T^{k+1}_k(t^{k+1},t^k) ,   
\end{align}
where the last equality follows from the first part of Proposition \ref{Markov1}.

\subsection{Second construction of symplectic Plancherel growth} \label{continuoustime}

\bp \label{plancherelgrowth}
Let $\psi_\gamma(x) = e^{\gamma(x-1)}$ and let $\psi \in C^1[-1,1]$ satisfy $\psi(1) = 1$. Then 
$$
P_{n,paths}^{\psi} e^{\gamma Q_{n}}  =P_{n,paths}^{\psi} A^{\psi_\gamma}_n,
$$
where the  $\mathbb{J}_{n,paths} \times \mathbb{J}_{n,paths}$-matrix $Q_n$ is the infinitesimal generator of (uncoordinatized) symplectic Plancherel growth considered up to level $n$.
\ep
\bpf
Let $\mf{B}_n$ be the Banach space defined as the completion of the subspace of $l^1(\mb{J}_{n,paths})$ consisting of measures of the form \eqref{seqmeasure} corresponding to functions $\psi$ nonzero at $1$. The stochastic operators $\{ A^{\psi_\gamma}_n \}_{\gamma \geq 0}$ form a semigroup on $\mf{B}_n$ by \eqref{pathtrans}. The Feller property for $\{ A^{\psi_\gamma}_n \}_{\gamma \geq 0}$ follows from the same property of $\{ T^{\psi_\gamma}_n \}_{\gamma \geq 0}$ (the fifth part of Proposition \ref{Markov1}): the form  \eqref{seqmeasure} implies
$$
\ba
\sum_{u \in \mb{J}_{n,paths}} |P_{n,paths}^{\psi_1 \psi_2}(u) - P_{n,paths}^{\psi_1}(u)| & = \sum_{u^n \in \mb{J}_{n}} |P_{n}^{\psi_1 \psi_2}(u^n) - P_{n}^{\psi_1}(u^n)| \sum_{\substack{u \in \mb{J}_{n-1,paths} \\ u^{n-1} \prec u^n}} \prod_{k=1}^{n-1} T^{k+1}_k(u^{k+1},u^k) \\
& = \sum_{u^n \in \mb{J}_{n}} |P_{n}^{\psi_1 \psi_2}(u^n) - P_{n}^{\psi_1}(u^n)|,
\ea
$$
so we have  
$$
\Vert A^{\psi_\gamma}_n - Id \Vert_{\mf{B}_n} = \Vert T^{\psi_\gamma}_n - Id \Vert_{l^1(\mb{J}_n)} \overset{\gamma \to 0}{\to} 0
$$
Hence (see, e.g., Chapter 19, \cite{k1}), there exists an operator $\bar{Q}_n$ on $\mf{B}_n$ such that $P_{n,paths}^{\psi} e^{\gamma \bar{Q}_n} = P_{n,paths}^{\psi} A^{\psi_\gamma}_n$, so it suffices to show $\bar{Q}_n = \frac{d}{d\gamma} \big |_{\gamma=0} A^{\psi_\gamma}_n =Q_n$. Since the calculation for the second order approximation of $A^{\psi_\gamma}_n$ only differs from \cite{abjk1} for transitions involving jumps into the wall, we only present this case, the others being similar. 

Assume the system is in a state so that a wall jump is possible at an odd level $k \geq 1$, i.e., a state $t \in \mb{J}_{n,paths}$ that satisfies $t^k_{r_k} = \widetilde{t^k_{r_k}} = 0$, $t^k_{r_k-1} \geq 1$, so $\widetilde{t^k_{r_k-1}} \geq 2$, and if $k \geq 3$, $t^{k-1}_{r_{k-1}} = \widetilde{t^{k-1}_{r_{k-1}} } \geq 1$ (cf. \eqref{coordinates1}). Consider a transition to $u \in \mb{J}_{n,paths}$ that agrees with $t \in \mb{J}_{n,paths}$ except that $u^k_{r_k} = \widetilde{u^k_{r_k}}= 1$. We compute
\begin{equation} \label{component}
\frac{T^{\psi_\gamma}_k(t^{k},u^{k} ) T_{k-1}^{k}(u^{k}, u^{k-1})}{\Delta^{k}_{k-1}(u^{k},t^{k-1})} 
\end{equation}
to second order. Noting $\psi_\gamma(x)=1 + \gamma(x-1)+O(\gamma^2)$ for small $\gamma >0$, we get, similarly to \eqref{tridiagonal1},
\begin{align*}
T^{\psi_\gamma}_k(t^{k},u^{k} ) & = \det \left [ \left \langle J_{\widetilde{t}^n_i, -1/2} , J_{\widetilde{u}^n_j, -1/2} \psi_\gamma \right \rangle_{-1/2} \right ]_{i,j = 1}^{r_k} \\
& = \det 
\begin{bmatrix}
(1-\gamma) & \frac{\gamma}{2} \delta_{\widetilde{t^k_1}, \widetilde{t^k_2} + 1} & 0 & 0 & \cdots & 0  \\
\frac{\gamma}{2} \delta_{\widetilde{t^k_2}+1, \widetilde{t^k_1}} & (1-\gamma)  & \frac{\gamma}{2}\delta_{\widetilde{t^k_2}, \widetilde{t^k_3} + 1}  & 0 & \cdots & 0  \\
0 & \frac{\gamma}{2} \delta_{\widetilde{t^k_3}+1, \widetilde{t^k_2}} & (1-\gamma)  & \frac{\gamma}{2} \delta_{\widetilde{t^k_3}, \widetilde{t^k_4} + 1} & \cdots & 0  \\
\vdots & \vdots  & \ddots  & \ddots & \vdots & 0  \\
0 & \cdots  & \cdots  & \cdots &  (1-\gamma) &  0  \\
0 & \cdots  & \cdots  & \cdots & 0 & \gamma \cdot 2^{-(\beta +1/2)}  \\
\end{bmatrix} + O(\gamma^2) \\
& = \gamma \cdot 2^{-(\beta +1/2)} + O(\gamma^2).
\end{align*}
where $\beta = \pm 1/2$ is the second parameter of Jacobi polynomials (Section \ref{sectionnotations}), indicating whether we are dealing with the initial conditions of the first (-) or third (+) kinds of Chebyshev polynomials, which correspond to the orthogonal and symplectic case, respectively. Similar calculations in the denominator of \eqref{component} yield an order of $1 + O(\gamma)$. Since the dimensions occurring in the numerator and denominator can all be shown to cancel, we arrive at
$$
\ba
\frac{T^{\phi_\gamma}_k(t^{k},u^{k} ) T_{k-1}^{k}(u^{k}, u^{k-1})}{\Delta^{k}_{k-1}(u^{k},t^{k-1})} &= \frac{\gamma \cdot 2^{-(\beta +1/2)} + O(\gamma^2)}{ 1 + O(\gamma)} \\
& = \left [ \gamma \cdot 2^{-(\beta +1/2)} + O(\gamma^2)\right ] (1 + O(\gamma)) = \gamma \cdot 2^{-(\beta +1/2)} + O(\gamma^2).
\ea
$$
where the second equality uses the geometric series. This shows that the orthogonal case $\beta = -1/2$ involves rate $1$ wall jumps and the symplectic case $\beta = 1/2$ involves rate $1/2$ wall jumps, as required.
\epf

\section{Asymptotic analysis of symplectic Plancherel growth} \label{sectionasymptotics}

\subsection{Bulk limit: incomplete beta kernel}

In this section, we work in the coordinates $\mb{J}_{n} \ni \lambda \mapsto \widetilde{\lambda}$. Focusing on the kernel $K^\gamma$ at two limiting values, assume $\gamma \geq 0$ and $(s,n),(t,m) \in \Z_{\geq 0} \times \Z_{> 0}$ depend on $N$ in such a way that $\gamma \sim N\tau >0$, $s, t \sim N\nu>0$, and $r_{n}, r_m \sim N\eta$. Assume the differences $s - t$, $r_{n} - r_{m}$ are of constant order. Recall the double integral term from our correlation kernel $K^\gamma$:
$$
\ba
\frac{2^{\alpha_n+1/2}}{\pi} \frac{1}{2 \pi i} \int_{-1}^1 \oint \frac{e^{\gamma x}}{e^{\gamma u}} J_{s,\alpha_n}(x) J_{t,\alpha_m}(u) \frac{(x-1)^{r_n}}{(u-1)^{r_m}} \frac{(1-x)^{\alpha_n}(1+x)^{1/2}}{x-u} dudx,
\ea
$$
where the positively oriented $u$-contour encloses the unit circle. Deforming this $u$-contour to be a loop centered at $1$ and passing through $\cos \theta$, for some $\theta$ to be determined, will produce a residue that is considered below. With an aim toward using the identities \eqref{identity1} make the change of variables $u = (v+v^{-1})/2$ so that the $v$-contour is \emph{outside} the unit circle but connecting $e^{i\theta}$ to $e^{-i\theta}$. Similarly, with the change of variables $x = (z+z^{-1})/2$, the weight becomes
\begin{equation} \label{weight}
(1-x)^\alpha(1+x)^{1/2}dx \to m_\alpha(dz):=
\begin{cases}
 \frac{ (z^{1/2} + z^{-1/2})^2}{4iz} dz & \text{if}  \ \  \alpha = -1/2 \\
\frac{-(z-z^{-1})^2}{8iz} dz & \text{if}   \ \  \alpha = 1/2
\end{cases}.
\end{equation}
The additional factor of ``$1/2$" appears because the mapping is two to one. Taking into account the residue at $u = x$ from the first deformation, we are left with the kernel 
\begin{equation}\label{limit1}
\ba 
& K^\gamma(s,n;t,m) \\
& =\frac{2^{\alpha_n+1/2}}{2\pi^2 i} \oint_{|z|=1} \oint \frac{e^{\gamma (z+z^{-1})/2}}{e^{\gamma (v+v^{-1})/2}} J_{s,\alpha_n} \left (\frac{z+z^{-1}}{2} \right ) J_{t,\alpha_m} \left (\frac{v+v^{-1}}{2} \right) \frac{(\frac{z+z^{-1}}{2}-1)^{r_n}}{(\frac{v+v^{-1}}{2}-1)^{r_m}} \frac{(1-v^{-2})dv m_\alpha(dz)}{z+z^{-1}-(v+v^{-1})} 
\ea
\end{equation}
\begin{equation} \label{limit2}
\ba 
&+1_{(n\geq m)} \frac{2^{\alpha_n+1/2}}{\pi} \oint_{|z|=1}  J_{s,\alpha_n}  \left (\frac{z+z^{-1}}{2} \right ) J_{t,\alpha_m} \left (\frac{z+z^{-1}}{2} \right )  \left (\frac{z+z^{-1}}{2}-1 \right)^{r_n-r_m} m_\alpha(dz) \\
& +\frac{2^{\alpha_n+1/2}}{\pi} \oint_{e^{-i\theta}}^{e^{i\theta}} J_{s,\alpha_n} \left (\frac{z+z^{-1}}{2} \right ) J_{t,\alpha_m} \left (\frac{z+z^{-1}}{2} \right )  \left (\frac{z+z^{-1}}{2}-1 \right)^{r_n-r_m}  m_\alpha(dz).
\ea
\end{equation}
where the last $z$ contour connects $e^{\pm i\theta}$ clockwise along the unit circle.
\bl \label{lemma1}
The following limit holds:
$$
\ba 
\lim_{N \to \infty} \eqref{limit2}  =1_{(n \geq m)} & \frac{2^{r_m - r_n}(-1)^{\alpha_n - \alpha_m}}{2\pi i} \oint_{|z|=1} z^{(t+r_m+\alpha_m)-(s+r_n+\alpha_n)-1} (1-z)^{n-m} dz \\
& + \frac{2^{r_m - r_n}(-1)^{\alpha_n - \alpha_m}}{2\pi i} \oint_{e^{-i\theta}}^{e^{i\theta}} z^{(t+r_m+\alpha_m)-(s+r_n+\alpha_n)-1} (1-z)^{n-m} dz 
\ea
$$
where the contour connecting $e^{\pm i\theta}$ is negatively oriented and crosses $(-\infty,0)$.
\el 
\bpf
Both terms in \eqref{limit2} can be handled by the same argument, so assume $n \geq m$ and concentrate on the first of these. We will use the identities \eqref{identity1} but by shrinking or enlarging the unit circle, we may ignore terms involving $z^{\pm(t+s)}$. This leaves us with
\begin{equation} \label{hardcalc1}
\ba
\frac{2^{\alpha_n+1/2}}{\pi} J_{s,\alpha_n} & \left (\frac{z+z^{-1}}{2} \right ) J_{t,\alpha_m} \left (\frac{z+z^{-1}}{2} \right )  \left (\frac{z+z^{-1}}{2}-1 \right)^{r_n-r_m} m_\alpha(dz) \\
& =
\begin{cases}
(z^{t-s} + z^{s-t})  \left (\frac{z+z^{-1}}{2}-1 \right)^{r_n-r_m} \frac{dz}{8iz} \frac{2}{\pi} & \text{if} \ \ \alpha_n = \alpha_m= \frac{1}{2} \\
(z^{t-s-1/2} - z^{s-t+1/2})  \left (\frac{z+z^{-1}}{2}-1 \right)^{r_n-r_m} \frac{(z-z^{-1})}{z^{1/2}+z^{-1/2}} \frac{dz}{8iz}  \frac{2}{\pi} & \text{if} \ \ \alpha_n = - \alpha_m = \frac{1}{2} \\
(z^{t-s} + z^{s-t})  \left (\frac{z+z^{-1}}{2}-1 \right)^{r_n-r_m} \frac{dz}{4iz}   \frac{1}{\pi}  & \text{if} \ \ \alpha_n =  \alpha_m =- \frac{1}{2} \\
(z^{t-s+1/2} - z^{s-t-1/2})  \left (\frac{z+z^{-1}}{2}-1 \right)^{r_n-r_m} \frac{z^{1/2}+z^{-1/2}}{z-z^{-1}} \frac{dz}{4iz} \frac{1}{\pi} & \text{if} \ \ \alpha_n = - \alpha_m = - \frac{1}{2}
\end{cases} \\
& = 
\begin{cases}
z^{(t+r_m)-(s+r_n) -1} \left (1-z\right)^{2r_n-2r_m} \frac{dz}{2^{r_n-r_m} 2\pi i} & \text{if} \ \ \alpha_n = \alpha_m= \frac{1}{2} \\
- z^{(t+r_m)-(s+r_n)-3/2}  \left (1-z \right)^{2r_n-2r_m} \frac{1-z}{z^{1/2}} \frac{dz}{2^{r_n-r_m} 2\pi i}  & \text{if} \ \ \alpha_n = - \alpha_m = \frac{1}{2} \\
z^{(t+r_m)-(s+r_n) -1} \left (1-z\right)^{2r_n-2r_m} \frac{dz}{2^{r_n-r_m} 2\pi i} & \text{if} \ \ \alpha_n =  \alpha_m =- \frac{1}{2} \\
- z^{(t+r_m)-(s+r_n)-1/2}  \left (1-z \right)^{2r_n-2r_m} \frac{z^{1/2}}{1-z} \frac{dz}{2^{r_n-r_m} 2\pi i} & \text{if} \ \ \alpha_n = - \alpha_m = - \frac{1}{2}
\end{cases}
\ea
\end{equation}
where the second equality follows first from making the substitution $z \to z^{-1}$ to the integrals determined by the second summands in each case, then from using the identity $z+z^{-1} -2 = z^{-1}(z-1)^2$. Note that in the first of these steps, the transformation changes the orientation of the $z$ contour over the unit circle but an additional sign change occurs after reorienting. Letting $N \to \infty$ completes the proof.
\epf

\begin{figure}[h]
\centering
\includegraphics[width=0.45\textwidth]{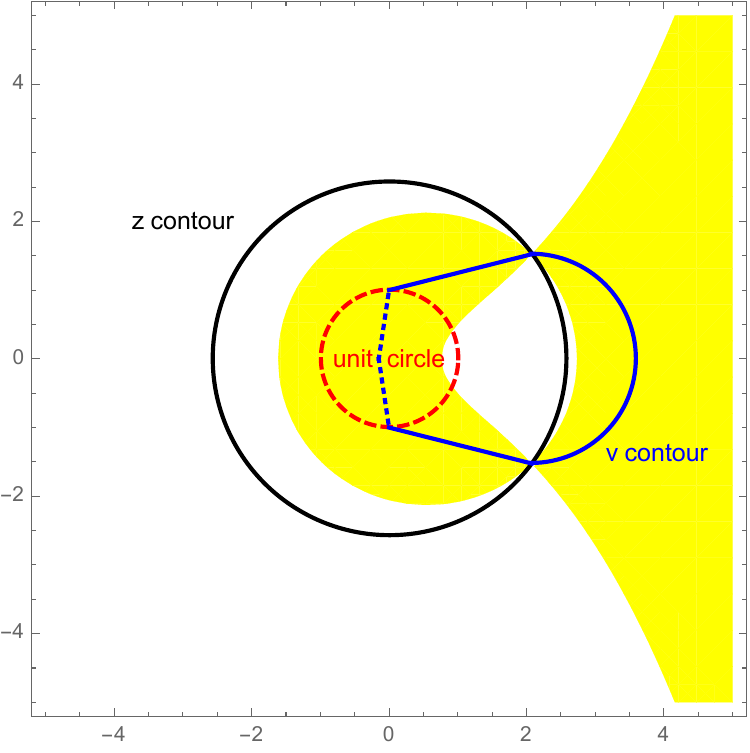} 
\caption{The yellow region signifies the region $\mf{R}((S(v) - S(z_0)))>0$ and the white elsewhere. The dotted blue line indicates the completion of the $v$-contour into a loop by the transformation $v \to v^{-1}$.} \label{liquidasymptotics}
\end{figure}

For the double integral term \eqref{limit1}, we again use the identities \eqref{identity1} to compute
\begin{equation} \label{hardcalc2}
\ba
 J_{s,\alpha_n} & \left (\frac{z+z^{-1}}{2} \right )  J_{t,\alpha_m} \left (\frac{v+v^{-1}}{2} \right) m_\alpha(dz)
\\
& = 
\begin{cases}
\frac{-dz}{8 i z} \left [ (z^sv^t) zv - (z^{-s}v^t)\frac{v}{z}+z^{-s} v^{-t}\frac{1}{zv} - z^s v^{-t} \frac{z}{v} \right ] \frac{z-z^{-1}}{v-v^{-1}} & \text{if} \ \ \alpha_n = \alpha_m= \frac{1}{2} \\
\frac{-dz}{8 i z} \left [ (z^sv^t) z \sqrt{v} - (z^{-s}v^t)\frac{\sqrt{v}}{z}-z^{-s} v^{-t}\frac{1}{z\sqrt{v}} +z^s v^{-t} \frac{z}{\sqrt{v}} \right ] \frac{z-z^{-1}}{v^{1/2}+v^{-1/2}} & \text{if} \ \ \alpha_n = - \alpha_m = \frac{1}{2} \\
\frac{dz}{4 i z} \left [ (z^sv^t) \sqrt{zv} + (z^{-s}v^t)\frac{\sqrt{v}}{\sqrt{z}}+z^{-s} v^{-t} \frac{1}{\sqrt{zv}} + z^s v^{-t} \frac{\sqrt{z}}{\sqrt{v}} \right ] \frac{z^{1/2}+z^{-1/2}}{v^{1/2}+v^{-1/2}}  & \text{if} \ \ \alpha_n =  \alpha_m =- \frac{1}{2} \\
\frac{dz}{4 i z} \left [ (z^sv^t) \sqrt{z}v + (z^{-s}v^t)\frac{v}{\sqrt{z}}-z^{-s} v^{-t}\frac{1}{\sqrt{z}v} - z^s v^{-t} \frac{\sqrt{z}}{v} \right ] \frac{z^{1/2}+z^{-1/2}}{v-v^{-1}} & \text{if} \ \ \alpha_n = - \alpha_m = - \frac{1}{2}
\end{cases} \\
& = 
\begin{cases}
\frac{-dz}{8 i z} \left [ - 2(z^{-s}v^t)\frac{v}{z}+2z^{-s} v^{-t}\frac{1}{zv} \right ] \frac{z-z^{-1}}{v-v^{-1}} & \text{if} \ \ \alpha_n = \alpha_m= \frac{1}{2} \\
\frac{-dz}{8 i z} \left [ - 2(z^{-s}v^t)\frac{\sqrt{v}}{z}-2z^{-s} v^{-t}\frac{1}{z\sqrt{v}}  \right ] \frac{z-z^{-1}}{v^{1/2}+v^{-1/2}} & \text{if} \ \ \alpha_n = - \alpha_m = \frac{1}{2} \\
\frac{dz}{4 i z} \left [  2(z^{-s}v^t)\frac{\sqrt{v}}{\sqrt{z}}+2z^{-s} v^{-t} \frac{1}{\sqrt{zv}} \right ] \frac{z^{1/2}+z^{-1/2}}{v^{1/2}+v^{-1/2}}  & \text{if} \ \ \alpha_n =  \alpha_m =- \frac{1}{2} \\
\frac{dz}{4 i z} \left [ 2 (z^{-s}v^t)\frac{v}{\sqrt{z}}-2z^{-s} v^{-t}\frac{1}{\sqrt{z}v} \frac{\sqrt{z}}{v} \right ] \frac{z^{1/2}+z^{-1/2}}{v-v^{-1}} & \text{if} \ \ \alpha_n = - \alpha_m = - \frac{1}{2}
\end{cases}
\ea
\end{equation}
where the second equality follows from the substitution $z \to z^{-1}$ applied to the integrals determined by the first and fourth summands in each case. Now we make the additional transformation $v \to v^{-1}$ applied to the integrals of the second terms in each case, which completes the $v$-contour into a positively oriented loop through the points $e^{\pm i \theta}$, as indicated by the dotted line in Figure \ref{liquidasymptotics}. Most importantly, we have a double integral over two loops with integrand
\begin{equation} \label{integrand1}
\ba
 \frac{2^{\alpha_n+1/2}}{2\pi^2 i}  J_{s,\alpha_n} & \left (\frac{z+z^{-1}}{2} \right )  J_{t,\alpha_m} \left (\frac{v+v^{-1}}{2} \right)  \frac{1-v^{-2}}{z+z^{-1}-(v+v^{-1})}  dv m_\alpha(dz)\\
 &=
 \frac{1}{2\pi i}  \frac{ z^{-s}v^t}{2 \pi i z}  \frac{1-v^{-2}}{z+z^{-1}-(v+v^{-1})} dv dz
  \cdot
\begin{cases}
\frac{v}{z} \frac{z-z^{-1}}{v-v^{-1}} & \text{if} \ \ \alpha_n = \alpha_m= \frac{1}{2} \\
 \frac{\sqrt{v}}{z} \frac{z-z^{-1}}{v^{1/2}+v^{-1/2}} & \text{if} \ \ \alpha_n = - \alpha_m = \frac{1}{2} \\
 \frac{\sqrt{v}}{\sqrt{z}} \frac{z^{1/2}+z^{-1/2}}{v^{1/2}+v^{-1/2}}  & \text{if} \ \ \alpha_n =  \alpha_m =- \frac{1}{2} \\
 \frac{v}{\sqrt{z}} \frac{z^{1/2}+z^{-1/2}}{v-v^{-1}} & \text{if} \ \ \alpha_n = - \alpha_m = - \frac{1}{2}
\end{cases}
\ea
\end{equation}
The part of the integrand depending on the hydrodynamic scaling parameter $N$ becomes 
\begin{equation} \label{integrandN}
\frac{e^{-N(S(z_0) - S(z))}}{e^{N(S(v)-S(z_0))}},
\end{equation}
where $z_0$ is to be determined and where, as in \eqref{steepestdescent}, 
$$
S(z) = S_{\tau,\nu,\eta}(z) : = \tau \frac{z+z^{-1}}{2}  - \nu \log z + \eta \log \left ( \frac{z+z^{-1}}{2} -1 \right ).
$$
To identify the correct steepest descent/ascent deformations, we first rely on Proposition 5.1 of \cite{abjk1} to determine $\theta$ and $z_0$:
\bp{(Proposition 5.1 of \cite{abjk1})} \label{analyticlines}
For $\tau/\eta \geq 0$, define
$$
q_{\pm}(\tau,\eta) = q_{\pm}\left(\frac{\tau}{\eta}\right) : = \sqrt{ - \frac{\tau^2}{2\eta^2} + \frac{5\tau}{\eta} + 1\pm \frac{\tau^2}{2\eta^2} \left ( 1 + \frac{4\eta}{\tau} \right )^{3/2} }
$$
but set $q_-$ to $0$ if $1/2 \leq \tau/\eta$. Let $z_{max}, z_{min}$ be the largest, smallest \emph{real} roots of $S'(z)$, and let $q_0$ be the complex root of $S'(z)$ in $\mb{H} \setminus \mb{D}$ . Also define 
$$
\begin{cases}
\mc{D}^{frozen} := \left \{ (\tau,\nu,\eta) \in \R_+^3 : \nu \geq \eta \cdot q_+(\tau,\eta) \right \} \\
\mc{D}^{liquid} := \left \{ (\tau,\nu,\eta) \in \R_+^3 :  \eta \cdot q_{-}(\tau,\eta) < \nu < \eta \cdot q_{+}(\tau,\eta) \right \}  \\
\mc{D}^{empty} := \left \{ (\tau,\nu,\eta) \in \R_+^3 :  \nu \leq \eta \cdot q_-(\tau,\eta) \right \}
\end{cases}.
$$
Then the quantity
$$
z_0 = z_0(\tau,\nu,\eta) :=
\begin{cases}
z_{min}, &  (\tau,\nu,\eta)  \in \mc{D}^{empty} \\
q_0, &(\tau,\nu,\eta)  \in \mc{D}^{liquid} \\
z_{max}, & (\tau,\nu,\eta)  \in \mc{D}^{frozen}
\end{cases}
$$
is critical for $S(z)$ and satisfies $z_{max} > 1$ in the first case; $|q_0| >1$ in the second; and $z_{min} < -1$ in the last.
\ep
Our work below will confirm that the definitions of the regions $\mc{D}^{liquid}$, $\mc{D}^{frozen}$, and $\mc{D}^{empty}$ given here in this proposition are consistent with the introduction's.

For $(\tau,\nu,\eta)  \in \mc{D}^{liquid}$, choose $\theta$ so that $e^{\pm i\theta}$ are points on the unit circle in the region of $z$ where the condition $\mf{R}((S(z_0) - S(z)))>0$ holds. Since $|z_0|>1$, we enlarge the $z$-unit disk to be a steepest \emph{descent} loop staying in the region $-\mf{R}((S(z_0) - S(z)))<0$ but passing through the critical points $z_0$ and $\overline{z_0}$. Similarly, make the steepest ascent $v$-loop pass through the points $z_0, \overline{z_0}$ and stay in the region $\mf{R}((S(v) - S(z_0)))>0$. Figure \ref{liquidasymptotics} exemplifies these deformations. The double integral \eqref{limit1} is now seen to vanish as $N \to \infty$, but leaves behind residues. By inspection, the residue at $v=z^{-1}$ involves the term $v^{s + t}$ with $|v| <1$ for Lebesgue almost all $v$ and hence vanishes as $N \to \infty$. The expression \eqref{integrand1} tells us the residues at $v=z$ are
\begin{equation} \label{hardresidues}
\ba
  \frac{1}{2\pi i} \int_{\overline{z_0}}^{e^{-i \theta}} z^{t-s-1} & \left (\frac{z+z^{-1}}{2}-1 \right )^{r_n-r_m}
\begin{cases}
1 & \text{if} \ \ \alpha_n = \alpha_m= \frac{1}{2} \\
 \frac{z-1}{z} & \text{if} \ \ \alpha_n = - \alpha_m = \frac{1}{2} \\
1  & \text{if} \ \ \alpha_n =  \alpha_m =- \frac{1}{2} \\
\frac{z}{z-1} & \text{if} \ \ \alpha_n = - \alpha_m = - \frac{1}{2}
\end{cases} \cdot dz  \\
& +  \frac{1}{2\pi i} \int_{e^{i\theta}}^{z_0} z^{t-s-1} \left (\frac{z+z^{-1}}{2}-1 \right )^{r_n-r_m}
\begin{cases}
1 & \text{if} \ \ \alpha_n = \alpha_m= \frac{1}{2} \\
 \frac{z-1}{z} & \text{if} \ \ \alpha_n = - \alpha_m = \frac{1}{2} \\
1  & \text{if} \ \ \alpha_n =  \alpha_m =- \frac{1}{2} \\
\frac{z}{z-1} & \text{if} \ \ \alpha_n = - \alpha_m = - \frac{1}{2}
\end{cases} \cdot dz
\ea
\end{equation}
We have implicitly used the fact that $ \frac{1-v^{-2}}{1+\frac{z^{-1}-v^{-1}}{z-v}} \to 1$ as $v \to z$.  Hence, we have proven
\bl \label{lemma2}
If $(\tau,\nu,\eta)  \in \mc{D}^{liquid}$, then
$$
\ba
\lim_{N \to \infty} \eqref{limit1}  =  & \frac{2^{r_m - r_n}(-1)^{\alpha_n-\alpha_m}}{2\pi i} \int_{\overline{z_0}}^{e^{-i \theta}} z^{(t+r_m+\alpha_m)-(s+r_n + \alpha_n) -1} \left (1-z \right )^{n-m} \\
&+ \frac{2^{r_m - r_n}(-1)^{\alpha_n-\alpha_m}}{2\pi i} \int_{e^{i\theta}}^{z_0} z^{(t+r_m + \alpha_m)-(s+r_n+\alpha_n)-1} \left (1-z \right )^{n-m }.
\ea 
$$
\el
Noting that the conjugating factor ``$2^{r_m - r_n}(-1)^{\alpha_n-\alpha_m}$" vanishes in the determinant, Lemmas \ref{lemma1} and \ref{lemma2} complete the proof of Theorem \ref{bulklimit} once we confirm the regions are consistently defined with the introduction.

\begin{figure}[h]
\centering
\includegraphics[width=0.45\textwidth]{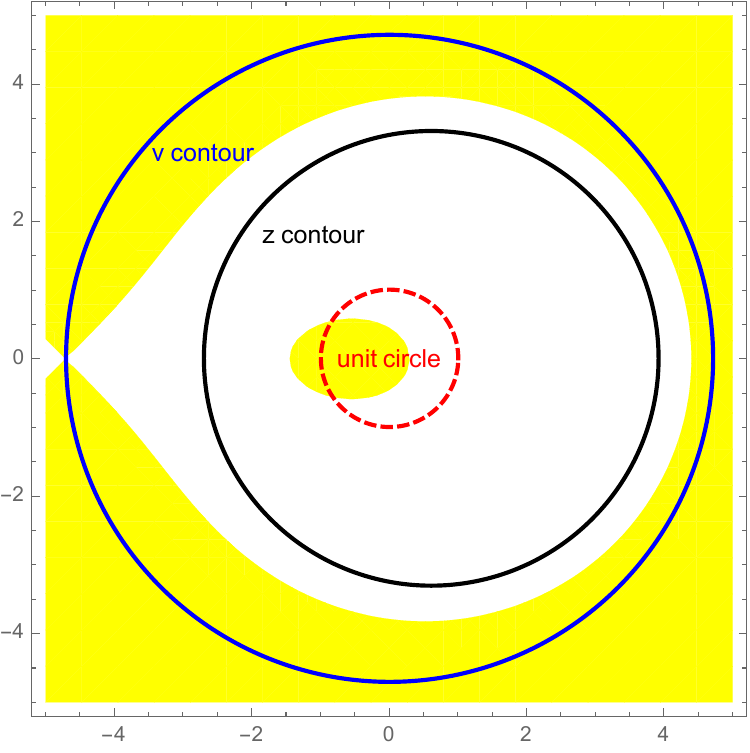} 
\includegraphics[width=0.45\textwidth]{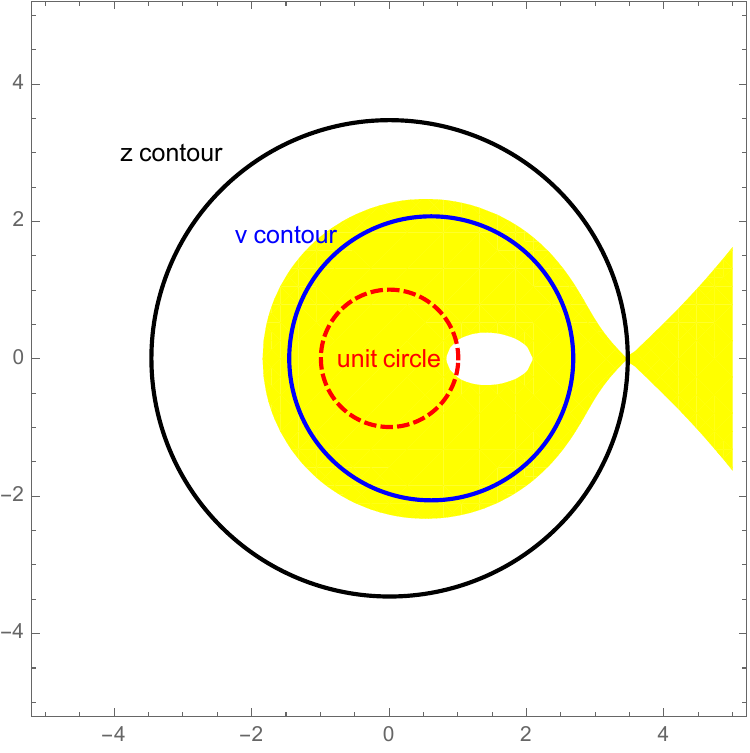}  
\caption{The plot exemplifies steepest descent deformations for the frozen (left) and empty (right) regions. The yellow region again signifies the region $\mf{R}((S(v) - S(z_0)))>0$ and the white elsewhere.} \label{emptyasymptotics}
\end{figure}
For the case $(\tau,\nu,\eta)  \in \mc{D}^{frozen}$, do not deform the initial $u$ contour and only make the substitutions $x=(z+z^{-1})/2$, $u=(v+v^{-1})/2$, as before. Steepest descent deformations as in the plot of Figure \ref{emptyasymptotics} do not produce residues at $v=z^{-1}$ since $z_0(\tau,\nu,\eta)=z_{min}<-1$ and also no residues at $v=z$ by considering the indicated regions of the plot. Hence, the double contour vanishes as $N \to \infty$ without leaving behind any terms and we are left with the triangular matrix governed by the first term in Lemma \ref{lemma1}:
$$
1_{(n \geq m)}  \frac{2^{r_m - r_n}(-1)^{\alpha_n - \alpha_m}}{2\pi i} \oint_{|z|=1} z^{(t+r_m+\alpha_m)-(s+r_n+\alpha_n)-1} (1-z)^{n-m} dz.
$$
Since the diagonal entries of this triangular matrix are $1$, the determinant of $K^\gamma$ converges to $1$. This confirms the definition of $\mc{D}^{frozen}$.

Finally, for $(\tau,\nu,\eta)  \in \mc{D}^{empty}$, again refrain from deforming the initial $u$ contour, but make the substitutions $x=(z+z^{-1})/2$, $u=(v+v^{-1})/2$, as before. We may perform the same calculation leading to \eqref{hardcalc2}, but now the transformation $v \to v^{-1}$ to the second summand at the end of \eqref{hardcalc2} leads to a loop inside the unit circle (currently the $z$ contour). This second summand disappears without leaving anything behind since we may shrink its $v$-loop arbitrarily to avoid residues at $v=z^{-1}$ while making $z$-deformations as in Figure \ref{emptyasymptotics}. But these deformations in the first summand of \eqref{hardcalc2} leave behind a full loop of residues at $v=z$ similar to \eqref{hardresidues} and Lemma \ref{lemma2}:
$$
-  \frac{2^{r_m - r_n}(-1)^{\alpha_n - \alpha_m}}{2\pi i} \oint_{|z|=1} z^{(t+r_m+\alpha_m)-(s+r_n+\alpha_n)-1} (1-z)^{n-m} dz,
$$  
Since we did not perform any deformations to the initial $u$ contour, we are asymptotically left with the sum of this last residue and the first term in Lemma \ref{lemma1}:
$$
-1_{(n < m)}  \frac{2^{r_m - r_n}(-1)^{\alpha_n - \alpha_m}}{2\pi i} \oint_{|z|=1} z^{(t+r_m+\alpha_m)-(s+r_n+\alpha_n)-1} (1-z)^{n-m} dz
$$
Since the diagonal entries of this triangular matrix are $0$, the determinant of $K^\gamma$ converges to $0$. This confirms the definition of $\mc{D}^{empty}$.

\subsection{Edge limits}

\subsubsection{Finite distance from the wall: Jacobi kernel}

Assume again that $\gamma \sim N \cdot \tau >0$ and $r_{n}, r_m \sim N \cdot \eta>0 $, but that $s, t$ are fixed and finite. Assume only the difference $n - m$ is of constant order. Let $A(z) := \tau z + \eta \log (1-z)$ and note $1-\eta/\tau$ is a zero of $A'(z)$. Write our kernel as 
$$
\ba
 \frac{2^{\alpha_n+1/2}}{\pi} \frac{1}{2 \pi i}& \int_{-1}^1 \oint \frac{e^{- N(A(1-\eta/\tau) - A(x))}}{e^{ N(A(u) - A(1-\eta/\tau))}}  J_{s,\alpha_n}(x) J_{t,\alpha_m}(u) \frac{(x-1)^{r_n}}{(u-1)^{r_m}} \frac{(1-x)^{\alpha_n}(1+x)^{1/2}}{x-u} dudx\\
& +1_{(n\geq m)} \left \langle J_{s,\alpha_n}, (x-1)^{r_n-r_m} J_{t,\alpha_m} \right \rangle_{\alpha_n}. 
\ea
$$
\begin{figure}[h]
\centering
\includegraphics[width=0.45\textwidth]{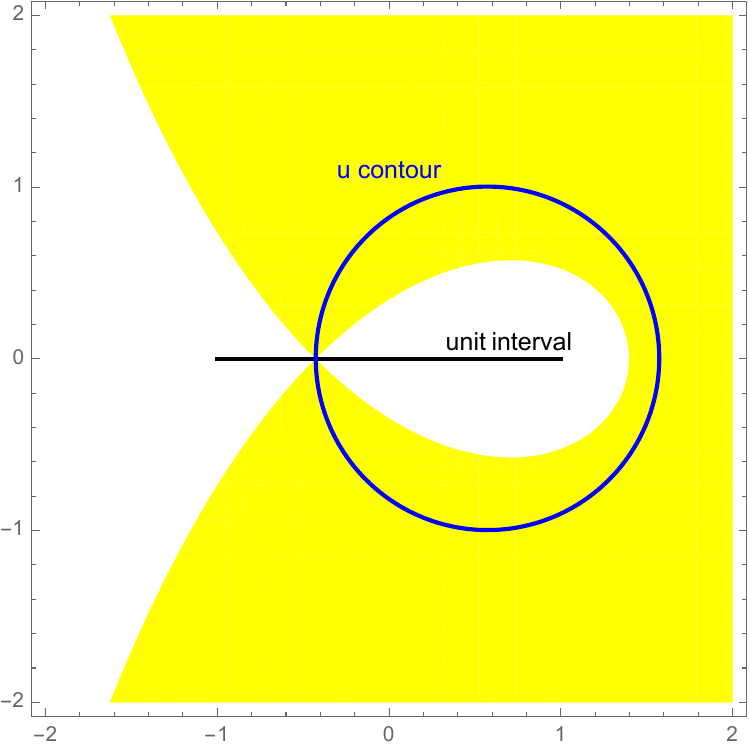}
\caption{The plot exemplifies steepest descent deformations when targeting a finite distance from the wall. The yellow region signifies where $\mf{R}(A(u)-A(1-\eta/\tau))>0$ and the white elsewhere.} \label{jacobiasymptotics}
\end{figure}
Now deform the $u$-contour, as in Figure \ref{jacobiasymptotics}, to be a steepest ascent loop remaining in the region $ \mc{R}(A(u) - A(1-\eta/\tau)) > 0$ and passing through the critical point $1-\eta/\tau$. If $1-\eta/\tau< -1$, the unit interval already lies in the region $ \mc{R}(A(1-\eta/\tau) - A(x)) >0$, so the double integral term tends to zero without picking up residues. As for the frozen region above, the kernel converges to a triangular matrix whose diagonal entries are $1$'s. But if $1-\eta/\tau> -1$, then we acquire a residue at $u = x$ given by
$$
 -  \frac{2^{\alpha_n+1/2}(-1)^{r_n-r_m}}{\pi}  \int_{-1}^{1-\eta/\tau} J_{s,\alpha_n}(x) J_{t,\alpha_m}(x) (1-x)^{r_n-r_m+\alpha_n} (1+x)^{1/2} dx
$$
Noting that $(-1)^{r_n-r_m}$ is a conjugating factor completes the proof of Theorem \ref{jacobilimit}.

\subsubsection{Corner of phase transition: Pearcey kernel}

Assume that $\gamma \sim N/2$, and $s \sim N^{1/4} \nu_1>0$, $t \sim N^{1/4} \nu_2>0$, and lastly that $r_{n} - N\sim \sqrt{N} \eta_1>0$, $r_{m} - N\sim \sqrt{N} \eta_2>0$

\bl \label{asymp}
Let $x = N^{-1/2}x'-1 \in [-1,1]$. Then for large $N$, 
$$
N^{-1/4} (-1)^s  J_{s, \alpha}(x) \sim \frac{\sin[\nu_1 \sqrt{2 x'} ] }{2^\alpha \sqrt{x'}}.
$$
\el
\bpf
Write $x = \cos \theta$ and note the power series
$$
\cos^{-1}(z-1) = \pi -  \int_0^z \frac{1}{\sqrt{2x}\sqrt{1-x/2}} dx = \pi - \sqrt{2} \sqrt{z} + \frac{\sqrt{2}}{12} z^{3/2} + \cdots,
$$
If $\alpha = -1/2$, then by \eqref{identity1}, $J_{s,-1/2}(x) = \frac{ \cos (s+1/2)\theta}{\cos(\theta/2)}$ so that
$$
\ba
(-1)^s J_{s,-1/2}(x) & = (-1)^s \frac{\cos[(s+1/2)\pi - s N^{-1/4} \sqrt{2x'} + o(1)]}{\cos [ \pi/2 - N^{-1/4} \sqrt{x'/2} + o(1)]} \\
& = (-1)^s \frac{\cos[(s+1/2)\pi+ o(1)] \cos[ \nu_1 \sqrt{2x'}] + \sin[(s+1/2)\pi+ o(1)] \sin[ \nu_1 \sqrt{2x'}]  }{\cos[\pi/2+ o(1)] \cos[ N^{-1/4} \sqrt{x'/2}] + \sin[\pi/2+ o(1)] \sin[ N^{-1/4}  \sqrt{x'/2}]} \\
& \sim \frac{\sin[ \nu_1 \sqrt{2x'}]  }{\sin[ N^{-1/4}  \sqrt{x'/2}]},
\ea
$$
and we may similarly use $J_{s,1/2}(x)= \frac{\sin(s+1)\theta}{\sin\theta}$ if $\alpha = 1/2$. Using $\lim_{y \to 0} \frac{y}{\sin y} = 1$ completes the proof.
\epf

Recall for any $\omega = (\alpha, \beta, \gamma)$, the complementary point process $ (\widetilde{\mc{X}}^\omega  )^c$ of $\widetilde{\mc{X}}^\omega$ is also determinantal with $k$th correlation function
$$
\widetilde{\rho}^{\omega,k}_\Delta(z_1,\ldots,z_k) := \mathbb{P} \left ( (\widetilde{\mc{X}}^\omega  )^c \supset \{z_1,\ldots, z_k \} \right ) = \det[ K_\Delta^\omega(z_i, z_j) ]_{i,j=1}^k,
$$
where the kernel $K^\omega_\Delta$ is readily computable (using orthogonality for the Chebyshev polynomials) as
\begin{equation} \label{compkernel}
\ba
K^\omega_\Delta((s,n),(t,m)) :=  - \frac{2^{\alpha_n+1/2}}{\pi} \frac{1}{2 \pi i} & \int_{-1}^1 \oint \frac{E^\omega(x)}{E^\omega(u)} J_{s,\alpha_n}(x) J_{t,\alpha_m}(u) \frac{(x-1)^{r_n}}{(u-1)^{r_m}} \frac{(1-x)^{\alpha_n}(1+x)^{1/2}}{x-u} dudx\\
& -1_{(n>m)} \left \langle J_{s,\alpha_n}, (x-1)^{r_n-r_m} J_{t,\alpha_m} \right \rangle_{\alpha_n}.
\ea
\end{equation}
\begin{figure}[h]
\centering
\includegraphics[width=0.45\textwidth]{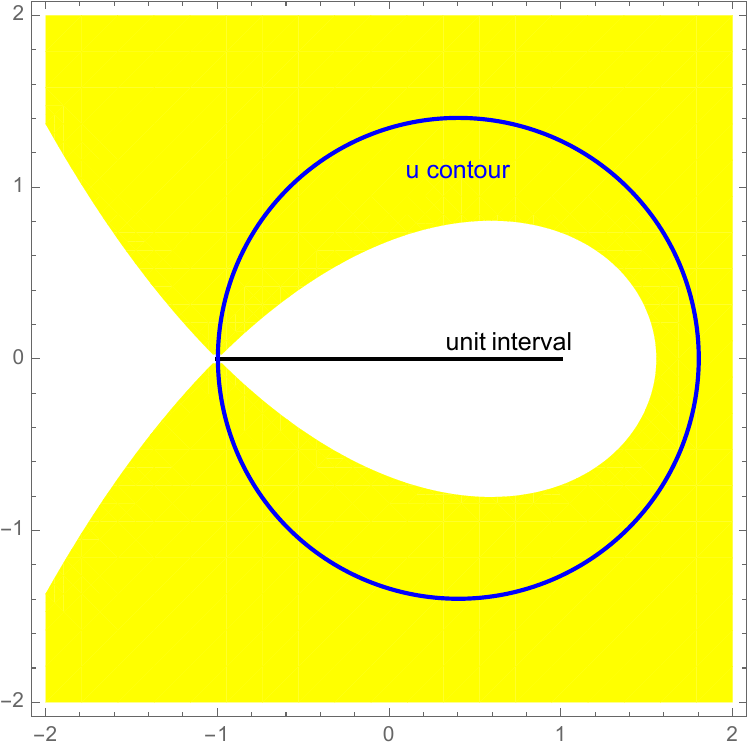}
\caption{Steepest descent deformations when targeting the corner of the phase transition. The yellow region signifies the region $\mf{R}(A(u)-A(-1))>0$ and the white elsewhere.} \label{pearceyasymptotics}
\end{figure}
Deforming the $u$-contour in the double integral term of \eqref{compkernel} (with $\omega = (0,0,\gamma)$) as in Figure \ref{pearceyasymptotics}, we endeavor to find the contribution at $-1$. Making the substitutions $x' = N^{1/2}(x+1)$ and $u' = N^{1/2}(u+1)$, we have, for large $N$,
$$
\ba
\frac{(1-x)^{\alpha_n}(1+x)^{1/2}}{x-u} dudx & = N^{1/2} \frac{(2-x'N^{-1/2})^{\alpha_n}\sqrt{x'} N^{-1/4}}{x'-u'} \cdot N^{-1} du'dx' \\
& \sim N^{-3/4} \frac{2^{\alpha_n}\sqrt{x'}}{x'-u'} \cdot du'dx'.
\ea
$$
Putting this together with Lemma \ref{asymp}, we get
$$
2^{\alpha_m} (-1)^{s} J_{s,\alpha_n}(x) (-1)^{t} J_{t,\alpha_m}(u)\frac{(1-x)^{\alpha_n}(1+x)^{1/2}}{x-u} dudx \sim N^{-1/4} \sin[\nu_1 \sqrt{2 x'} ] \sin[\nu_2 \sqrt{2 u'} ] \frac{du'dx'}{\sqrt{u'}(x'-u')}.
$$
The remainder of the integrand satisfies
$$
\ba
(-1)^{r_m-r_n} 2^{r_m-r_n} \frac{e^{\gamma x}}{e^{\gamma u}}  \frac{(x-1)^{r_n}}{(u-1)^{r_m}} & = \frac{ e^{-N(A(-1)-A(x))+ \eta_1 \sqrt{N}(\log(1-x) - \log 2)}}{e^{N(A(u) - A(-1))+ \eta_2 \sqrt{N}(\log(1-u) - \log 2)}} \sim \frac{ e^{-(x')^2/8 - \eta_1 x'/2}}{ e^{-(u')^2/8 - \eta_2 u'/2}},
\ea 
$$
where $A(z): = z/2+ \log(1-z)$ and the last approximation follows from, respectively, second and first order Taylor expansions about $-1$. Hence, putting everything together, we have
$$
\ba
  N^{1/4} 2^{\alpha_m - \alpha_n} (-2)^{r_m-r_n}  & (-1)^{t-s} \frac{2^{\alpha_n+1/2}}{\pi} \frac{1}{2 \pi i}  \int_{-1}^1 \oint \frac{e^{\gamma x}}{e^{\gamma u}} J_{s,\alpha_n}(x) J_{t,\alpha_m}(u) \frac{(x-1)^{r_n}}{(u-1)^{r_m}} \frac{(1-x)^{\alpha_n}(1+x)^{1/2}}{u-x} dudx \\
& \to  \frac{\sqrt{2}}{ \pi}  \frac{1}{2 \pi i}  \oint_{-i\infty}^{i\infty} \int_0^\infty e^{ \frac{(u')^2-(x')^2}{8} + \frac{\eta_2 u' - \eta_1 x'}{2}}\sin[\nu_1 \sqrt{2 x'} ] \sin[\nu_2 \sqrt{2 u'} ] \frac{dx'du'}{\sqrt{u'}(u'-x')}
\ea
$$
Turning to the other term of $K^\gamma_\Delta$, note the following two asymptotic relations:
$$
(-2)^{r_m-r_n} (x-1)^{r_n-r_m} = \left (1-\frac{x'}{2 \sqrt{N}} \right )^{(\eta_1 - \eta_2) \sqrt{N}} \to e^{\frac{-(\eta_1-\eta_2)x'}{2}}.
$$
$$
2^{\alpha_m} (-1)^{s} J_{s,\alpha_n}(x) (-1)^{t} J_{t,\alpha_m}(x) (1-x)^{\alpha_n}(1+x)^{1/2} dx \sim N^{-1/4} \sin[\nu_1 \sqrt{2 x'} ] \sin[\nu_2 \sqrt{2 x'} ] \frac{dx'}{\sqrt{x'}}.
$$
Then the single integral term satisfies
$$
\ba
N^{1/4} & 2^{\alpha_m - \alpha_n}(-2)^{r_m-r_n}  (-1)^{t-s} \left ( -1_{(n>m)} \left \langle J_{s,\alpha_n}, (x-1)^{r_n-r_m} J_{t,\alpha_m} \right \rangle_{\alpha_n} \right ) \\
&  \to 1_{(\eta_1 > \eta_2)} \left ( - \frac{\sqrt{2}}{\pi} \int_0^\infty  e^{\frac{-(\eta_1-\eta_2)x'}{2}}   \sin[\nu_1 \sqrt{2 x'} ] \sin[\nu_2 \sqrt{2 x'} ] \frac{dx'}{\sqrt{x'}} \right ) \\
& = 1_{(\eta_1 > \eta_2)} \left ( \frac{1}{\sqrt{2}\pi} \int_0^\infty  e^{\frac{-(\eta_1-\eta_2)x'}{2}}  \left [ \cos[ (\nu_1+\nu_2) \sqrt{2 x'} ]  - \cos[ (\nu_1 - \nu_2) \sqrt{2 x'} ] \right ] \frac{dx'}{\sqrt{x'}} \right )\\
& = 1_{(\eta_1 > \eta_2)} \frac{1}{\sqrt{\pi (\eta_1 - \eta_2)}} \left ( \exp \left [ - \frac{  (\nu_1 + \nu_2)^2}{\eta_1 - \eta_2} \right ] - \exp \left [ - \frac{ (\nu_1 - \nu_2)^2}{\eta_1 - \eta_2} \right ] \right ).
\ea
$$
where the last equality follows from standard Gaussian computations. Noting that the conjugating factors will disappear in the determinant completes the proof of Theorem \ref{pearceylimit}.

\appendix 

\section{Algebraic argument for determinantal correlations} \label{appendix}

For completeness, this section covers the linear algebraic details on which Proposition \ref{checkconditions} relies; we follow a combination of Theorem 4.2 of \cite{abpf2} and Lemma 3.4 of \cite{bfs1}, which in turn rely on the Eynard-Mehta theorem in the manner of \cite{br1}. 

Fix $n \geq 1$. For $1 \leq m \leq n$ and $N \in \Z_{\geq 0}$, consider point sets $\mathfrak{X}^{(m)} : = \{ x^m_k \}_{k=1}^{r_m}$ of variables $x^m_k$ ranging in $\{0,\ldots, N\}$, where we suppress dependence on $N$. Define $\mf{Y} := \cup_{m=1}^n \mathfrak{X}^{(m)}$ and recall our work in Section \ref{checkconditions} furnishes a measure on point configurations in $\mf{Y}$ up to level $n$ by
\begin{equation} \label{detmeasure2}
\ba
 \text{const} \cdot \prod_{m=1}^{n-1} \det \left[ \phi_m(x^m_i, x^{m+1}_j ) \right]_{i,j=1}^{r_{m+1}} \cdot \det \left [ \Psi^n_{r_n-j}(x^{n}_i) \right]_{i,j=1}^{r_n}.
 \ea
\end{equation}
Note \eqref{detmeasure2} uses our notational convention $x^m_{r_{m+1}} \equiv -1$ for $m$ even (recall the proof of Lemma \ref{indication}), so add the \emph{virtual variables} $\mc{V} := \{ v_i \equiv -1 \}_{i=1}^{r_n}$ as a notational device to form the larger space $\mathfrak{X} := \mc{V} \cup \mf{Y}$. To separate these variables out, let $W_m$, $1 \leq m \leq n-1$, be the $\mathfrak{X}^{(m)} \times \mathfrak{X}^{(m+1)}$ matrix with entries $W_m(x^m_i, x^{m+1}_j):= \phi_m(x^m_i, x^{m+1}_j)$, and similarly, for $1 \leq m \leq n$, define a $\mc{V} \times \mf{X}^{(m)}$ matrix $E_{m-1}$ to have entries $E_{m-1}(v_i, x^m_j) = \phi_{m-1}(v_i,x^m_j) 1_{(i=r_m)}1_{(\text{m odd})}$. Define a $\mf{X}^{(k)} \times \mc{V}$ matrix $\mbf{\Psi}^{k}$ to have entries $\mbf{\Psi}^{k}(x^{k}_i,v_j):=\Psi^k_{r_k-j}(x^{k}_i) = \langle J_{x^{k}_i,\alpha_k}, (x-1)^{r_k-i} R^{E^\omega}_{i-r_k} \rangle_{\alpha_k}$, as in Proposition \ref{comprule}. 

Now the measure \eqref{detmeasure2} is proportional to a minor of the $\mf{X} \times \mf{X}$ matrix 
 $$
 L:=
 \begin{bmatrix}
 0 & E_0 & 0 & E_2 & \cdots & E_{n-1}  \\
0 & 0 & -W_1 & 0 &  \ddots & 0 \\
0 & 0 & 0 & -W_2 & \ddots & 0 \\
\vdots & \ddots & \ddots &  \ddots &\ddots & \vdots \\
\vdots & \ddots & \ddots &  \ddots &\ddots & -W_{n-1}\\
\mbf{\Psi}^{n} & 0 & 0 & \cdots & 0 & 0
\end{bmatrix}
$$
The matrix $L$ supplies a point process (called an \emph{L-ensemble}) on the larger space $\mf{X}$ with distribution $\mb{P}_{\mf{X}}(X) := \frac{\det L_X}{\det[I_{\mf{X}} + L]}$, $X \in 2^{\mf{X}}$, where $L_X$ is the $X \times X$ submatrix and $I_{\mf{X}}:=I_{\mf{X} \times \mf{X}}$. Its correlation function is determinantal with kernel $K_{\mf{X}} = I_{\mf{X}}-(I_{\mf{X}}+ L)^{-1}$ (see  \cite{dj1}). Write $I_{\mf{Y}}+L= :\begin{bmatrix} 0 & A \\ C & D \end{bmatrix}$, where as indicated, $A$ is the $\mc{V} \times \mf{Y}$ block, $C$ is the $\mf{Y} \times \mc{V}$ block, and $D$ the $\mf{Y} \times \mf{Y}$ block. One can see by inspection that the $\mf{Y} \times \mf{Y}$ matrix $D^{-1}$ has $(i,j)$-block $W_{[i,j)}: = (W_1 \cdots W_{j-1}) 1_{(i<j)} +  I_{\mf{X}^{(j)} \times \mf{X}^{(j)}} 1_{(i=j)}$. Notice
$$
AD^{-1} = \left ( E_{0}W_{[1,1)}, \sum_{k=1}^{2} E_{k-1}W_{[k,2)}, \sum_{k=1}^{3} E_{k-1}W_{[k,3)}, \ldots, \sum_{k=1}^{n} E_{k-1}W_{[k,n)} \right )
$$
Define also an $r_n \times r_n$ matrix $M_n:=AD^{-1}C = \sum_{k=1}^{n} E_{k-1}W_{[k,n)}\mbf{\Psi}^{n}$.
\bp
The point process of particle configurations up to level $n$ governed by \eqref{detmeasure2} is the same as the $\mf{Y}$ conditional $L$-ensemble defined by the prescription $\mb{P}_{\mf{Y}}(Y) := \frac{\det L_{Y \cup \mf{Y}^c}}{\det[I_{\mf{Y}} + L_{\mf{Y}}]}$, $Y \in 2^{\mf{X}}$, which is determinantal with correlation kernel 
\begin{equation} \label{kernelEM}
K = I_{\mf{Y}} -(I_{\mf{Y}}+ L)^{-1}|_{\mf{Y}\times \mf{Y}} = I_{\mf{Y}} - D^{-1} + D^{-1}CM_n^{-1}AD^{-1}.
\end{equation}
More explicitly, its $\mf{X}^{(i)} \times \mf{X}^{(j)}$-block is given by
\begin{equation} \label{kernelblock}
K_{(i,j)}:= K_{\mf{X}^{(i)} \times \mf{X}^{(j)}} = -W_{[i,j)} 1_{(i<j)} + \left [ W_{[i,n)} \mbf{\Psi}^{n} \right ]  \left [ M^{-1} \left ( \sum_{k=1}^{j} E_{k-1}W_{[k,j)} \right ) \right ].
\end{equation}
\ep
\bpf
For the the basic facts on conditional ensembles and in particular the first equality of \eqref{kernelEM}, see \cite{dj1}. The second equality is the content of the Eynard--Mehta theorem, and is proved by computing the inverse of the $2 \times 2$ block expression of $I_{\mf{Y}}+L=\begin{bmatrix} 0 & A \\ C & D \end{bmatrix}$ in order to compute $(I_{\mf{Y}}+ L)^{-1}|_{\mf{Y}\times \mf{Y}}$.
\epf
Now let $N \to \infty$ so that matrix multiplication induces convolution over $\Z_{\geq 0}$,  which is defined by $(f*g)(x,y): = \sum_{z \geq 0} f(x,z) g(z,y)$ for bivariate functions $f,g$ and by $(f*g)(x):= \sum_{z \geq 0} f(x,z) g(z)$ if $g$ is univariate. We need to compute \eqref{kernelblock}. Let $\phi^{[i,n)}:=W_{[i,n)}$ and note $\phi^{[i,i)}:=I_{\mf{X}^{(i)}}$. Then for each $1 \leq j \leq n$,
\begin{equation} \label{rowspace}
 \left ( \sum_{k=1}^{j} E_{k-1}W_{[k,j)}  \right )(v_s,x^j_t) = \sum_{k=1}^{j} 1_{(s=r_k, \  k \ odd)} (\phi_{k-1}*\phi^{[k,j)})(v_k, x^j_t)  = (\phi_{q_s-1}*\phi^{[q_s,j)})(v_s,x^j_t),
\end{equation}
where $q_k:=2k-1$. Take a $\mc{V} \times \mf{X}^{(j)}$ matrix ${\bf \Phi}^{j}$ with entries denoted ${\bf \Phi}^{j}(v_s,x^j_t)=:\Phi^{j}_{r_j-s}(x^j_t) 1_{(s \leq r_j)}$, whose nonzero rows determine a basis for the row space  $\{ (\phi_{q_s-1}*\phi^{[q_s,j)})(v_s,x) | \ x \in \mf{X}^{(j)} \}_{s=1}^j$ of \eqref{rowspace}, and which is biorthogonal to the ${\bf \Psi}^{i}$ in the sense 
$$
\sum_{k \geq 0}  \Psi^{m}_{r_m - s}(k)\Phi^{m}_{r_m-t}(k) = \delta_{st}, \ \ \ 1 \leq m \leq n.
$$
Define an $r_j \times r_j$, $1 \leq j \leq n$, matrix $B^j$ to have entries
$$
B^j_{st}  :=(\phi_{q_s-1}*\phi^{[q_s,j)}*\Psi^{j}_{r_j-t})(v_s), \ \ \ 1 \leq s,t \leq r_j.
$$
Add $0$'s to extend $B^j$ to an $r_n \times r_n$ matrix $B^j_0$, which acts as a change of basis in the sense 
$$
B^j_0 {\bf \Phi}^{j} = [\phi_{0}*\phi^{[1,j)}, \phi_{2}*\phi^{[3,j)}, \ldots , \phi_{q_{r_j}-1}*\phi^{[q_{r_j},j)}, 0,0 ,\ldots ]'
$$ 
(note we identify $v_i \leftrightarrow i$). Note further that for $1 \leq t \leq r_{j}$,
$$
B^{j}_{r_{j},t} =(\phi_{q_{r_{j}}-1}*\phi^{[q_{r_{j}},j)}*\Psi^{j}_{r_j-t})(v_{r_j})=(\phi_{q_{r_{j}}-1}*\Psi^{q_{r_j}}_{r_j-t})(v_{r_j}) = \sum_{x \geq 0} \Psi^{q_{r_j}}_{r_j-t}(x) = \delta_{r_{j},t}
$$ 
where the second equality uses $W_{[i,n)} \mbf{\Psi}^{n} = \mbf{\Psi}^{i}$ (by Proposition \ref{comprule}) and the last uses orthogonal decomposition and $J_{k,-1/2}(1)=1$. This implies in particular that if $B^{j-1}$ is upper triangular for $j \geq 2$, then so is $B^j$. Hence, $B^n$ is upper triangular and since $B^n = B^n_0 = M_n$ by definition, we have
$$
 M_n^{-1} \left ( \sum_{k=1}^{j} E_{k-1}W_{[k,j)}  \right ) =
 \begin{bmatrix}
 B^j & * \\
 0 & * 
 \end{bmatrix}^{-1} B^j_0 {\bf \Phi}^{j}  
 = 
  \begin{bmatrix}
 (B^j)^{-1} & *' \\
 0 & *' 
 \end{bmatrix} 
 \begin{bmatrix}
 B^j & 0 \\
 0 & 0
 \end{bmatrix}
 {\bf \Phi}^{j} = {\bf \Phi}^{j}.
$$
We have thus computed the $(i,j)$-block \eqref{kernelblock} of the kernel to be
$$
-W_{[i,j)} 1_{(i<j)} + \left [ W_{[i,n)} \mbf{\Psi}^{n} \right ]  \left [ M^{-1} \left ( \sum_{k=1}^{j} E_{k-1}W_{[k,j)} \right ) \right ] = -\phi^{[i,j)}1_{(i<j)} + {\bf \Psi}^{i} {\bf \Phi}^{j}.
$$
(See \cite{bfs1} to find a slightly more complicated expression for the kernel when the upper triangularity of $M_n$ fails.) The entrywise statement reads as follows.
\bp \label{algebraicargument}
For any $(s,i), (t,j) \in \Z_{\geq 0} \times \Z_{> 0}$, the point process $\widetilde{\mc{X}}^\omega$ determined by \eqref{detmeasure2} has determinantal correlations with kernel
$$
K((s,i), (t, j)) = - \phi^{[i,j)}(s,t)1_{(i<j)} + \sum_{k=1}^{r_j} \Psi^{i}_{r_{i} - k}(s) \Phi^{j}_{r_{j} - k}(t).
$$
\ep

\bibliographystyle{plain}
\bibliography{sympgrowth}

\end{document}